\documentclass[11pt]{article}
\usepackage[centering,includeheadfoot,margin=2 cm]{geometry}
\usepackage{graphics,enumitem,epsfig,textcomp}
\usepackage{amsfonts,amsmath,amssymb,euscript,color}
\usepackage{subfigure}

\begin{document}

\def\ALERT#1{{\large\bf $\clubsuit$#1$\clubsuit$}}

\numberwithin{equation}{section}
\newtheorem{defin}{Definition}
\newtheorem{theorem}{Theorem}
\newtheorem{notice}{Notice}
\newtheorem{lemma}{Lemma}
\newtheorem{cor}{Corollary}
\newtheorem{example}{Example}
\newtheorem{remark}{Remark}
\newtheorem{conj}{Conjecture}
\def\begproof{\noindent{\bf Proof: }}
\def\endproof{\par\rightline{\vrule height5pt width5pt depth0pt}\medskip}
\def\div{\nabla\cdot}
\def\rot{\nabla\times}
\def\sign{{\rm sign}}
\def\arsinh{{\rm arsinh}}
\def\arcosh{{\rm arcosh}}
\def\diag{{\rm diag}}
\def\const{{\rm const}}
\def\d{\,\mathrm{d}}
\def\eps{\varepsilon}
\def\theta{\vartheta}
\def\N{\mathbb{N}}
\def\R{\mathbb{R}}
\def\C{\hbox{\rlap{\kern.24em\raise.1ex\hbox
      {\vrule height1.3ex width.9pt}}C}}
\def\P{\hbox{\rlap{I}\kern.16em P}}
\def\Q{\hbox{\rlap{\kern.24em\raise.1ex\hbox
      {\vrule height1.3ex width.9pt}}Q}}
\def\M{\hbox{\rlap{I}\kern.16em\rlap{I}M}}
\def\Z{\hbox{\rlap{Z}\kern.20em Z}}
\def\({\begin{eqnarray}}
\def\){\end{eqnarray}}
\def\[{\begin{eqnarray*}}
\def\]{\end{eqnarray*}}
\def\part#1#2{\frac{\partial #1}{\partial #2}}
\def\partk#1#2#3{\frac{\partial^#3 #1}{\partial #2^#3}} 
\def\mat#1{{D #1\over Dt}}
\def\dx{\nabla_x}
\def\dv{\nabla_v}
\def\grad{\nabla}
\def\Norm#1{\left\| #1 \right\|}
\def\pmb#1{\setbox0=\hbox{$#1$}
  \kern-.025em\copy0\kern-\wd0
  \kern-.05em\copy0\kern-\wd0
  \kern-.025em\raise.0433em\box0 }
\def\bar{\overline}
\def\lbar{\underline}
\def\fref#1{(\ref{#1})}
\def\half{\frac{1}{2}}
\def\oo#1{\frac{1}{#1}}

\def\tot#1#2{\frac{\d #1}{\d #2}} 
\def\laplace{\Delta}
\def\d{\,\mathrm{d}}
\def\N{\mathbb{N}}
\def\R{\mathbb{R}}
\def\T{\mathbb{T}}
\def\supp{\mbox{supp }}
\def\eps{\varepsilon}

\def\W{W}
\def\w{w}
\def\M{\mathcal{M}}
\def\O{\mathcal{O}}
\def\P{\mathcal{P}}

\def\comment#1{\textcolor{blue}{\bf [Comment: #1]}}


\centerline{{\huge Individual based and mean-field modelling of direct aggregation}}


\vskip 5mm

\centerline{
{\large Martin Burger}\footnote{Institut f\"ur Numerische und Angewandte Mathematik, Westf\"alische Wilhelms-Universit\"at M\"unster, Einsteinstr. 62, 48149 M\"unster, Germany; e-mail: {\it martin.burger@wwu.de}}
\qquad
{\large Jan Ha{\v s}kovec}\footnote{
King Abdullah University of Science and Technology, Thuwal, KSA;
e-mail: {\it jan.haskovec@kaust.edu.sa}}
\quad
{\large Marie-Therese Wolfram}\footnote{Faculty of Mathematics, Universit\"at Wien,
Nordbergstrasse 15, A-1090 Wien, Austria and DAMTP, University of Cambridge, Wilberforce Road, Cambridge CB3 0WA, Great Britain;
e-mail: {\it marie-therese.wolfram@univie.ac.at}}
}
\vskip 6mm


\noindent{\bf Abstract.}
We introduce two models of biological aggregation,
based on randomly moving particles with individual stochasticity
depending on the perceived average population density in their neighbourhood.
In the first-order model the location of each individual
is subject to a density-dependent random walk, while in the second-order model
the density-dependent random walk acts on the velocity variable,
together with a density-dependent damping term.
The main novelty of our models is that we do not assume any explicit
aggregative force acting on the individuals; instead, aggregation
is obtained exclusively by reducing the individual stochasticity in response to higher perceived density.
We formally derive the corresponding mean-field limits, leading to nonlocal degenerate diffusions.
Then, we carry out the mathematical analysis of the first-order model,
in particular, we prove the existence of weak solutions and
show that it allows for measure-valued steady states.
We also perform linear stability analysis and identify
conditions for pattern formation. Moreover, we discuss the role
of the nonlocality for well-posedness of the first-order model.
Finally, we present results of numerical simulations
for both the first- and second-order model on the individual-based
and continuum levels of description.
\vskip 5mm


\newpage
\section{Introduction}\label{sec:Introduction}
Animal aggregation is the process of finding a higher density of
animals at some place compared to the overall mean density. Its formation may be
triggered by some environmental heterogeneity that is attractive to animals (the aggregate
forms around the environmental template), by physical currents that trap the organisms
through turbulent phenomena, or by social interaction between animals~\cite{Grunbaum-Okubo-94, Mimura-Yamaguti-82}.
Aggregation may serve diverse functions such as reproduction, formation of local microclimates,
anti-predator behaviour (see for instance~\cite{Turchin-Kareiva-89} for a study
of reducing the risk of predation to an individual by aggregation
in \emph{Aphis varians}), collective foraging and much more (see, e.g., \cite{Krause-Ruxton-02}
for a relatively recent survey).
Aggregation also plays an important role as an evolutionary step towards social
organization and collective behaviours~\cite{Krebs-Davies-84}.
These aspects explain the continuing interest in understanding not only the function of animal aggregation,
but also the underlying mechanisms.
The development of quantitative mathematical models is an essential part in this quest.
While such models first of all help to link individual behavioural
mechanisms to spatio-temporal group patterns, they also often
aim at explaining the observed dynamic efficiency of animal groups to adapt to environmental variation,
and, moreover, provide a valuable tool to study in more detail the system dynamics
and their robustness to variations in individual
behavioural parameters or external parameters, see, e.g., \cite{JJGBT}.

Going back to the pioneering work of Skellam~\cite{Skellam51},
continuum spatiotemporal population dynamics are traditionally modelled
by reaction-convection-diffusion PDEs and systems thereof.
In these models, diffusion typically describes the avoidance of crowded areas
by the individuals, and as such acts as an anti-aggregative force,
working against the typically aggregative effect of convection,
see the surveys~\cite{Okubo80} and~\cite{Murray}.
Our work goes into the opposite direction and shows that biological aggregations
can be a consequence of solely random, diffusive motion of individuals,
who respond to the local population density observed in their neighborhood
by increasing or decreasing the amplitude of their random motion.
This kind of behaviour was observed in insects, for instance
the pre-social German cockroach Blattella germanica.
These animals are known to be attracted to dark, warm and humid places~\cite{Rust}.
However, the works by Jeanson et al.~\cite{Jeanson03,Jeanson05} have shown that
cockroach larvae also aggregate in the absence of any environmental template or heterogeneity.
In this case aggregation is the result of social interactions
and happens as a self-organized process. The individual
based model developed and parameterized in~\cite{Jeanson05} identified a
simple decisive mechanism that can be summarised in the following way: cockroaches
do not rest for a long time in places with few conspecifics, and once moving they stop
preferentially in places of high cockroach density.
A mathematically better tractable version of this decisive mechanism
is the model of individuals performing random walks with density dependent coefficients.
In the continuum limit, this leads to degenerate diffusion models,
which, however, depending on their parameters and data, might have ill-posed regimes.
In particular, Turchin~\cite{Turchin89} derived a 1D model for a population with density
$u(x,t)$ of the form
\(   \label{Turchin}
    \part{u}{t} = \part{}{x} \left[ \phi'(u) \part{u}{x} \right] \,,
\)
where $\phi(u) = (\mu/2) u - k_0 u^2 + (2k_0/2\omega) u^3$
with the positive coefficients $\mu$, $k_0$ and $\omega$.
Turchin used the above equation to describe the aggregative movement
of \emph{Aphis varians}, a herbivore of fireweeds (\emph{Epilobium angustifolium}).
He also discussed the possible ill-posedness of some initial and
boundary value problems associated with \eqref{Turchin}.
Depending on the actual profile of~$\phi$, Turchin also classified
two different types of aggregation.
Essentially the same model has been derived independently in~\cite{Ang-Schm}
as a model of cell motility with volume filling and cell-to-cell adhesion.
The authors showed that the diffusivity can become negative if the cell adhesion
coefficient is sufficiently large and related this to the presence of spatial oscillations
and development of plateaus in their numerical solutions of the discrete
model. Moreover, they used a combination of stability analysis of the discrete equations
and steady-state analysis of the limiting PDE to gain better understanding
of the qualitative predictions of the model.
Another work studying an equation of the type~\eqref{Maini} is~\cite{Padron98},
where existence and uniqueness of solutions of the initial value problem
with homogeneous Neumann boundary conditions in a bounded domain of $\R^n$ was shown
and some aspects of the aggregating behaviour were studied analytically.

In~\cite{Maini10}, following~\cite{Maini06}, the authors provide an alternative derivation of~\eqref{Turchin}
for low population density, using a biased random walk approach.
Then they study the existence of traveling wave solutions
for a purely negative or zero diffusion equation with a logistic rate of growth $g(u)$,
\(   \label{Maini}
   \part{u}{t} = \part{}{x} \left[ \phi'(u) \part{u}{x} \right] + g(u)\,,
\)
discuss the well- and ill-posedness of certain boundary conditions
associated with some purely negative diffusion equations with logistic-like kinetic part
and provide some numerical examples.

One possibility to overcome the difficulty of the possible ill-posedness of~\eqref{Turchin} and~\eqref{Maini}
is to introduce a nonlocality, i.e., to substitute the term $\laplace\phi(u)$ by
$\laplace J$ with
\[
    J(x,t) = \int_\Omega K(x,y) \phi(u(y,t)) \d y
\]
with a nonnegative kernel $K(x,y)$.
A particular choice made in~\cite{Grindrod88} is to define $K(x,y)$
as the Green function of $(I-\lambda\laplace)$ for a constant $\lambda > 0$.
Equation~\eqref{Maini} becomes then
\[
   \part{u}{t} = \laplace [I-\lambda\laplace]^{-1}\phi(u) + g(u) \,,
\]
which is equivalent to
\[
   \part{u}{t} = \laplace \left(\phi(u) - \lambda g(u) + \lambda\part{u}{t} \right) + g(u) \,.
\]
This equation was studied in~\cite{Padron03} and can be considered as a model of aggregating population
with a migration rate determined by $\phi$ and total birth and mortality rates characterized by $g$.
In~\cite{Padron03} it was shown that the aggregating mechanism induced by $\phi(u)$
allows for survival of a species in danger of extinction and performed numerical simulations
suggesting that, for a particular version of $\phi(u)$, the solutions stabilize asymptotically
in time to a not necessarily homogeneous stationary solution.
Other two works going in this direction are~\cite{Bertozzi},
which we discus later (Section~\ref{sec:Stability}),
and the model of home range formation in wolves due to scent marking~\cite{Lewis},
\[
    \part{u}{t} &=& \laplace\left(\frac{D_0 u}{1+p/\alpha}\right) \,,\\
    \part{p}{t} &=& u(\gamma + m(p)) - \mu p \,.
\]
Here $u(x,t)$ is the population density of wolves, $p(x,t)$
the density of their scent marks and $D_0$, $\gamma$, $\alpha$ and $\mu$ positive parameters.
The increasing function $m(p)$ describes enhanced scent mark rates in the
presence of existing scent marks. A nonlocal version of the model with
$m$ depending on the averaged version of $p$ is also considered.
The authors of~\cite{Lewis} show that
the model produces distinct home ranges; in this case the pattern formation
results from the positive feedback interaction
between the decreased diffusivity of wolves in the presence
of high scent mark densities and the increased production of new scent marks
in these locations.

In our paper we introduce and study two models where formation of aggregates
results from random fluctuations in the population density
and is supported merely by reducing the amplitude of the individual random motion
in response to higher perceived density.
This leads to nonlocal individual-based and PDE models,
where the nonlocality stems from calculating
the perceived density as a weighted average over a finite sampling radius.
This is usual in modelling biological interactions,
see for instance~\cite{Kawasaki, Grunbaum-Okubo-94, Bertozzi, Mogilner-Edelstein-Keshet} and many more.
Our models were inspired by the above mentioned observations
of German cockroach~\cite{Jeanson03, Jeanson05},
but does not aim to be a realistic description of their social behaviour.
Instead, we consider our work as a proof of concept,
where the main characteristic of our models 
is that we do not assume any deterministic interaction between the individuals
that would actively push them to aggregate.
This approach was followed for instance in stochastic run-and-tumble
models of chemotaxis~\cite{Schnitzer}. Closely related to our work is
the series of papers~\cite{Hernandez-Lopez, Lopez1, Lopez2}.
However, only the first-order model under specific conditions was studied there.
The new aspects contributed by our work are the generality of the models (first- and second-order),
more rigorous derivation of the mean-field limits based on the generalized BBGKY-hiearchy,
and rigorous well-posedness and stability analysis of the first-order model.

Our paper is structured as follows: In Section~\ref{sec:Models}
we introduce two stochastic individual based models,
where every individual is able to sense the average population density in its neighborhood,
and respond in terms of increased or decreased stochasticity of its movement.
In particular, the diffusivity is reduced in response to higher perceived density.
We present a first-order model, where the location of each individual
is subject to a density-dependent random walk, and a second-order model,
where the density-dependent random walk acts on the velocity variable,
together with a density-dependent damping term.
The advantage of the second-order model is that it is possible
to introduce a cone of vision, which depends on the direction
of each individual's movement.
In Section~\ref{sec:MeanField} we formally derive the corresponding mean-field limits, leading
to a nonlocal, nonlinear diffusion equation for the first-order model,
and a nonlinear Fokker-Planck kinetic equation for the second-order model.
Moreover, we show that the diffusive limit of the kinetic equation
leads to the first-order diffusion equation derived before.
Then, in Section~\ref{sec:Analysis}, we show the existence of weak
and measure-valued solutions of the first order
mean-field model and perform linear stability analysis of the uniform steady states,
finding regimes that correspond to the sought-for pattern formation (direct aggregation).
Finally, in Section~\ref{sec:Numerics} we present the results of numerical simulations
of our models.
We performed Langrangian simulations of the first- and second-order agent-based models
in a periodic 2D domain,
Eulerian simulations of the first-order mean-field model
in 1D and 2D periodic domains, and, finally, of the second-order mean-field model
in a spatially 1D periodic domain with 1D velocity.

\section{The stochastic individual based models}\label{sec:Models}
We consider the set of $N\in\N$ individuals
with time-dependent positions $x_i(t)\in\R^d$, $i=1, \dots, N$,
with $d\geq 1$.
Every individual is able to sense the average density of other individuals in its neighborhood,
given by
\(   \label{rho_i}
    \vartheta_i(t) = \frac{1}{N} \sum_{j\neq i} \W(x_i-x_j) \,,
\)
where $\W(x) = w(|x|)$ with $\w: \R^+ \to \R^+$ is a bounded, nonnegative and nonincreasing weight,
integrable on $\R^d$.
A generic example of $\w$ is the characteristic function
of the interval $[0,R]$, corresponding to the sampling radius $R > 0$.
The average density $\vartheta_i$ is then simply the fraction
of individuals located within the distance $R$ from the $i$-th individual:
\[
    \vartheta_i(t) = \frac{1}{N}\#\{j; |x_i-x_j| \leq R\} \,.
\]
For the passage to the mean-field limit $N\to\infty$ and the analysis
of the resulting equations, we will have to pose certain smoothnes assumptions
on $\W$, which actually exclude the choice of $\w$ to be a characteristic function.
However, from the modelling point of view this is not a concern.

Let us now introduce our two models:
\begin{itemize}
\item
In the first-order model the individual locations
are subject to average density-dependent random walk,
\(
   \d x_i = G(\vartheta_i) \d B_i^t \,,\qquad i = 1, \dots, N \,, \label{model1}
\)
where $B_i^t$ are independent $d$-dimensional Brownian motions
and $G:\R^+ \to\R^+$ is a bounded, nonnegative and nonincreasing function.
For instance, in the numerical simulations of Section~\ref{sec:Numerics}
we use $G(s) = \exp(-s/3)$. However, we as well allow for degeneracy,
where $G(s) = 0$ for all $s \geq s_0 > 0$.
\item
In the second-order model the individuals are described
not only by their locations $x_i(t)\in\R^d$, but also by the velocities $v_i(t)\in\R^d$, $i=1,\dots,N$.
The advantage of this description, in contrast to the first-order model,
is that every individual has a well defined direction of movement.
Therefore, the weight $\W$ in the calculation of the average densities $\vartheta_i$
can also depend on the relative angle with respect to the individual's direction of movement,
\(   \label{rho_i2}
    \vartheta_i(t) = \frac{1}{N} \sum_{j\neq i} \W(x_i-x_j,v_i) \,,
\)
with $\W(x,v) = w\left(|x|,\frac{x\cdot v}{|x| |v|}\right)$.
This is important from the modeling point of view,
since we can define not only the sampling radius $R>0$,
but also a restricted cone of vision with angle $\alpha\in(0,\pi]$;
we set then $w(s,z) = \chi_{[0,R]}(s)\chi_{[\cos\alpha,1]}(z)$.

The velocity in our model is subject to a density-dependent random walk
and a density-dependent damping term:
\(
    \d x_i &=& v_i \d t \,, \label{model2-1} \\
    \d v_i &=& -H(\vartheta_i) v_i \d t + G(\vartheta_i) \d B_i^t \label{model2-2} \,.
\)
The function $G$ is as before, while $H:\R^+ \to \R^+$ is a nonnegative and nondecreasing function.
The damping term $-H(\vartheta_i) v_i$ is introduced in order to slow down the agents' motion
when they approach a crowded place. Obviously, this mechanism is not only necessary
to allow for formation of aggregates, but also very natural from the modeling point of view.
\end{itemize}

\section{The mean-field limits}\label{sec:MeanField}

\subsection{The first-order model}\label{subsec:FirstOrder}
We start with the derivation of the mean-field limit
of the first-order model~\fref{model1}.
Unfortunately, the standard framework of BBGKY hierarchies cannot be applied,
since due to the structure of the model it is impossible to obtain a hierarchy
where the evolution of the $k$-th marginal is expressed in terms
of a finite number of higher-order marginals.
Therefore, we have to apply the recently developed
technique of introduction of auxiliary variables~\cite{Burger}
and their elimination after the mean-field limit passage.

In particular, under the additional assumption $\W\in C^1(\R^d)$,
we extend the model by introducing the average densities $\vartheta_i$ given by~\fref{rho_i}
as a new set of independent variables, governed by the system of stochastic differential equations
\(   \label{drho}
   \d \vartheta_i = \frac1N \sum_{j\neq i} \grad\W(x_i-x_j)\left(G(\vartheta_i)\d B^t_i - G(\vartheta_j)\d B^t_j\right) \,.
\)
Let us point out that the random walks $B^t_i$, $B^t_j$ in~\fref{drho} are correlated with
those in~\fref{model1}.
Using the It\^ o formula, we turn to the
equivalent formulation of the stochastic system~\fref{model1}, \fref{drho}
in terms of the corresponding Fokker-Planck equation
\[
   \part{f^N}{t} &=& \frac12 \sum_{i=1}^N \laplace_{x_i} \left( G(\vartheta_i)^2 f^N \right)
            + \sum_{i=1}^N \part{}{\vartheta_i}\grad_{x_i}\cdot
      \left(\frac1N \sum_{j\neq i} \grad\W(x_i-x_j) G(\vartheta_i)^2 f^N \right)    \\
     &-& \frac1N \sum_{i=1}^N \sum_{j\neq i} \sum_{k\neq i} \part{}{\vartheta_j}\grad_{x_i}\cdot
      \left(\grad\W(x_k - x_i) G(\vartheta_i)^2 f^N \right) \\
    &+& \frac{1}{2N^2} \sum_{i=1}^N \sum_{j\neq i} \sum_{k\neq i} \partk{}{\vartheta_i}{2}
     \left( \grad\W(x_i-x_k)\cdot\grad\W(x_i-x_l) G(\vartheta_i)^2 f^N \right) \\
    &+& \frac{1}{2N^2} \sum_{i=1}^N \sum_{j\neq i} \partk{}{\vartheta_i}{2}
     \left( |\grad\W(x_i-x_j)|^2 G(\vartheta_j)^2 f^N \right) \\
    &+& \frac{1}{2N^2} \sum_{i=1}^N \sum_{j\neq i} \sum_{k\neq i,j} \frac{\partial^2}{\partial\vartheta_i\partial\vartheta_j}
     \left( \grad\W(x_i-x_k)\cdot\grad\W(x_j-x_k) G(\vartheta_k)^2 f^N \right) \\
    &-& \frac{1}{N^2} \sum_{i=1}^N \sum_{j\neq i} \sum_{k\neq i,j} \frac{\partial^2}{\partial\vartheta_i\partial\vartheta_j}
     \left( \grad\W(x_i-x_k)\cdot\grad\W(x_i-x_j) G(\vartheta_i)^2 f^N \right) \,,
\]
where $f^N = f^N(t, x_1, \vartheta_1, \dots, x_N, \vartheta_N)$ is the $N$-particle distribution function.
Defining the $k$-particle marginals $f^N_k$ by
\[
    f^N_k(t,x_1,\vartheta_1,\dots,x_k,\vartheta_k) = \int_{\R^{d(N-k)}}\int_{(0,\infty)^{N-k}}
        f^N(t,x_1,\vartheta_1,\dots,x_N,\vartheta_N) \d \vartheta_{k+1}\dots\vartheta_N \d x_{k+1}\dots\d x_N \,,
\]
we obtain the so-called BBGKY-hierarchy for the system $(f^N_k)_{k=1}^N$.
In particular, the equation for $f_1 = f_1(t,x,\vartheta)$ reads
\(
   \part{f^N_1}{t} &=& \frac12 \laplace_{x} \left( G(\vartheta)^2 f^N_1 \right)
      + \part{}{\vartheta} \grad_{x}\cdot \left(
         G(\vartheta)^2 \int_{\R^d}\int_0^\infty \grad\W(x-y) f^N_2(x,\vartheta,y,\sigma) \d\sigma\d y \right)
          \nonumber\\
      &+& \frac12 \partk{}{\vartheta}{2} \left( G(\vartheta)^2 \int_{\R^d}\int_{\R^d}\int_0^\infty\int_0^\infty
         \grad\W(x-y)\grad\W(x-z) f^N_3(x,\vartheta,y,\sigma,z,\tau) \d\sigma\d\tau\d y\d z \right)
           \nonumber\\
      &+& \frac{1}{2N} \partk{}{\vartheta}{2} \left( \int_{\R^d}\int_0^\infty
        G(\vartheta)^2 |\grad\W(x-y)|^2 f^N_2(x,\vartheta,y,\sigma) \d\sigma\d y \right) \,,\label{BBGKY}
\)
where we, for the sake of legibility, dropped the indices at $x_1$ and $\vartheta_1$.
Now, we pass formally to the limit $N\to\infty$,
assuming that $\lim_{N\to\infty} f^N_k = f_k$ for all $k\geq 1$.
Moreover, we admit the usual \emph{molecular chaos} assumption
about vanishing particle correlations as $N\to\infty$,
\[
    f_k(t, x_1, \vartheta_1,\dots,x_k,\vartheta_k) = \prod_{i=1}^k f_1(t,x_i,\vartheta_i) \qquad
      \mbox{for all } k \geq 2 \,.
\]
Then, one obtains from~\fref{BBGKY} the one-particle equation
\(
   \part{f_1}{t} = \frac12\laplace_x \left(G(\vartheta)^2 f_1\right)
          + \part{}{\vartheta}\grad\cdot \left(G(\vartheta)^2 (\grad\W\circledast f_1) f_1\right)
          + \frac12 \partk{}{\vartheta}{2} \left(G(\vartheta)^2 |\grad\W\circledast f_1|^2 f_1 \right) \label{1part} \,,
\)
where the operator $\circledast$ is defined as
\[
    \grad\W\circledast f_1(t,x) = \int_{\R^d}\int_0^\infty \grad\W(x-y) f_1(t,y,\vartheta) \d\vartheta\d y \,.
\]
Finally, we reduce~\fref{1part} to obtain the standard mean-field
description of the system~\fref{model1} by removing the auxiliary variable $\vartheta$.
Indeed, a relatively lengthy, but straight-forward formal calculation
shows that~\fref{1part} posesses weak solutions of the form
\[
    f_1(t,x,\vartheta) = \varrho(t,x) \delta(\vartheta - \W\ast \varrho(t,x)) \,,\qquad
    \mbox{with }
    \W\ast \varrho(t,x) = \int_{\R^d} \W(x-y) \varrho(t,y) \d y
\]
and $\varrho=\varrho(t,x)$ satisfies the nonlinear diffusion equation
\(
   \part{\varrho}{t} = \frac{1}{2} \laplace_x\left(G(\W\ast \varrho)^2 \varrho\right) \,. \label{mmodel1}
\)

\subsection{The second-order model}\label{subsec:SecondOrder}
With the same procedure as before we derive the formal mean-field limit
of the model~\fref{model2-1}--\fref{model2-2}.
We omit the details here and immediately give the resulting
kinetic Fokker-Planck equation for the particle distribution function
$f = f(t,x,v)$,
\(
   \part{f}{t} + v\cdot\grad_x f = \grad_v\cdot\left( H(\W\circledast f)v f + \frac12\grad_v
       \left(G(\W\circledast f)^2 f\right) \right) \,,\label{mmodel2}
\)
where, with a slight abuse of notation, the convolution operator $\circledast$ is defined as
\[
    \W\circledast f(t,x,v) = \int_{\R^d}\int_{\R^d} W(x-y,v) f(t,y,w) \d w\d y \,.
\]

Let us make the following observation:
If the weight $\W$ does not depend on $v$,
such that $\W\circledast f$ does not depend on $v$ as well
and $\W\circledast f(t,x) = \W\ast\varrho(t,x)$,
we can write the (non-closed) system for 
the mass, momentum and energy densities
\(   \label{moments}
    \varrho (t,x) = \int_{\R^d} f(t,x,v) \d v \,,\qquad
    \varrho u(t,x) = \int_{\R^d} f(t,x,v) v \d v \,,\qquad
    \varrho E(t,x) = \frac12 \int_{\R^d} f(t,x,v) |v|^2 \d v \,,
\)
associated with the solution $f$ of~\fref{mmodel2} as
\(
    \part{\varrho}{t} + \grad_x\cdot(\varrho u) &=& 0 \,, \label{moments1} \\
    \part{(\varrho u)}{t} + \grad_x\cdot\left( \int_{\R^d} f(t,x,v) v\otimes v \d v \right) &=&
        - H(\W\ast\varrho) \varrho u \,, \label{moments2}  \\
    \part{(\varrho E)}{t} + \grad_x\cdot\left( \frac12 \int_{\R^d} f(t,x,v) |v|^2 v \d v \right) &=&
        - 2 H(\W\ast\varrho) \varrho E + \frac{d}{2} G(\W\ast\varrho)^2 \varrho \,. \label{moments3} 
\)
We observe that only mass is conserved, while
momentum and energy are not (neither locally nor globally).
Indeed, the momentum is dissipated due to the ``friction'' term, whose strength
depends non-locally on $\varrho$ due to $H(\W\ast\varrho)$.
The energy is, on one hand, dissipated due to the same term,
on the other hand is created due to the diffusive term in~\fref{mmodel2}, at the rate
$\frac{d}{2E} G(\W\ast\varrho)^2$.
In equilibrium, we have $H(\W\ast\varrho) \varrho u = 0$, which means that we either have
empty regions ($\varrho = 0$) or regions with positive density, but zero velocity.
In these populated regions the equilibrium energy is given by
\[
    E = E[\varrho] = \frac{d}{4} \frac{G(\W\ast\varrho)^2}{H(\W\ast\varrho)} \,.
\]

\subsection{Diffusive limit of the second-order model~\fref{mmodel2}}\label{sec:DiffusiveLimit}
We show that the first-order model~\fref{mmodel1} is obtained from~\fref{mmodel2}
in the formal diffusive limit, under the assumption that the weight $W$ does not depend
on $v$, which we make for all of this section.
Let us recall that then $W\circledast f$ does not depend on $v$ as well
and $W\circledast f(t,x) = W\ast\varrho(t,x)$, with $\varrho$ defined by~\fref{moments}.

We start by observing that the equilibria of the collision operator of~\fref{mmodel2}
are given by the local Maxwellians
\(   \label{M}
   \M[\varrho](v) = \left(\frac{H(\W\ast\varrho)}{2\pi G(\W\ast\varrho)^2}\right)^{d/2} \varrho
         \exp\left(-\frac{H(\W\ast\varrho)}{G(\W\ast\varrho)^2} |v|^2 \right) \,.
\)
Therefore, at equilibrium the mean velocity $u$ vanishes if $\varrho\neq 0$,
and the ``statistical temperature'' $T$, defined by
\[
    \frac{d}{2}\varrho T = \frac12 \int_{\R^d} \M[\varrho](v) |v|^2 \d v \,,
\]
depends non-locally on $\varrho$ and is given by
\(   \label{T}
    T = T[\W\ast\varrho] = \frac{G(\W\ast\varrho)^2}{2H(\W\ast\varrho)}   \,.
\)
Consequently, in crowded regions (where $\W\ast\varrho$ is high) the temperature
is low (``freezing of the aggregates''), while the unhabitated regions are ``hot''.

We introduce the diffusive scaling with the small parameter $\varepsilon>0$ to~\fref{mmodel2},
\[
    \eps^2 \part{f}{t} + \eps v\cdot\grad_x f = \grad_v\cdot\left( H(\W\ast\varrho)v f + \frac12\grad_v
       \left(G(\W\ast\varrho)^2 f\right) \right) \,,
\]
and consider the Hilbert expansion in terms of $\eps$, $f = \M[\varrho] + \eps g$,
where $\M[\varrho]$ is given by~\fref{M}.
Moreover, we perform the Taylor expansion of $H(\W\ast\varrho)$, which we formally write as
\[
    H(\W\ast\varrho) = H_0 + \eps H_1 + \O(\eps^2) \,,
\]
with
\[
    H_0 = H(\W\ast\varrho_\M) \,,\qquad
    H_1 = H'(\W\ast\varrho_\M) \W\ast\varrho_g \,,
\]
and similarly for $G(\W\ast\varrho) = G_0 + \eps G_1 + \O(\eps^2)$. Here, $\varrho_\M$ and $\varrho_g$ denote the velocity averages of the lowest order terms in the expansion, i.e.
\[
    \varrho_\M = \int_{\R^d} \M[\varrho] \d v, \quad 
     \varrho_g = \int_{\R^d} g \d v.
\]
Then, collecting terms of order $\eps$, we obtain
\[
    v\cdot\grad_x\M = \grad_v\cdot\left( H_0 v g + H_1 v\M + \frac12\grad_v\left(G_0 g + G_1\M\right)\right) \,,
\]
which yields, after multiplication by $v$ and integration,
\(    \label{f1}
    \frac{d}{2}\grad_x\left(\varrho_\M T[\W\ast\varrho_\M]\right) = \grad_x\cdot\int_{\R^d} \M v\otimes v\d v = 
      - H_0 \int_{\R^d} g v\d v \,,
\)
with $T[W\ast\varrho_\M]$ given by~\fref{T}.
Collecting terms of order $\eps^2$ and integrating with respect to $v$, we obtain
\[
    \part{\varrho_\M}{t} + \grad_x\cdot\int_{\R^d} g v \d v = 0 \,,
\]
and using~\fref{f1}, we finally obtain the nonlinear diffusion equation
\(   \label{diff_scaling}
    \part{\varrho_\M}{t} - \frac{d}{2}\grad_x\cdot\left(
            \frac{1}{H(\W\ast\varrho_\M)} \grad_x\left(T[\W\ast\varrho_\M]\varrho_\M \right)\right) = 0 \,.
\)
Observe that with the choice $H\equiv\mbox{const.}$, \fref{diff_scaling}
reduces to~\fref{mmodel1}, possibly up to a linear rescaling of time.

\section{Mathematical analysis of the first-order model}\label{sec:Analysis}
In this section we show the existence of weak solutions
of the first-order nonlinear diffusion equation \eqref{mmodel1} and some asymptotic properties of these solutions.
For simplicity, we consider the full space setting $\Omega=\R^d$, $d\geq 1$, but our
analysis can be easily adapted to the case of a bounded domain $\Omega$
with homogeneous Neumann boundary conditions, 
or to the case of periodic boundary conditions, as in the numerical examples
of Section~\ref{sec:Numerics}.

To simplify the notation, we set $F(z):= \frac{1}2 G(z)^2$, thus we rewrite~\eqref{mmodel1} as
\(
   \part{\varrho}{t} =  \laplace\left(F(\W\ast \varrho) \varrho\right) \,,
   \label{mmodelF}
\)
subject to the initial condition
\(
   \varrho(t=0) = \varrho_0 \,. \label{mmodelF_IC}
\)
For the rest of this Section, and without further notice,
we make the following reasonable assumptions on $F$ and $W$:
\begin{itemize}
\item $F: \R^+ \rightarrow \R^+$ is a bounded, nonnegative and nonincreasing function.
Let us note that we allow for degeneracy in~\eqref{mmodel1},
i.e., there might exist an $s_0 > 0$ such that $F(s) = 0$ for all $s>s_0$.

\item $F$ is continuously differentiable with globally Lipschitz continuous derivative
and $G=\sqrt{2F}$ is globally Lipschitz continuous.

\item $W \in W^{1,\infty}(\R^d) \cap H^2(\R^d)$.
\end{itemize}

\begin{defin}\label{def:weak_sol}
We call $\varrho\in L^\infty(0,T; \P(\Omega))$,
where $\P(\Omega)$ denotes the set of probability measures on $\Omega$,
a weak solution of~\eqref{mmodelF} subject to the initial condition \eqref{mmodelF_IC}
with $0 \leq \varrho_0 \in L^2(\Omega) \cap \P(\Omega)$,
if $\grad (F(\W\ast\varrho)\varrho) \in L^2(0,T; L^2(\Omega))$
and for every smooth, compactly supported test function $\varphi\in C^\infty_c([0,T)\times\Omega)$
we have
\(
   \int_0^\infty \int_\Omega \varrho\part{\varphi}{t} \d x\d t + \int_\Omega \varrho_0 \varphi(t=0) \d x =
       \int_0^\infty\int_\Omega \grad(F(\W\ast\varrho)\varrho)\cdot\grad\varphi \d x\d t \,.
   \label{mmodelF_weak}
\)
\end{defin}

The proof of existence of weak solutions in the sense of Definition~\ref{def:weak_sol}
will be performed in two steps. First, we consider an approximating, uniformly parabolic equation,
and prove the existence of its solutions. Then, we remove the approximation in a limiting procedure,
to obtain global in time distributional solutions.

\subsection{Approximation}\label{subsec:Approximation}
In order to obtain a uniformly positive diffusion coefficient, we use the approximation
$F_{\eps}(z) := F(z) + \eps$ for $\eps>0$,
and analyze the approximating equation
$\part{\varrho}{t} =  \laplace\left(F_{\eps}(\W\ast \varrho) \varrho\right)$,
which we write in the Fokker-Planck form
\(
   \part{\varrho}{t} =  \grad\cdot\bigl(\varrho\grad F_{\eps}(\W\ast \varrho) +  F_{\eps}(\W\ast \varrho)\grad\varrho\bigr)
   \label{mmodelFepsilon}
\)
subject to the inital condition
\(
   \varrho(t=0) = \varrho_0 \,.   \label{mmodelFepsilon_IC}
\)

\begin{theorem} \label{thm:approx}
For every $\eps > 0$, $T>0$ and $\varrho_0 \in L^2(\Omega) \cap {\cal P}(\Omega)$ there exists
a nonnegative weak solution
\[
   \varrho^\eps \in L^2(0,T;H^1(\Omega)) \cap H^1(0,T;H^{-1}(\Omega)) \cap L^\infty(0,T;{\cal P}(\Omega))
\]
of (\ref{mmodelFepsilon}) in the sense of Definition~\ref{def:weak_sol}
(with $F_\eps$ in place of $F$).
\end{theorem}

\begproof
The solution is constructed as a fixed point of the map
$\Theta: u \mapsto \varrho$, defined on the convex, bounded set
$$
   {\cal B}_R:= \{ u \in L^2((0,T)\times\Omega)\,;~ \Norm{u}_{L^\infty(0,T; L^2(\Omega))} \leq R \mbox{ and } u(t,\cdot) \in \P(\Omega) \text{ for a.e. } t \in (0,T)\}
$$
with $R>0$ to be specified later,
where $\varrho = \Theta(u)$ is the unique solution of the Fokker-Planck equation
\(  \label{Fokker-Planck}
   \part{\varrho}{t} =  \grad\cdot\bigl(\varrho\grad F_{\eps}(\W\ast u) +  F_{\eps}(\W\ast u)\grad\varrho\bigr) \,.
\)
Due to the assumed properties of $W$ and $F$, the diffusion coefficient satisfies
$F_{\eps}(\W\ast u) \in L^\infty((0,T)\times\Omega)$ and is bounded below by $\eps$,
and the convection coefficient $\grad F_{\eps}(\W\ast u) \in L^\infty((0,T)\times\Omega)$.
Thus, \eqref{Fokker-Planck} is a linear, uniformly parabolic equation with bounded coefficients,
and the standard theory~\cite{Evans} implies the existence and uniqueness of
the global, nonnegative weak solution 
\(
   \varrho \in L^\infty(0,T;L^2(\Omega)) \cap L^2(0,T;H^1(\Omega)) \cap 
   H^1(0,T;H^{-1}(\Omega)),   \label{varrho}
\)
such that $\varrho(\cdot,t) \in {\cal P}(\Omega)$ for almost every $t\geq 0$.
Therefore, $\Theta$ is well defined on ${\cal B}_R$ and, as can be easily proven,
is continuous there. 

In order to apply the Schauder fixed point theorem (cf.~\cite{Brezis}),
we need to show that $\Theta$ maps ${\cal B}_R$ into a compact subset of itself.
The key point is the following a-priori estimate, obtained by
using $\varrho$ as a test function for~\eqref{Fokker-Planck}
(note that due to~\eqref{varrho}, $\varrho$ is indeed an admissible test function):
\[
   \int_0^s \int_\Omega \varrho \part{\varrho}{t} \d x \d t  + \int_\Omega \varrho_0^2 \d x
   + \int_0^s \int_\Omega F_{\eps} |\grad\varrho|^2 \d x \d t =
     - \int_0^s \int_\Omega \varrho \grad F_{\eps} \cdot \grad \varrho \d x \d t \,,
\]
for $s\in(0,T)$, where we introduced the shorthand notation $F_\eps := F_\eps(\W\ast u)$.
Now we use the identity
\[
   \int_0^s \int_\Omega \varrho \part{\varrho}{t} \d x \d t =
     \frac12 \int_0^s \frac{d}{dt} \left(\int_\Omega \varrho^2 \d x\right) \d t =
     \frac12 \int_\Omega \varrho(x,s)^2 \d x - \frac12 \int_\Omega \varrho_0^2 \d x 
\]
and the Young inequality
\[
    \left| \int_0^s \int_\Omega \varrho \grad F_{\eps} \cdot \grad \varrho \d x \d t \right| \leq 
    \frac{1}2 \int_0^s \int_\Omega F_{\eps} |\grad\varrho|^2 \d x \d t + 
    2 \int_0^s \int_\Omega \left|\nabla \sqrt{F_{\eps}}\right|^2 \varrho^2 \d x \d t 
\]
to obtain
\[
   \int_\Omega \varrho(x,s)^2 \d x + \frac12 \int_0^s \int_\Omega F_{\eps} |\grad\varrho|^2 \d x \d t \leq
   \int_\Omega \varrho_0^2 \d x  +  \int_0^s \int_\Omega \left|\grad\sqrt{F_{\eps}}\right|^2 \varrho^2 \d x \d t \,.
\]
Finally, since $F_\eps \geq \eps$ and $|F'(y)|\leq C$
for $0 \leq y \leq \sup_{x\in\Omega} \W\ast u(x) \leq \Norm{\W}_{L^\infty}$,
\[
   \left|\grad\sqrt{F_\eps}\right|^2 = 
   \left|\nabla \sqrt{F_{\eps}(W*u)}\right|^2 = 
   \frac{|F'(\W\ast u)|^2}{4 F_\eps(\W\ast u)} |\nabla \W\ast u|^2 \leq 
   \frac{C^2}{4\eps} \Norm{\W}^2_{W^{1,\infty}} =: C_\eps,
\]
and we conclude
\[
   \int_\Omega \varrho(x,s)^2~dx + \frac{\eps}{2} \int_0^s \int_\Omega |\nabla \varrho|^2 \d x \d t
       \leq  \int_\Omega \varrho_0^2 \d x  +  C_\eps \int_0^s \int_\Omega \varrho^2 \d x \d t \,.
\]
The Gronwall inequality yields an a-priori bound on $\varrho$ in $L^\infty(0,T;L^2(\Omega)$,
which implies that $\Theta$ indeed maps ${\cal B}_R$ into itself, for a suitable $R>0$.
Moreover, reinserting into the right-hand side, we obtain an a-priori bound in $L^2(0,T;H^1(\Omega))$,
and an analysis of the right-hand side of \eqref{Fokker-Planck} immediately yields
a~bound in $H^1(0,T;H^{-1}(\Omega))$.
The proof is concluded by an application of the Aubin-Lions Lemma (cf.~\cite{Brezis}),
which states that the space $L^2(0,T;H^1(\Omega)) \cap H^1(0,T;H^{-1}(\Omega))$
is compactly embedded into $L^2((0,T)\times\Omega)$,
so that $\Theta$ maps ${\cal B}_R$ into a compact subset of itself.
\endproof

Most a-priori estimates on $\varrho^\eps$
that have been obtained along the lines of the above proof are not uniform in $\eps$
and can thus not be used to pass to the limit $\eps\to 0$.
However, we can tackle this problem with a more careful analysis,
where the key idea is to look for estimates on $\sqrt{F_\eps} \varrho^\eps$,
rather than on $\varrho^\eps$.

\begin{lemma} \label{lemma:estimates}
There exists a constant $C$ independent of $\eps>0$ small enough,
such that the solutions $\varrho^\eps$ of (\ref{mmodelFepsilon}),
constructed in Theorem~\ref{thm:approx}, satisfy
\[
   \Norm{\varrho^\eps}_{L^\infty(0,T;\P(\Omega))} \leq 1 \,,\qquad
   \Norm{\varrho^\eps}_{L^\infty(0,T;L^2(\Omega))} \leq C \,,\qquad
   \Norm{\part{\varrho^\eps}{t}}_{L^2(0,T;H^{-1}(\Omega))} \leq C \,,
\]
and
\[
   \Norm{\W\ast\varrho^\eps}_{L^\infty(0,T;W^{1,\infty}(\Omega))} \leq C \,,\qquad
   \Norm{\W\ast\part{\varrho^\eps}{t}}_{L^\infty(0,T;L^2(\Omega))} \leq C \,.
\]
\end{lemma}

\begproof
For the sake of better legibility, we will drop the superscript at $\varrho^\eps$ in the proof.
Since $\varrho$ is constructed as a nonnegative weak solution
of the Fokker-Planck equation~\eqref{mmodelFepsilon},
the first estimate follows immediately due to the mass conservation
\[
   \Norm{\varrho(t,\cdot)}_{L^1(\Omega)} = 
   \int_\Omega \varrho(t,x)~dx =
   \int_\Omega \varrho_0(x) \d x = 1 \,.
\]

Proceeding in a similar way as in the proof of Theorem~\ref{thm:approx},
we conclude for $s \in (0,T)$
\[
   \frac12 \int_\Omega \varrho(s,x)^2 \d x +
     \int_0^s \int_\Omega F_{\eps} |\grad\varrho|^2 \d x \d t +
     \int_0^s \int_\Omega \varrho \grad F_\eps\cdot\grad\varrho \d x \d t = \int_\Omega \varrho_0^2 \d x \,,
\]
where we again use the shorthand notation $F_\eps := F_\eps(\W\ast\varrho)$.
Using the identity
\[
   \varrho \grad F_\eps\cdot\grad\varrho + F_\eps |\grad\varrho|^2 =
   \left|\grad\left(\sqrt{F_\eps}\varrho\right)\right|^2
      - \left|\varrho\grad\sqrt{F_\eps}\right|^2  \,,
\]
we obtain
\[
   \frac12 \int_\Omega \varrho(s,x)^2 \d x +
      \int_0^s \int_\Omega \left|\grad \left(\sqrt{F_\eps}\varrho\right)\right|^2 \d x \d t
   =  \frac12 \int_\Omega \varrho_0^2 \d x +
      \int_0^s \int_\Omega \varrho^2 \left|\grad\sqrt{F_\eps}\right|^2 \d x \d t \,.
\]
Now, since $|G'(y)|\leq C$
for $0 \leq y \leq \sup_{x\in\Omega} \W\ast u(x) \leq \Norm{\W}_{L^\infty}$,
we have
\[
   \left| \grad\sqrt{F_{\eps}} \right| =
   \left| \grad\sqrt{F_{\eps}(W\ast\varrho)} \right| =
   |G'(W\ast\varrho)| |\nabla W\ast\varrho| \leq
   C \Norm{W}_{W^{1,\infty}} \,.
\]
Hence,
\[
   \frac12 \int_\Omega \varrho(s,x)^2 \d x +
   \int_0^s \int_\Omega \left|\grad \left(\sqrt{F_\eps} \varrho\right)\right|^2 \d x \d t
   \leq  \frac12 \int_\Omega \varrho_0^2 \d x +
         C^2 \Norm{W}^2_{W^{1,\infty}} \int_0^s \int_\Omega \varrho^2 \d x \d t \,,
\]
and the Gronwall inequality yields a uniform estimate for $\varrho$ in
$L^\infty(0,T;L^2(\Omega))$, which subsequently implies a uniform estimate for 
$\sqrt{F_\eps(W\ast\varrho)} \varrho$ in $L^2(0,T;H^1(\Omega))$
and, subsequently, also for $\part{\varrho}t$ in
$L^2(0,T;H^{-1}(\Omega)$.

Finally, we derive the estimates for $\W\ast\varrho$.
Since $\W\in W^{1,\infty}(\R^d)$, we immediately obtain
\[
   \Norm{\W\ast\varrho}_{L^\infty(0,T; W^{1,\infty}(\Omega))} \leq
   \Norm{\W}_{W^{1,\infty}(\R^d)} \sup_{t\in[0,T]} \int_\Omega \left|\varrho(t,x)\right| \d x =
   \Norm{\W}_{W^{1,\infty}(\R^d)} \,.
\]
Moreover, we have
\[
   \W\ast\part{\varrho}{t} = \W\ast\laplace(F_\eps \varrho) = 
   \laplace \W\ast (F_\eps \varrho)   \,,
\]
and with  $\W\in H^2(\R^d)$, $\varrho \in L^\infty(0,T; {\cal P}(\Omega))$
and the uniform boundedness of $F_\eps$ in $L^\infty((0,T)\times\Omega)$,
we obtain a uniform bound for $\W\ast(\part{\varrho}t)$
in $L^\infty(0,T; L^2(\Omega))$.
\endproof

\subsection{Global existence}\label{sec:GlobalExistence}
From the above approximation and a-priori estimates
we can easily pass to the limit $\eps\to 0$
and conclude the existence of a weak solution for \eqref{mmodelF}:

\begin{theorem} \label{thm:GlobalExistence}
Let $\varrho_0 \in L^2(\Omega) \cap {\cal P}(\Omega)$ and $T>0$. Then there exists a solution
\[
   \varrho \in L^\infty(0,T;L^2(\Omega)) \cap L^\infty(0,T;{\cal P}(\Omega)) \cap H^1(0,T;H^{-1}(\Omega))
\]
of \eqref{mmodelF}--\eqref{mmodelF_IC} in the sense of the weak formulation~\eqref{mmodelF_weak},
such that in addition
\[
   \sqrt{F(\W\ast\varrho)} \varrho \in L^2(0,T;H^1(\Omega)) \,.
\]
\end{theorem}

\begproof
We use the uniform estimates of Lemma~\ref{lemma:estimates}
and the Banach-Alaoglu theorem~\cite{Brezis} to extract a subsequence $\varrho^{\eps_n}$ such that
\begin{eqnarray*}
   \varrho^{\eps_n} \rightharpoonup^* \varrho && \quad \text{weakly-* in } L^\infty(0,T;L^2(\Omega)) \,, \\
   \part{\varrho^{\eps_n}}t \rightharpoonup \part{\varrho}t && \quad  \text{weakly in } L^2(0,T;H^{-1}(\Omega)) \,, \\
   \sqrt{F_{\eps_n}(\W\ast\varrho^{\eps_n}) } \varrho^{\eps_n} \rightharpoonup u && \quad \text{weakly in } L^2(0,T;H^{1}(\Omega)) \,,
\end{eqnarray*}
for some $u\in L^2(0,T;H^{1}(\Omega))$.
Moroever, we use a variant of the Aubin-Lions Lemma~\cite{Simon}
to conclude the compact embedding of $L^\infty(0,T;W^{1,\infty}(\Omega)) \cap W^{1,\infty}(0,T;L^2(\Omega))$
into $L^\infty(0,T;C(\Omega))$.
Thus, we can extract a further subsequence, again denoted by $\varrho^{\eps_n}$, such that as $\eps_n\to 0$,
\[
   \W\ast\varrho^{\eps_n} \rightarrow v \quad \text{strongly in } L^\infty(0,T;C(\Omega)) \,.
\]
Since further
\[
   \W\ast\varrho^{\eps_n} \rightharpoonup^* \W\ast\varrho \quad \text{weakly-* in } L^\infty(0,T;L^2(\Omega)) \,,
\]
we conclude $v = \W\ast\varrho$ by the uniqueness of the limit.
Due to the continuity properties of $F$ and $F'$, we also have
\begin{eqnarray*}
   {F_{\eps_n}(\W\ast\varrho^{\eps_n}) }  \rightarrow {F(\W\ast\varrho)} && \quad  \text{strongly in }L^\infty(0,T;C(\Omega)) \,, \\
   \sqrt{F_{\eps_n}(\W\ast\varrho^{\eps_n}) }  \rightarrow \sqrt{F(\W\ast\varrho)} && \quad  \text{strongly in }L^\infty(0,T;C(\Omega)) \,, \\
   {F_{\eps_n}'(\W\ast\varrho^{\eps_n}) }  \rightarrow {F'(\W\ast\varrho)} && \quad  \text{strongly in }L^\infty(0,T;C(\Omega)) \,.
\end{eqnarray*} 
Consequently,
\[
   \sqrt{F_{\eps_n}(\W\ast\varrho^{\eps_n}) } \varrho^{\eps_n} \rightharpoonup \sqrt{F(\W\ast\varrho)} \varrho \qquad \text{weakly-* in } L^\infty(0,T;L^2(\Omega)) \,,
\]
and, again, by the uniqueness of the limit we identify $u = \sqrt{F(\W\ast\varrho)}\varrho$.
Finally, we use the identity
\[
   \grad \left(F_{\eps_n} \left(\W\ast\varrho^{\eps_n}\right) \varrho^{\eps_n} \right) =
    F_{\eps_n}' \left(\W\ast\varrho^{\eps_n}\right) \varrho^{\eps_n} \,\grad\left(\W\ast\varrho^{\eps_n}\right)
 + \sqrt{F_{\eps_n}\left(\W\ast\varrho^{\eps_n}\right)} \,\grad \left(\sqrt{F_{\eps_n}\left(\W\ast\varrho^{\eps_n}\right)} \varrho^{\eps_n} \right)
\]
and the above limits to conclude
\[
   \grad \left(F_{\eps_n} \left(\W\ast\varrho^{\eps_n}\right) \varrho^{\eps_n} \right)
    \rightharpoonup \grad \left(F \left(\W\ast\varrho\right) \varrho\right) \quad \text{weakly in } L^2(0,T;L^2(\Omega)) \,.
\]
Hence, we can pass to the limit $\eps_n\to 0$ in the weak formulation~\eqref{mmodelF_weak}
and conclude that $\varrho$ is a weak solution of \eqref{mmodelF}--\eqref{mmodelF_IC}.
\endproof

Finally, we want to reduce the regularity of the initial value towards solely probability measures.
This is relevant since we shall see below that indeed Dirac $\delta$ distributions can be stationary solutions.
For this sake we define the very weak (distributional) notion of the solution:

\begin{defin}\label{def:very_weak_sol}
We call $\varrho\in L^\infty(0,T; \P(\Omega))$
a very weak (distributional) solution of~\eqref{mmodelF} subject to the initial condition \eqref{mmodelF_IC}
with $\varrho_0 \in \P(\Omega)$,
if for every smooth, compactly supported test function $\varphi\in C^\infty_c([0,T)\times\Omega)$
we have
\(
   \int_0^\infty \int_\Omega \part{\varphi}{t} \varrho \d x\d t + \int_\Omega \varphi(t=0) \varrho_0 \d x =
       - \int_0^\infty\int_\Omega (F(\W\ast\varrho) \laplace\varphi) \varrho \d x\d t \,,
   \label{mmodelF_very_weak}
\)
where we denote by $\varrho\d x$ the integration with respect to the probability measure $\varrho(t, \cdot)$.
\end{defin}

\begin{theorem}
Let $\varrho_0 \in {\cal P}(\Omega)$ and $T>0$.
Then there exists a very weak solution
\[
   \varrho \in L^\infty(0,T;{\cal P}(\Omega)) 
\]
of \eqref{mmodelF}--\eqref{mmodelF_IC} in the sense of~\eqref{mmodelF_very_weak}.
\end{theorem}

\begproof
Let us consider a sequence $\left(\varrho^n_0\right)_{n\in\N} \subseteq L^2(\Omega)\cap\P(\Omega)$
such that $\varrho^n_0 \to \varrho_0$ in $\P(\Omega)$.
Moreover, let us denote by $\varrho^n$ the corresponding weak solutions
of~\eqref{mmodelF}--\eqref{mmodelF_IC}, which are trivially
also the very weak solutions in the sense of~\eqref{mmodelF_very_weak}.
Conservation of mass provides immediately a uniform bound on $\varrho^n$
in $L^\infty(0,T;\P(\Omega))$, such that for a subsequence, again denoted by $\varrho^n$,
we have
\[
    \varrho^n \rightharpoonup^* \varrho \quad\mbox{weakly-* in } L^\infty(0,T;\P(\Omega)) \,.
\]
To show that $\varrho$ is a solution to~\eqref{mmodelF_very_weak},
we need to prove that $F(\W\ast\varrho^n)$ converges to
$F(\W\ast\varrho)$ strongly in $L^1(0,T;C(\Omega))$.
First, let us note that, due to the assumptions on $\W$,
$\W\ast\varrho^n$ is uniformly bounded in $L^\infty(0,T;W^{1,\infty}(\Omega))$,
and the same holds for $F(\W\ast\varrho^n)$.
Consequently, the sequence $F(\W\ast\varrho^n)\varrho^n$
is uniformly bounded in $L^\infty(0,T; \P(\Omega))$.
From~\eqref{mmodelF} it immediately follows that
$\part{\varrho^n}{t} = \laplace(F(\W\ast\varrho^n)\varrho^n)$
is uniformly bounded in $L^\infty(0,T; (C_c^2(\Omega))^*)$,
where $(C_c^2(\Omega))^*$ denotes the dual space to $C_c^2(\Omega)$,
the space of twice continuously differentiable functions
with compact support in $\Omega$.
Consequently, the sequence $\part{}{t} F(\W\ast\varrho^n)
= F'(\W\ast\varrho^n)\W\ast\part{\varrho^n}{t}$
is uniformly bounded in $L^\infty(0,T; (C_c^2(\Omega))^*)$,
and the generalization of the Aubin-Lions lemma by Simon~\cite{Simon}
implies then the strong convergence of $F(\W\ast\varrho^n)$
to $F(\W\ast\varrho)$ in $L^1(0,T; C(\Omega))$.
\endproof

\subsection{Stationary Solutions}\label{subsec:StatSols}
The numerical simulations provided in Section~\ref{sec:Numerics} suggest that for $G>0$
(for instance, $G(s) = e^{-s}$) and bounded domain with periodic boundary conditions,
the stationary solutions of~\eqref{mmodel1} consist of one or more localized, but not compactly supported
aggregates (clumps), see the bottom right panel of Fig.~\ref{fig:num_mm1} and Fig.~\ref{fig:num_mm3}.
However, we were not able to characterize these stationary aggregates analytically.
Instead, we provide a few examples of other types of stationary solutions,
posed either in the full space setting $\Omega=\R^d$ or on a torus $\Omega=\T^d$
with periodically extended $\W$. These examples are rather trivial, however,
still provide an interesting insight into the relatively rich structure
of the solutions of~\eqref{mmodel1}.

\begin{itemize}
\item
The most trivial type of stationary solution is the constant state $\varrho\equiv c > 0$.
Clearly, this has finite mass only if $\Omega = \T^d$.
\item
If there exists an $s_0 > 0$ such that $G(s) \equiv 0$ for all $s\geq s_0$,
then any profile $\varrho$ such that $\varrho \geq s_0 \left(\int_\Omega \W(x)\d x\right)^{-1}$
almost everywhere on $\Omega$
is a stationary solution to the distributional formulation~\eqref{mmodelF_very_weak}.
Indeed, such a solution satisfies $\W\ast\varrho \geq s_0$,
so that $G(\W\ast\varrho) \equiv 0$.
However, again, this solution has infinite mass if $\Omega = \R^d$.
\item
If there exists an $s_0 > 0$ such that $G(s) \equiv 0$ for all $s\geq s_0$,
and, moreover, $\W$ is continuous on $\Omega$, we construct the atomic measure
\[
    \varrho(x) = \sum_{i=1}^N c_i \delta(x-x_i)
\]
for some $N\in\N$, $x_i\in\Omega$ and $c_i > 0$
such that $c_i \W(0) > s_0$ for all $i=1,\dots,N$.
Then $\varrho$ is a distributional stationary solution to~\eqref{mmodelF_very_weak}.
Indeed, for any $i=1,\dots,N$ we have
\[
   \W\ast\varrho(x_i) = \sum_{j=1}^N c_j \W(x_i - x_j) \geq c_i \W(0) > s_0 \,.
\]
By the continuity of $\W$ and, consequently, $\W\ast\varrho$,
we have $\W\ast\varrho > s_0$ on some neighbourhood of $x_i$.
Therefore, $G(\W\ast\varrho) \equiv 0$ on some open set containing $\bigcup_{i=1}^N x_i$ and, finally,
$G(\W\ast\varrho) \varrho \equiv 0$ everywhere on $\Omega$.
\end{itemize}

\subsection{Linear stability analysis}\label{sec:Stability}
In this section we perform a linear stability analysis of constant density states
of the nonlinear diffusion equation~\fref{mmodel1}.
As discussed previously, we have constant steady state solutions
either in the full-space setting $\Omega=\R^d$ (however, with infinite mass),
or on a bounded domain $\Omega\subset\R^d$ with periodic boundary conditions.
Without loss of generality, we assume the normalization $\int_\Omega W(x) \d x = 1$.
Let us make the perturbation ansatz $\varrho = \varrho_0 + \eps\tilde\varrho$,
where $\varrho_0>0$ is a constant state such that $G(\varrho_0) > 0$ and $G'(\varrho_0) < 0$.
Then $\W\ast\varrho = \varrho_0 + \eps\W\ast\tilde\varrho$ and
assuming sufficient smoothness of $G$, we have
\[
    G(\W\ast\varrho)^2 = G(\varrho_0)^2 + 2\eps G(\varrho_0) G'(\varrho_0) \W\ast\tilde\varrho + O(\eps^2) \,.
\]
Inserting the ansatz into~\fref{mmodel1} and collecting terms of order $\O(\eps)$, we arrive at
\[
    \part{\tilde\varrho}{t} - \frac{G(\varrho_0)}{2} \left(
              G(\varrho_0)\laplace\tilde\varrho + 2G'(\varrho_0)\varrho_0 \laplace(\W\ast\tilde\varrho)
         \right) = 0 \,.
\]
Performing the Fourier transform, we obtain
\(   \label{Fourier}
    \part{\hat{\varrho}}{t} + |\xi|^2 \frac{G(\varrho_0)}{2}
         \left( G(\varrho_0) + 2G'(\varrho_0)\varrho_0 \hat{\W}\right) \hat{\varrho} = 0 \,,
\)
where we denoted $\hat{\varrho} = \hat{\varrho}(t,\xi)$ the Fourier transform of $\tilde\varrho$.
Consequently, with the assumption $G(\varrho_0) > 0$ and $G'(\varrho_0) < 0$,
those wavenumbers $\xi$ of $\tilde\varrho$ are stable for which
\(   \label{stability}
    \mbox{Re }\hat{\W}(\xi) < - \frac{G(\varrho_0)}{2G'(\varrho_0)\varrho_0} \geq 0 \,.
\)
Since $\W\in L^1(\R^d)$, we have $\hat{\W} \in C^0(\R^d)$, and, therefore,
all wavenumbers of $\tilde\varrho$ larger than a certain threshold
will be stable.
On the other hand, wavenumbers violating~\fref{stability}
will lead to pattern formation, as we show in our numerical examples,
Section~\ref{sec:Numerics}.
Moreover, a quick inspection of~\fref{Fourier}
leads to the expectation that the stable high wavenumbers
will be smoothed-out on a faster time scale,
compared to the slower emergence of patterns
due to the unstable lower wavenumbers.
Finally, considering the scalings $\W_r(x) := r^{-d}\W_1(x/r)$
of a fixed kernel $\W_1\in L^1(\R^d)$,
we have $\hat{\W_r}(\xi) = \hat{W_1}(r\xi)$
and, therefore, we expect that the wavelength of the patterns
(size of aggregates) will scale with $r$.
This can be as well clearly observed in our numerical examples.

It is interesting to consider the formal extremal case $\W=\delta_0$,
i.e., $\W\ast\tilde\varrho = \tilde\varrho$.
Then $\hat{\W} \equiv 1$ and we conclude that
the constant state $\varrho_0$ is stable if and only if
\(   \label{stability2}
   (G(\varrho_0)^2 \varrho_0)' = G(\varrho_0)^2 + 2 G(\varrho_0) G'(\varrho_0)\varrho_0 > 0  \,.
\)
In fact, violation of~\fref{stability2} with $\W=\delta_0$ leads to
an ill-posed problem, since then~\fref{mmodel1}
looks like a backward heat equation, which can be seen
if we expand the derivatives and write it in the Fokker-Planck form as
\(   \label{ill-posed}
    \part{\varrho}{t} = \frac12 \grad\cdot\left( G(\varrho) [2G'(\varrho)\varrho + G(\varrho)]\grad\varrho\right) \,.
\)
The equation is parabolic (and thus well-posed)
only if the diffusivity $2G(\varrho) G'(\varrho)\varrho + G^2(\varrho)$ is strictly positive,
which is our stability condition~\fref{stability2}.
Therefore, only imposing an initial condition $\varrho_0$
uniformly satisfying~\fref{stability2} leads to a well-posed diffusion equation for all $t\geq 0$,
since~\fref{stability2} is preserved due to the maximum principle.
On the other hand, if $G>0$, the nonlocality $\W\in L^1(\R^d)$ always stabilizes the equation
in the sense of~\fref{stability}.
Indeed, writing it in the Fokker-Planck form
\[
   \part{\varrho}{t} = \frac12 \grad\cdot\left(Individual based and mean-field modelling of direct aggregation
          \left[2G(W\ast\varrho)G'(W\ast\varrho)\grad W\ast\varrho\right]\varrho + G^2(W\ast\varrho)\grad\varrho\right) \,,
\]
we see that the second-order term appears with the positive diffusivity $G^2(W\ast\varrho)$.
The first-order term describes then the transport of $\varrho$ along
the generalized gradient $\grad W\ast\varrho$ and is responsible for the
aggregative effect.

Clearly, the ill-posedness of~\eqref{ill-posed}
can be avoided by merely introducing a nonlocality in the first-order transport term,
while the diffusivity may stay local (and possibly degenerate).
Such a model was derived and studied in~\cite{Bertozzi},
which with our notation is written as
\[
    \part{\varrho}{t} = \grad\cdot\left( - \varrho\grad\W\ast\varrho + \varrho^2\grad\varrho \right) \,.
\]
This equation was constructed as a model for biological aggregations in which individuals
experience long-range social attraction and short range dispersal.
Let us note that here, in contrast to our model, the diffusivity is an increasing
function of the density $\varrho$.
In~\cite{Bertozzi} it was shown that it produces strongly nonlinear states
with compact support and steep edges that correspond to localized biological aggregations, or clumps.
Similarly as can be observed in our numerical simulations in Section~\ref{sec:Numerics},
these clumps are approached through a dynamic coarsening process.

Another insight into the stabilizing effect of the nonlocality is provided
by the introduction of a formal expansion of $\W$.
Taking $\W$ as the standard mollifier and $\W_\eps(s) := \eps^{-d}\W(s/\eps)$,
such that $\W_\eps\to \delta_0$ as $\eps\to 0$, we expand
\[
    \W_\eps\ast\varrho &=& \eps^{-d} \int_{\R^d} \W\left(\frac{x-y}{\eps}\right)\varrho(y) \d y \\
            &=& \int_{\R^d} \W(-z) \varrho(x+\eps z) \d z \\
            &=& \int_{\R^d} \W(-z) \left( \varrho(x) + \eps\grad\varrho(x)z 
                  + \frac12 \eps^2 z^t\grad^2\varrho(x) z \right) \d z + \O(\eps^3) \,.
\]
Now, due to the symmetry $\W(-z) = \W(z)$ and the normalization $\int_{\R^d} \W(z) \d z = 1$,
we have
\[
    \W_\eps\ast\varrho = \varrho + \frac{\eps^2}{2}\beta\laplace\varrho + \O(\eps^4) \,,
\]
where the constant $\beta>0$ is such that $\int_{\R^d} \W(z) z_i z_j \d z = \beta\delta_{ij}$.
Inserting this into~\fref{mmodel1}, we obtain
\[
    \partial_t\varrho &=& \frac12 \laplace\left(
         G\left(\varrho + \frac{\eps^2}{2}\beta\laplace\varrho + \O(\eps^4)\right)^2\varrho \right) \\
        &=& \frac12 \laplace \left( G(\varrho)^2\varrho
            + \eps^2\beta G(\varrho) G'(\varrho)\varrho\laplace\varrho \right) + \O(\eps^4) \,.
\]
Up to the terms of fourth order in $\eps$, this is a Cahn-Hilliard-type equation
which is well-posed for every $\eps > 0$ since the term
$\eps^2\beta G(\varrho) G'(\varrho)\varrho$ is strictly negative if $G(\varrho) > 0$ and $G'(\varrho)<0$.
However, if $\eps = 0$, i.e. $\W=\delta_0$, this regularizing effect is lost.

\section{Numerical examples}\label{sec:Numerics}
In this section we present numerical examples for the first and second-order
individual based models~\eqref{model1}, \eqref{model2-1}--\eqref{model2-2} in 2D,
the first-order mean-field limit~\eqref{mmodel1} in 1D and 2D
and the second-order mean-field limit~\eqref{mmodel2} in 1D with 1D velocity.

\subsection{The first-order individual based model~\fref{model1}}
\begin{figure}
{\centering \begin{tabular}[h]{ccc}
\resizebox*{0.3\linewidth}{!}{\includegraphics{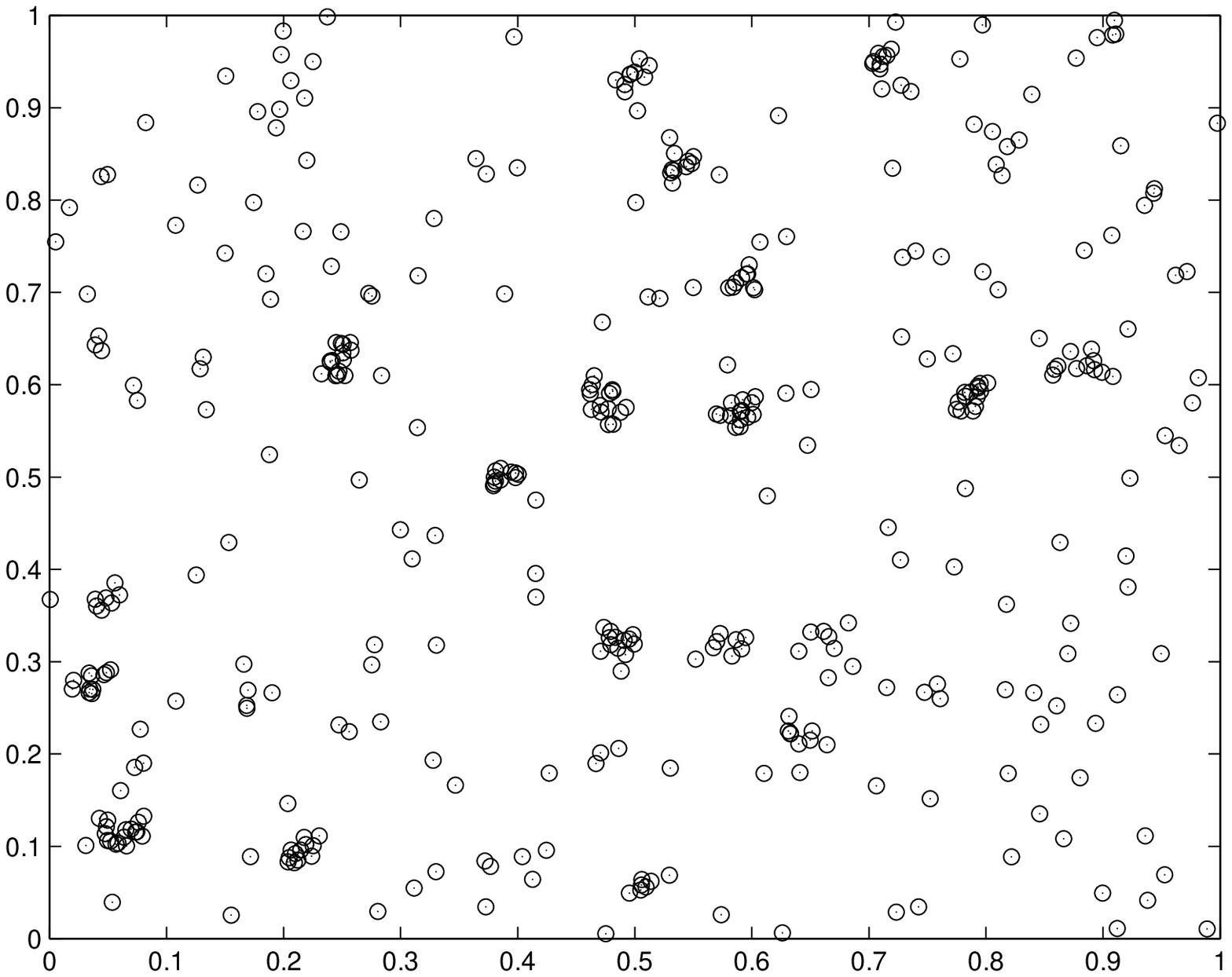}} &
\resizebox*{0.3\linewidth}{!}{\includegraphics{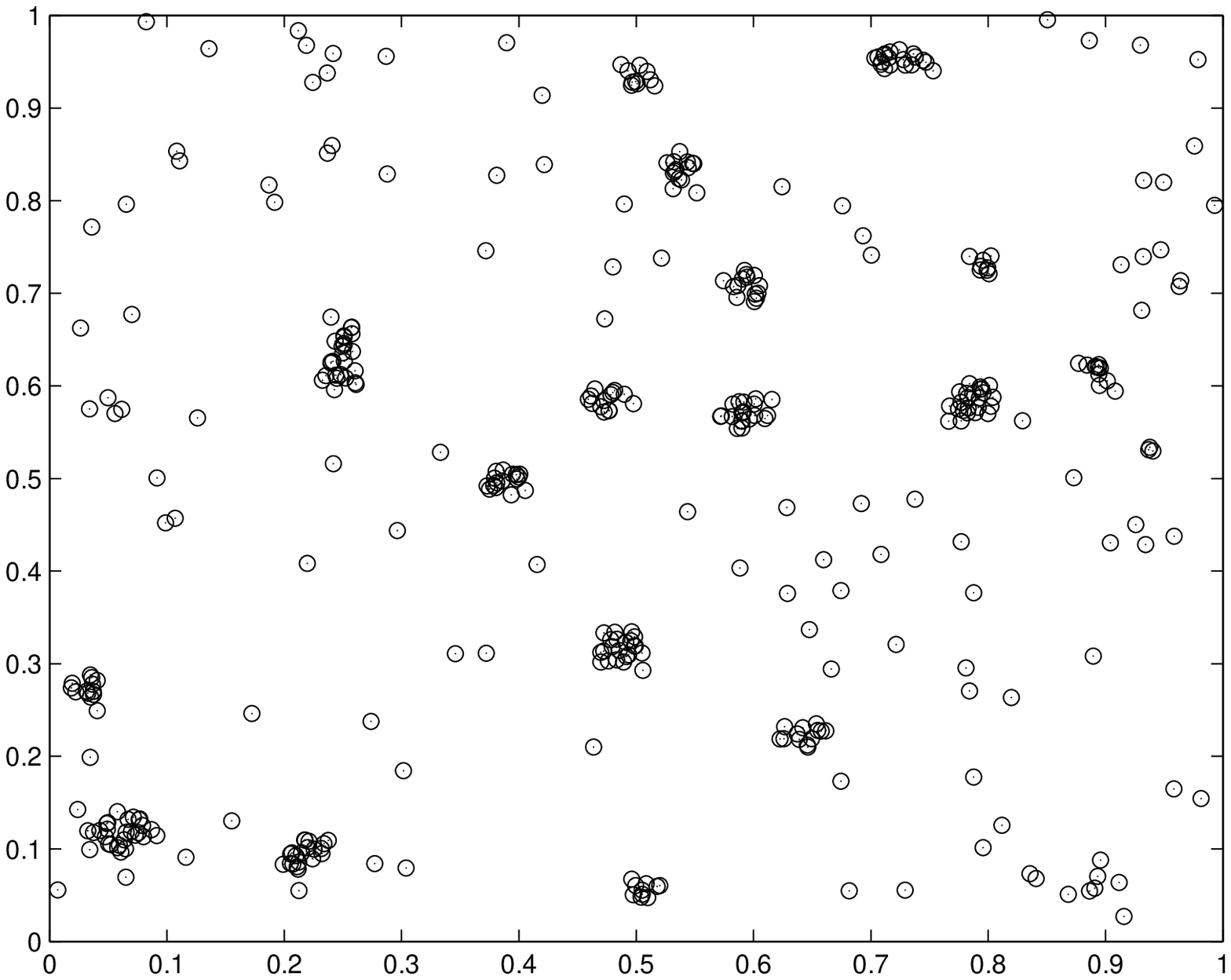}} &
\resizebox*{0.3\linewidth}{!}{\includegraphics{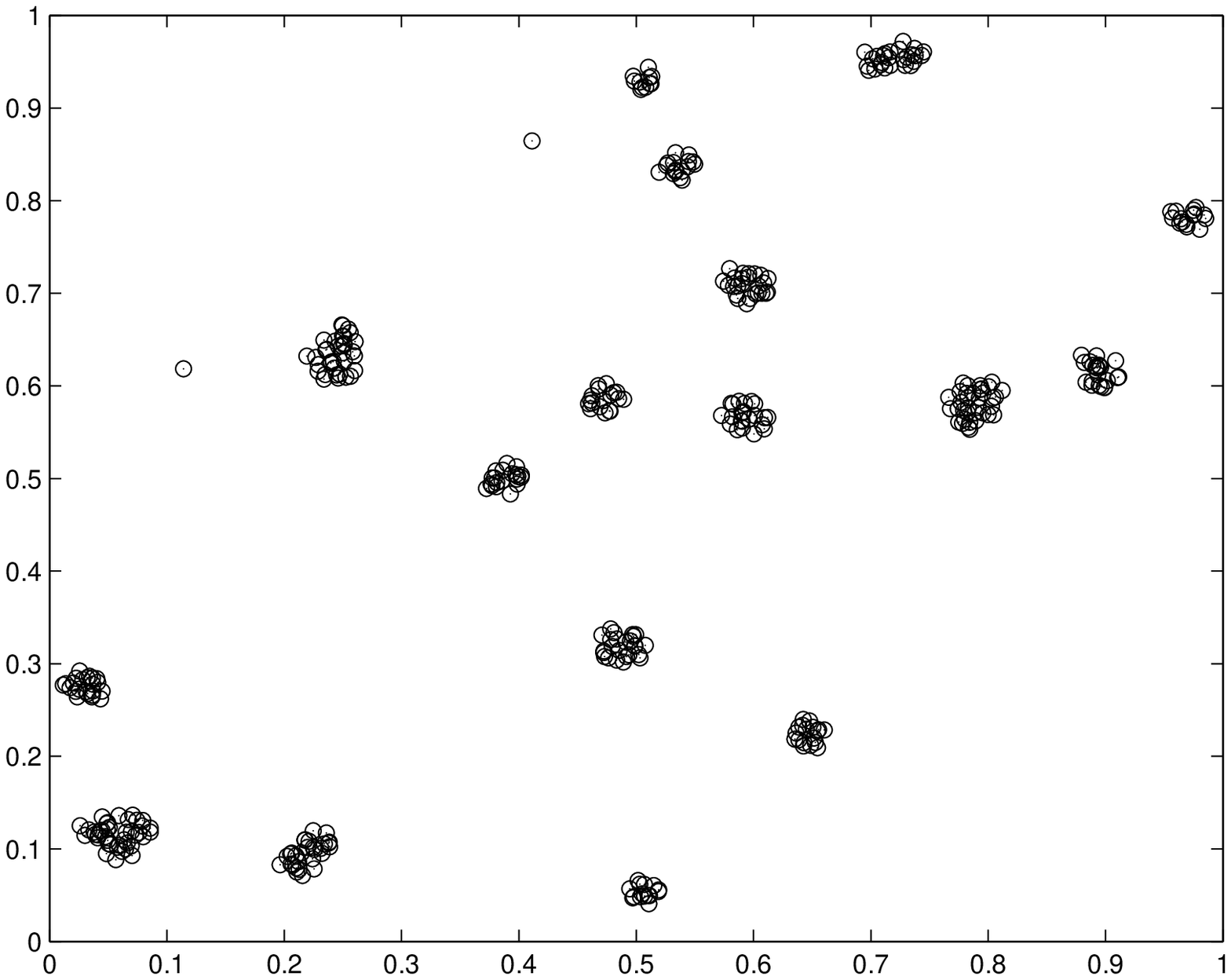}} \\
t=0.1 & t=0.2 & t=1.7
\end{tabular}\par}
\caption{The first-order individual based model with $N=400$ agents and sampling radius $R=0.025$,
subject to a random initial condition. 
\label{fig:num1}}
\end{figure}

\begin{figure}
{\centering \begin{tabular}[h]{ccc}
\resizebox*{0.3\linewidth}{!}{\includegraphics{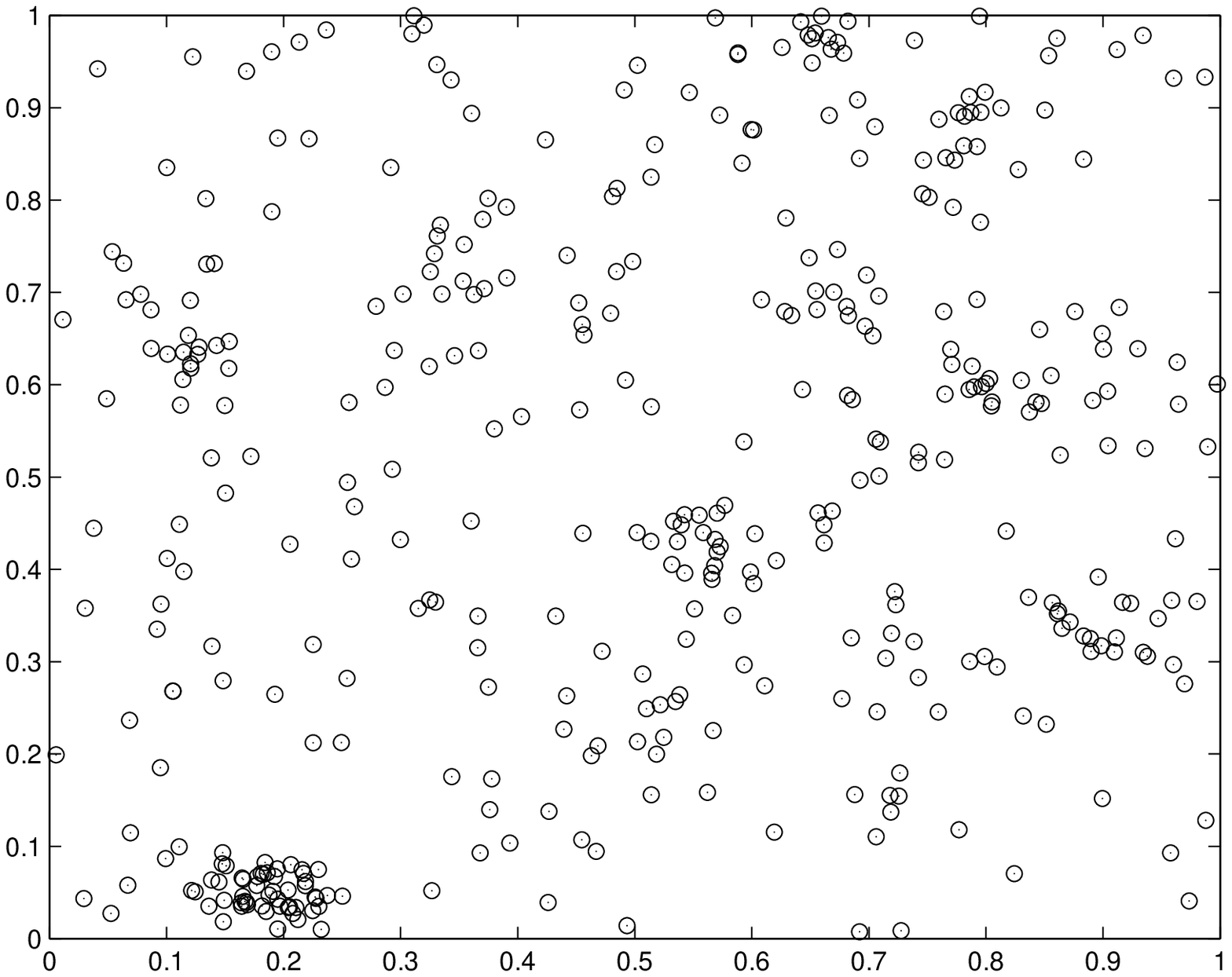}} &
\resizebox*{0.3\linewidth}{!}{\includegraphics{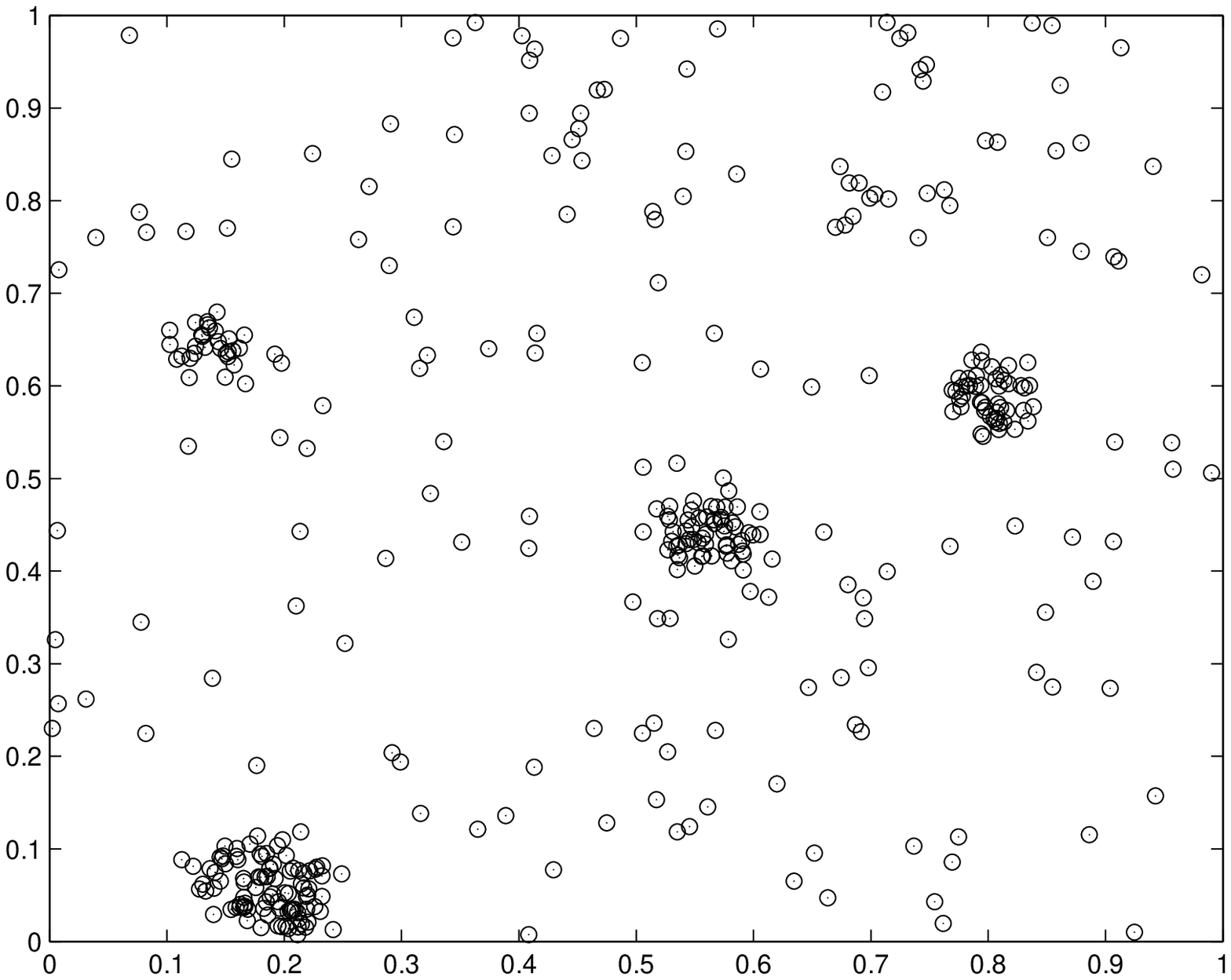}} &
\resizebox*{0.3\linewidth}{!}{\includegraphics{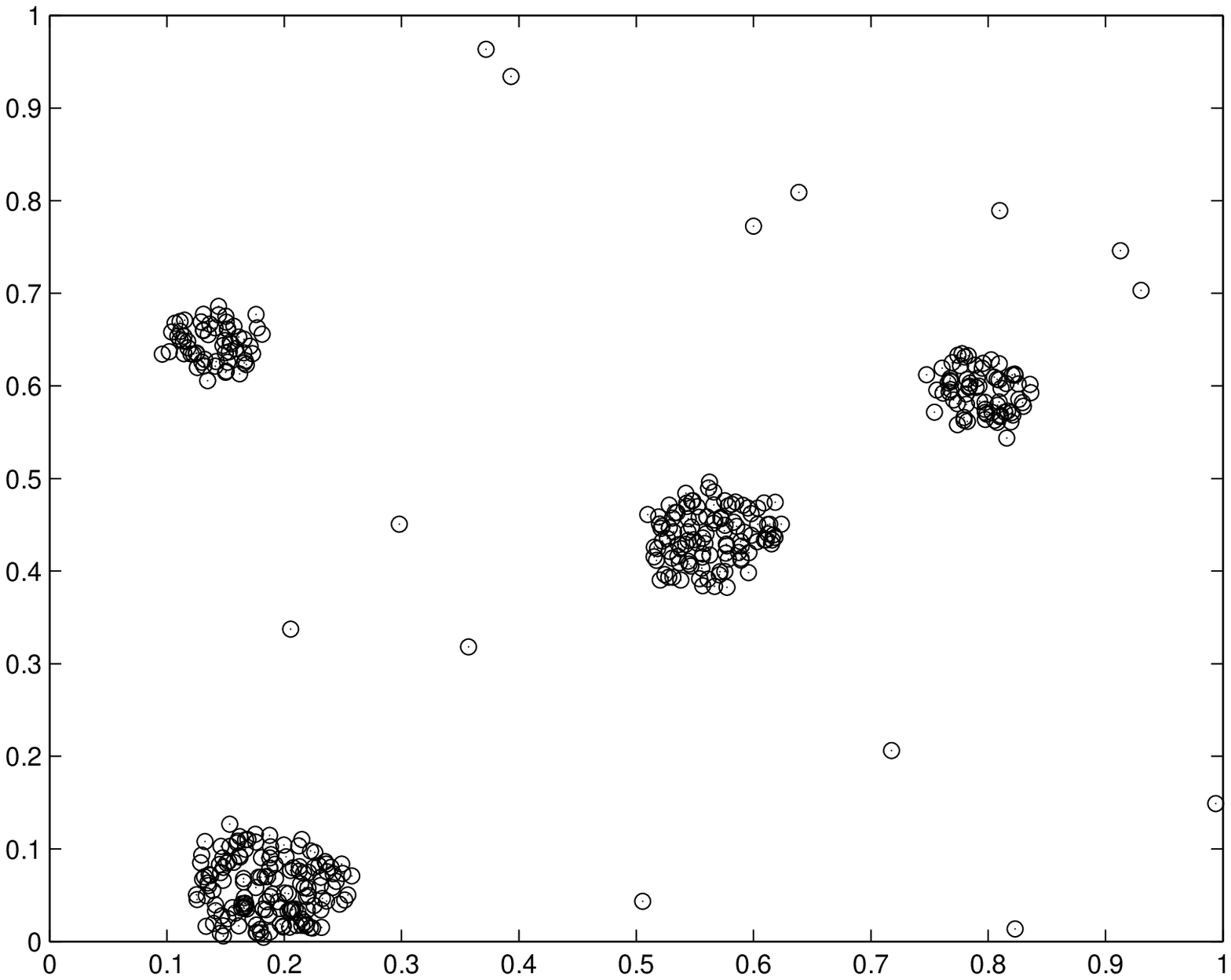}} \\
t=0.1 & t=0.25 & t=0.7
\end{tabular}\par}
\caption{The first-order individual based model with $N=400$ agents and sampling radius $R=0.05$,
subject to a random initial condition. 
\label{fig:num2}}
\end{figure}

\begin{figure}
{\centering \begin{tabular}[h]{ccc}
\resizebox*{0.3\linewidth}{!}{\includegraphics{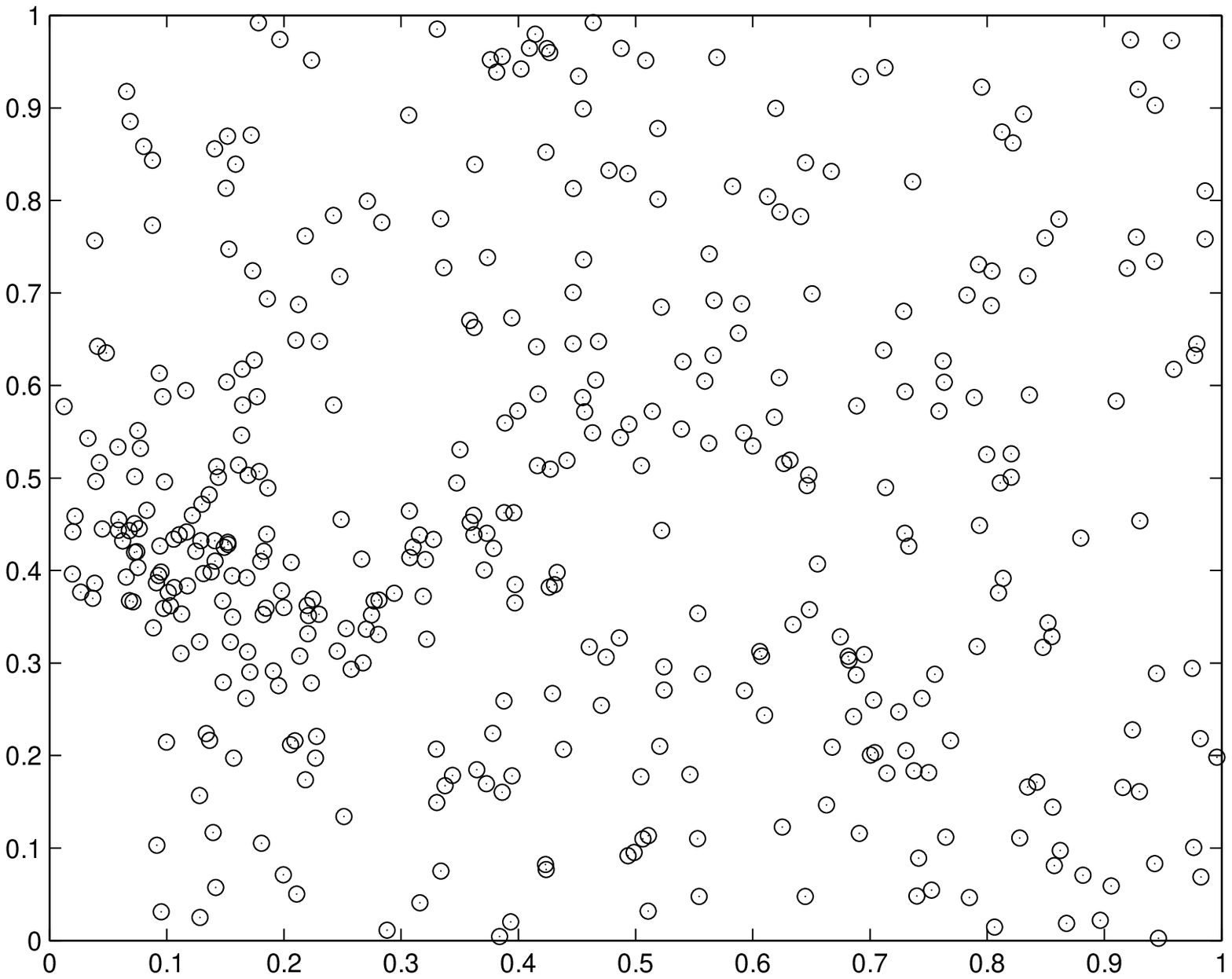}} &
\resizebox*{0.3\linewidth}{!}{\includegraphics{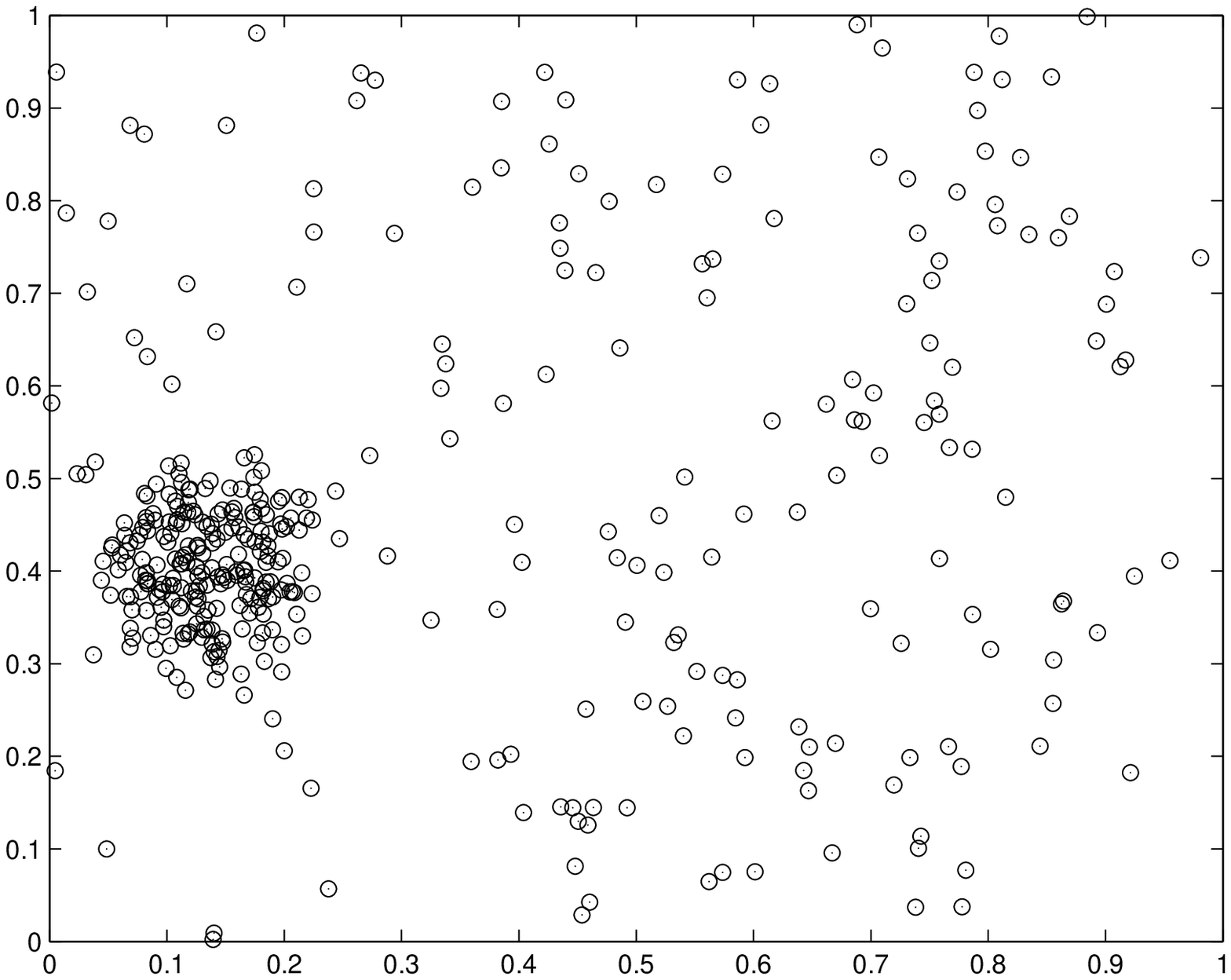}} &
\resizebox*{0.3\linewidth}{!}{\includegraphics{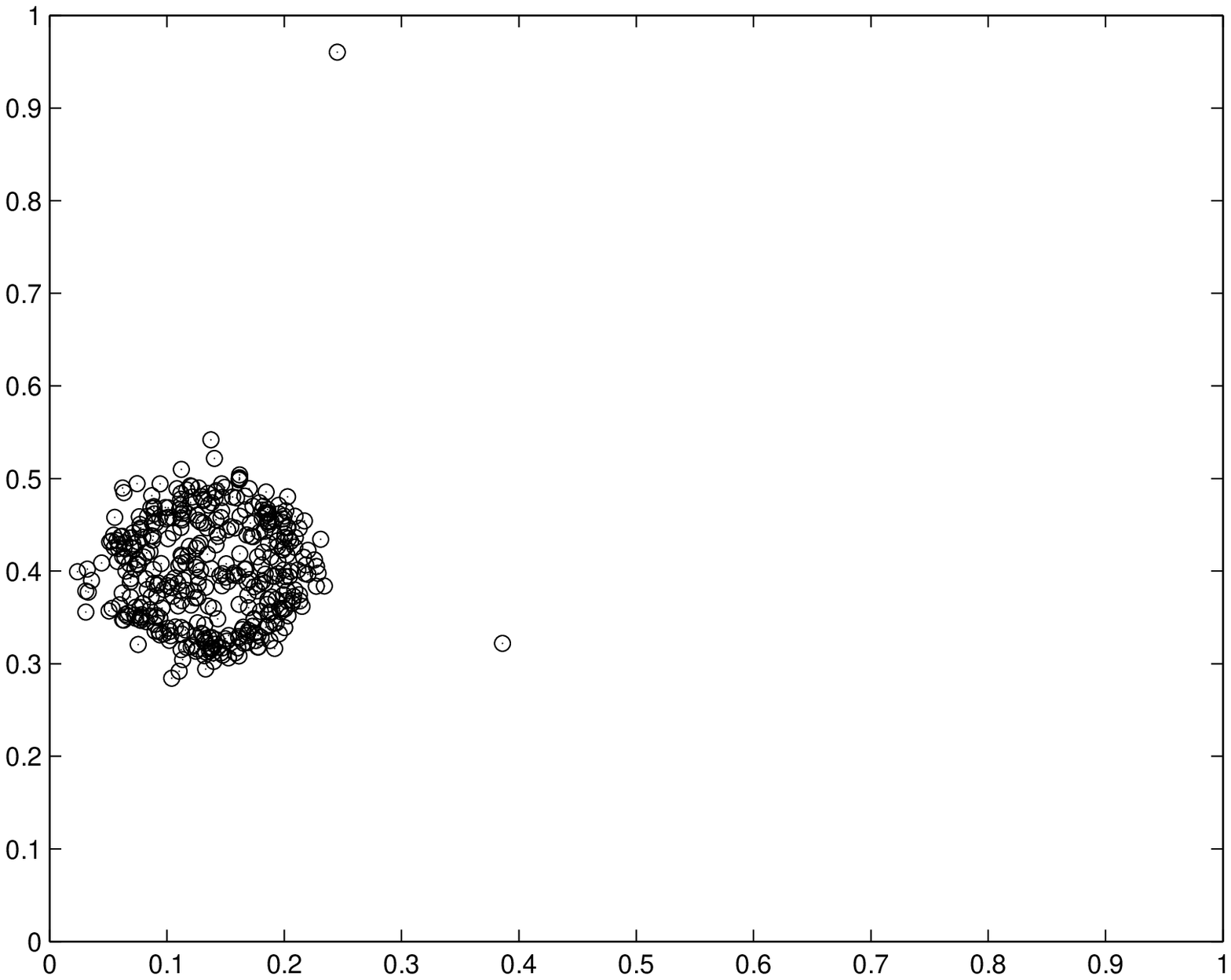}} \\
t=1.5 & t=2.0 & t=8.3
\end{tabular}\par}
\caption{The first-order individual based model with $N=400$ agents and sampling radius $R=0.1$,
subject to a random initial condition. 
\label{fig:num3}}
\end{figure}

We consider a system consisting of $N=400$ individuals
in a 2D domain $\Omega = (0,1)\times(0,1)$ with periodic boundary conditions.
The initial positions are generated randomly and independently for each individual
from the uniform distribution on $\Omega$; we took the same initial condition
for all the three experiments below. 
We choose $\W(x) = \w(|x|)$ with $\w$ the characteristic function
of the interval $[0,R]$, corresponding to the sampling radius $R$,
and for $R$ we take the values $0.025$ (Fig.~\ref{fig:num1}), $0.05$ (Fig.~\ref{fig:num2}) and $0.1$ (Fig.~\ref{fig:num3}).
For $G$ we make the choice $G(s) = \exp(-s/3)$.
The system of stochastic differential equations~\fref{model1}
is integrated in time using the Euler-Mayurama scheme
with time-step length $t=10^{-3}$.
We used the linear stability analysis in Section~\ref{sec:Numerics}
to make the ``right'' choice of $G$, such that we could observe pattern
formation. Indeed, if $G$ decreases too quickly, the system will
``freeze'' immediately, before any aggregates could be formed;
on the other hand, if $G$ does not decrease fast enough,
the system is ``overheated'' and does not allow aggregates to persist.

In Figs.~\ref{fig:num1}, \ref{fig:num2} and~\ref{fig:num3}
we observe that with a smaller sampling radius $R$,
a larger number of small aggregates is created on a faster time scale.
The choice $R=0.1$ (Fig.\ref{fig:num3}) leads to creation of one sigle
aggregate. This aggregate is approximately ring-shaped,
with higher density of particles around the circumference and lower density in the middle.
This can be explained by the fact that the aggregate grows by ``capturing'' particles
from its neighborhood, and once a particle is captured, its mobility is greatly reduced,
so that it only slowly makes its way towards the center of the aggregate.
Let us also mention that the right-most panels in Figs.~\ref{fig:num1}--\ref{fig:num3}
present quasi-steady states, where the aggregates are in a dynamic equilibrium
with a very few free-running particles. However, on a very long time scale,
the smaller aggregates typically disintegrate and their particles are caught by the larger ones.
Such a coarsening behavior is typical to diffusive aggregation systems, as for instance
the Cahn-Hilliard equation, or the nonlocal continuum model of~\cite{Bertozzi}.


\noindent
\vspace{0.3cm}

\noindent
\vspace{0.3cm}

\noindent
\vspace{0.3cm}

\subsection{The second-order individual based model~\fref{model2-1}--\fref{model2-2}}
Generally, the behavior of the second-order individual based model
is very similar to this of the first-order one, at least if we
do not impose a restricted cone of vision (i.e., $\W$ in~\fref{rho_i2} does not depend on $v$).
Indeed, with a suitable choice of parameters, one again observes the
formation of quasi-stable aggregates, whose number and size depend on the sampling radius.
The most striking difference with respect to the first-order model
is that the movement of the agents is smoother (their velocities are continuous)
and the shape of the aggregates is more complex
(we observed emergence of ellipsoidal aggregates, instead of the almost-circular ones
in the first-order model).

\begin{figure}
{\centering \begin{tabular}[h]{ccc}
\resizebox*{0.3\linewidth}{!}{\includegraphics{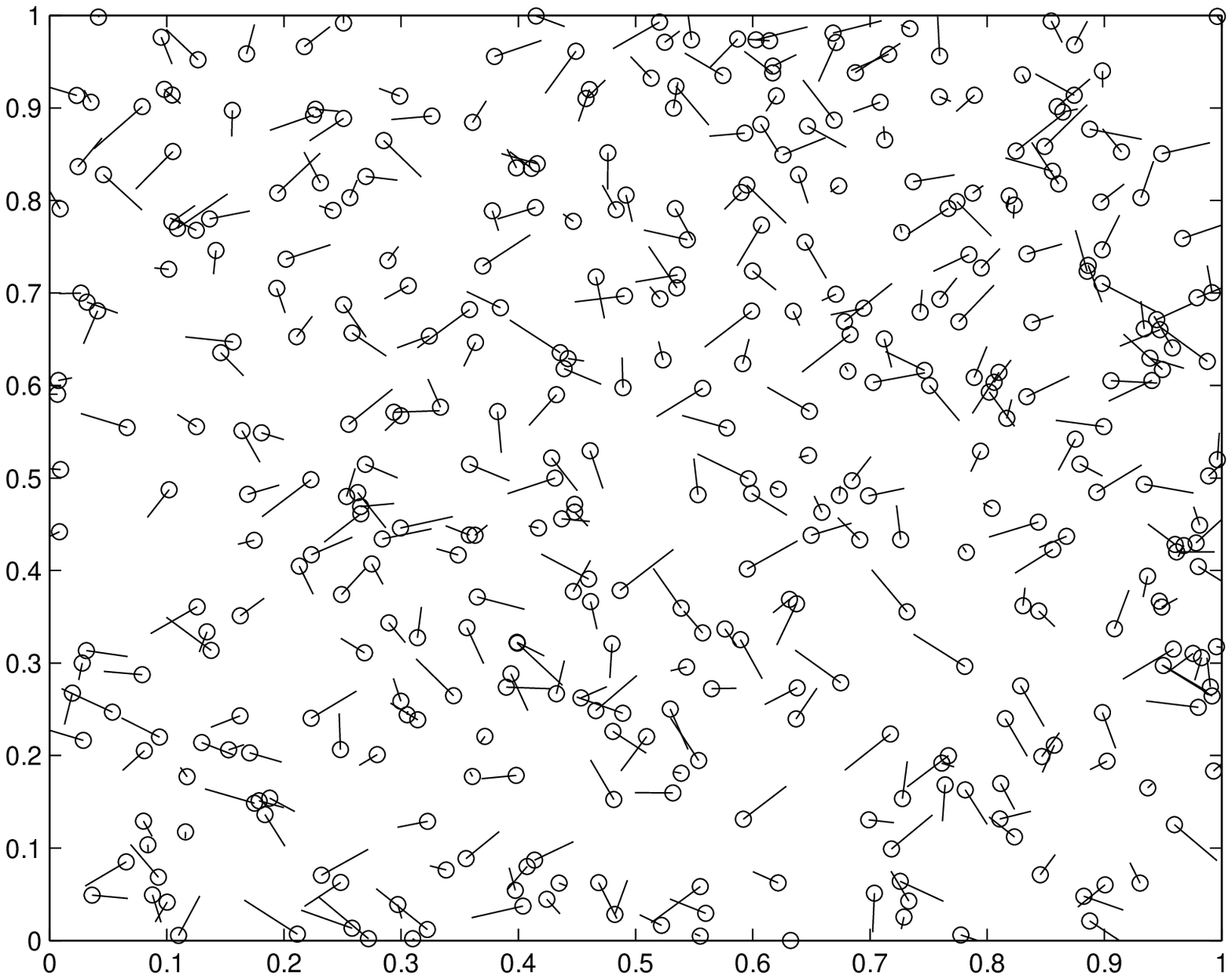}} &
\resizebox*{0.3\linewidth}{!}{\includegraphics{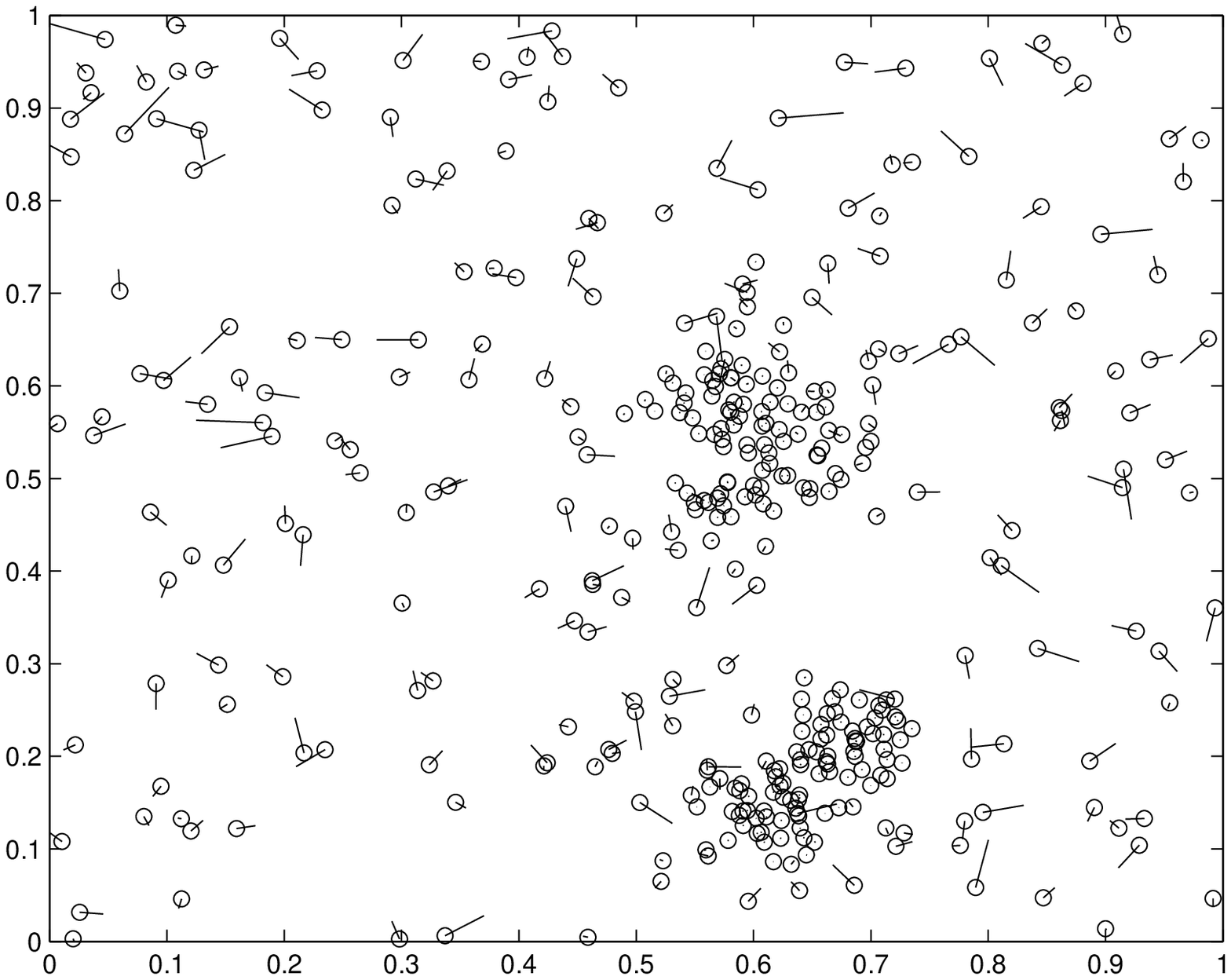}} &
\resizebox*{0.3\linewidth}{!}{\includegraphics{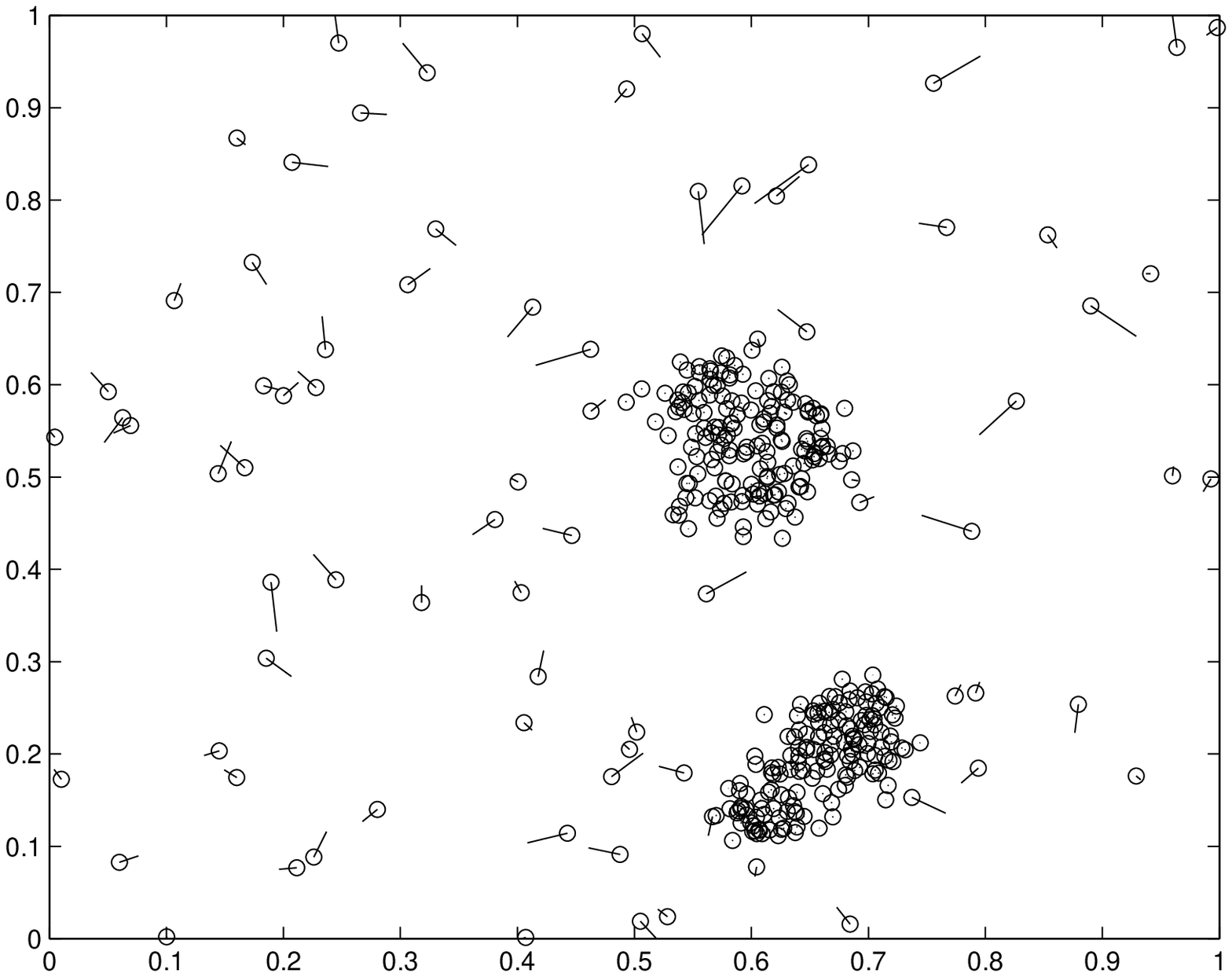}} \\
t=0.000 & t=0.044 & t=0.104
\end{tabular}\par}
\caption{The second-order individual based model with $N=400$ agents, sampling radius $R=0.05$
and 180\textdegree\ cone of vision, subject to a random initial condition.
\label{fig:num_im2}}
\end{figure}

The situation becomes slightly more interesting if we consider a restricted cone of vision
- we choose 180\textdegree, such that $\W(x,v) = w\left(|x|,\frac{x\cdot v}{|x| |v|}\right)$
with $w(s,z) = \chi_{[0,R]}(s)\chi_{[0,1]}(z)$. The sampling radius is chosen as $R=0.05$.
Again, we simulate with $N=400$ individuals
in a 2D domain $\Omega = (0,1)\times(0,1)$ with periodic boundary conditions,
and $G(s) = \exp(-s/3)$. For $H$ we set the constant $H\equiv 2$.
The system of stochastic differential equations~\fref{model2-1}--\fref{model2-2}
is integrated in time using the Euler-Mayurama scheme
with time-step length $t=10^{-3}$.
The initial positions of the agents are generated randomly
in the same way as for the first-order model,
while their initial velocities are generated independently and randomly
from the 2D normalized Gaussian distribution.
Three snapshots of the evolution of the system are shown in Fig.~\ref{fig:num_im2},
the velocities being marked by the linear elements for each individual.
The right-most panel in Fig.~\ref{fig:num_im2} shows the quasi-steady state
with two aggregates formed, in dynamic equilibrium with the ``free-running'' individuals.
It is obvious that, compared to the previous simulations, the aggregates are less densely packed
and their shapes are less circular. Also, the portion of the ``free-running'' particles is much higher,
which clearly is an effect of the restricted cone of vision
(the captured particles can leave the aggregate more easily).

\noindent
\vspace{0.3cm}

\subsection{The first-order mean-field model~\fref{mmodel1} in 1D}
\begin{figure}
{\centering \begin{tabular}[h]{ccc}
\resizebox*{0.3\linewidth}{!}{\includegraphics{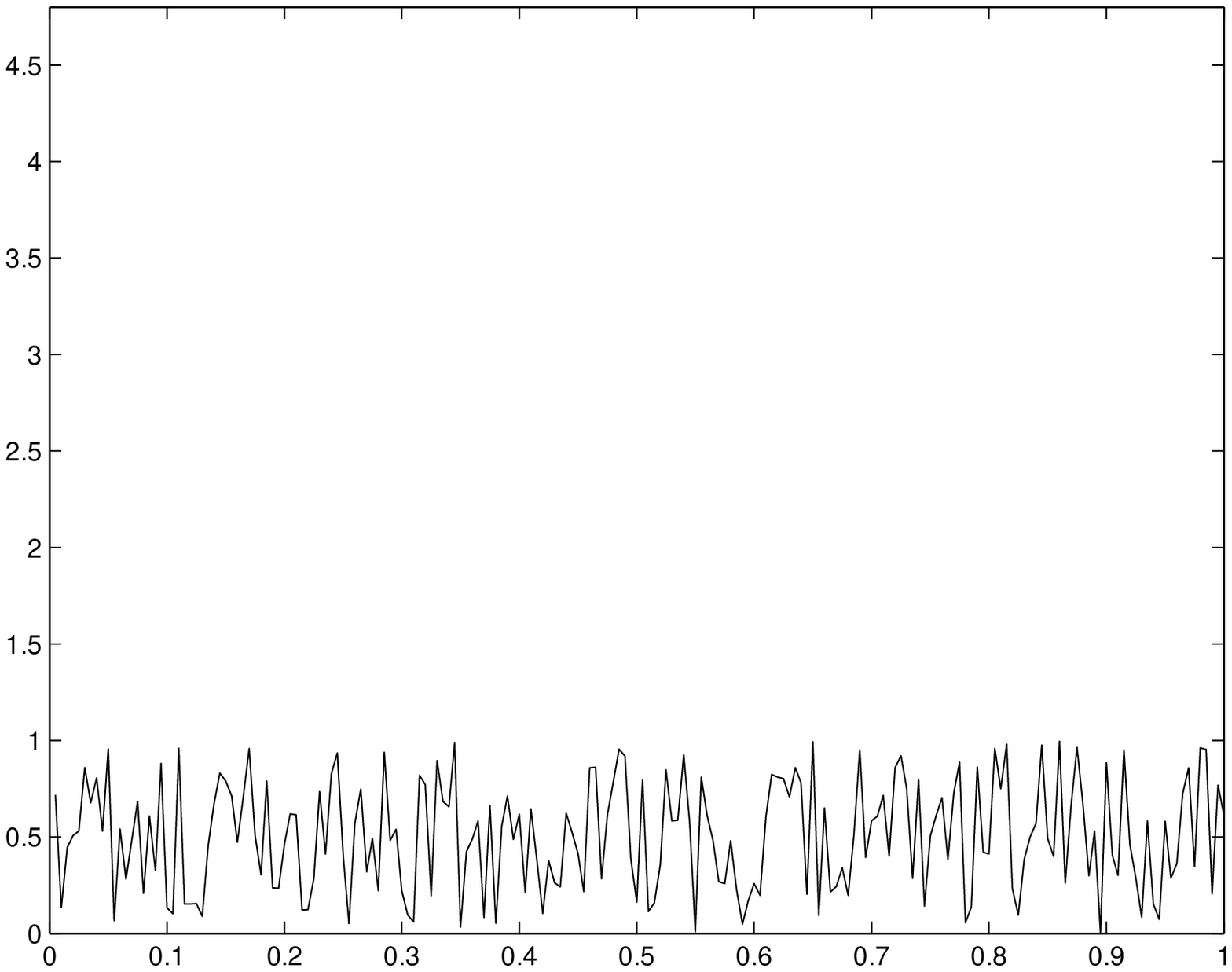}} &
\resizebox*{0.3\linewidth}{!}{\includegraphics{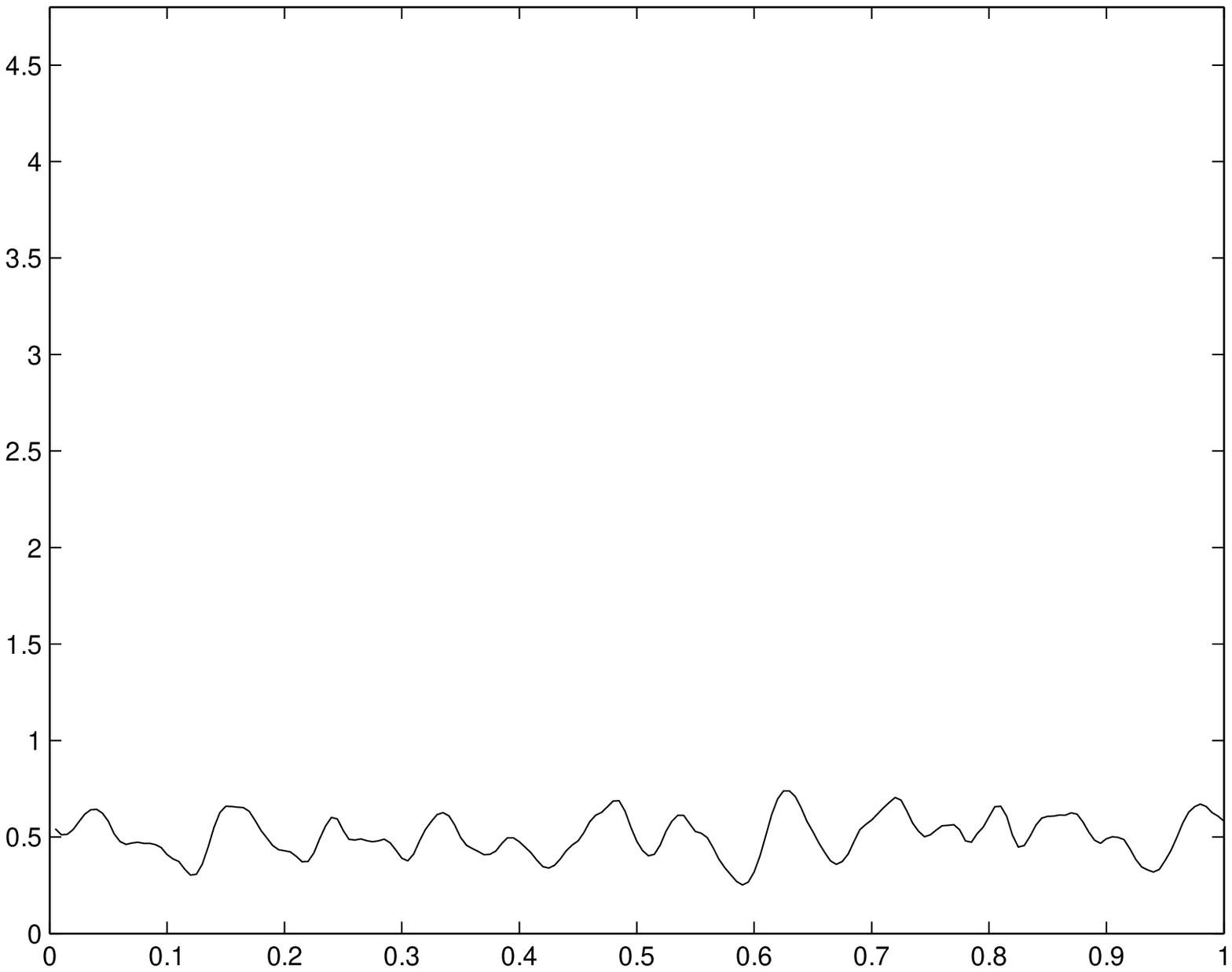}} &
\resizebox*{0.3\linewidth}{!}{\includegraphics{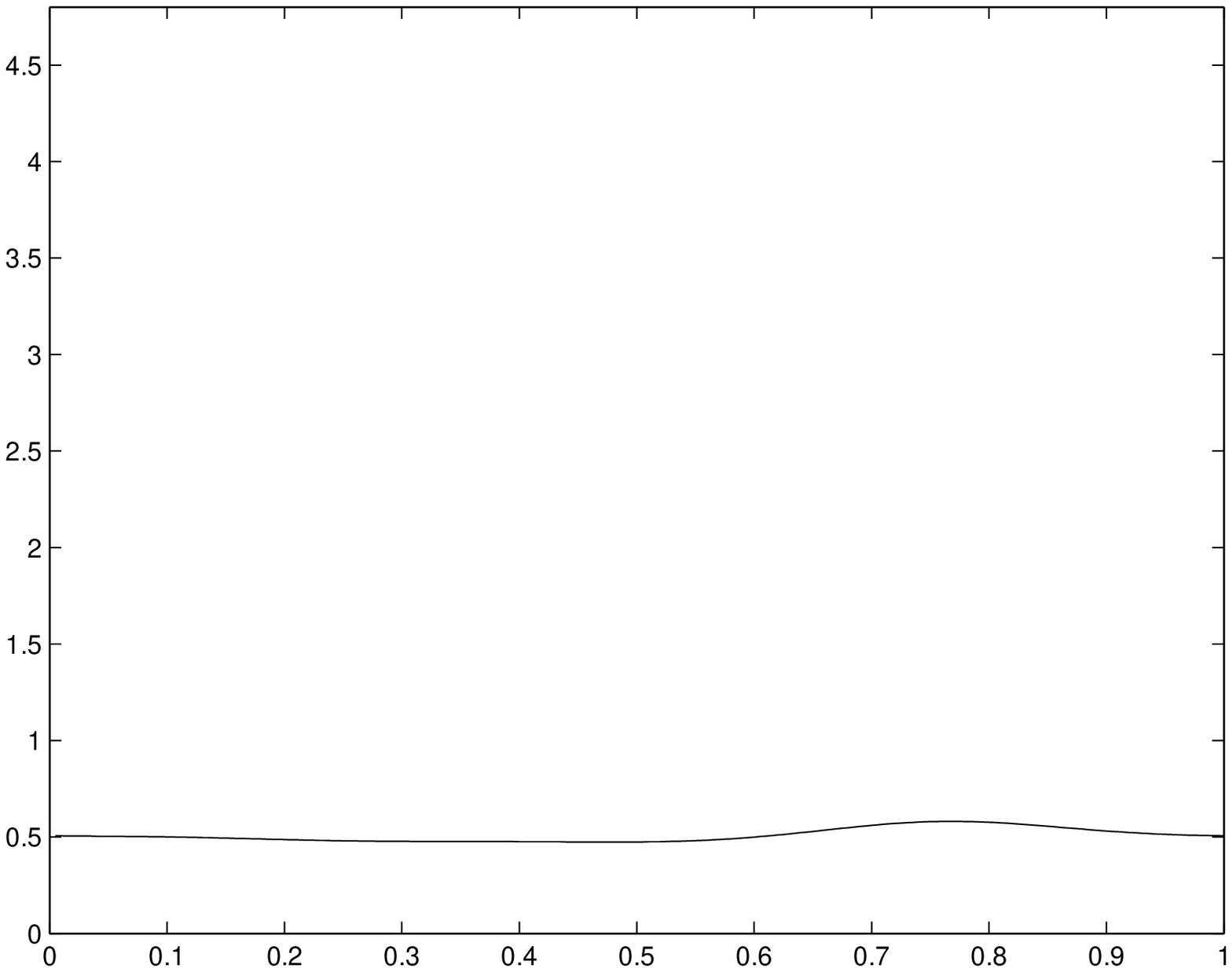}} \\
t=0.0000 & t=0.0008 & t=0.1300 \\
\resizebox*{0.3\linewidth}{!}{\includegraphics{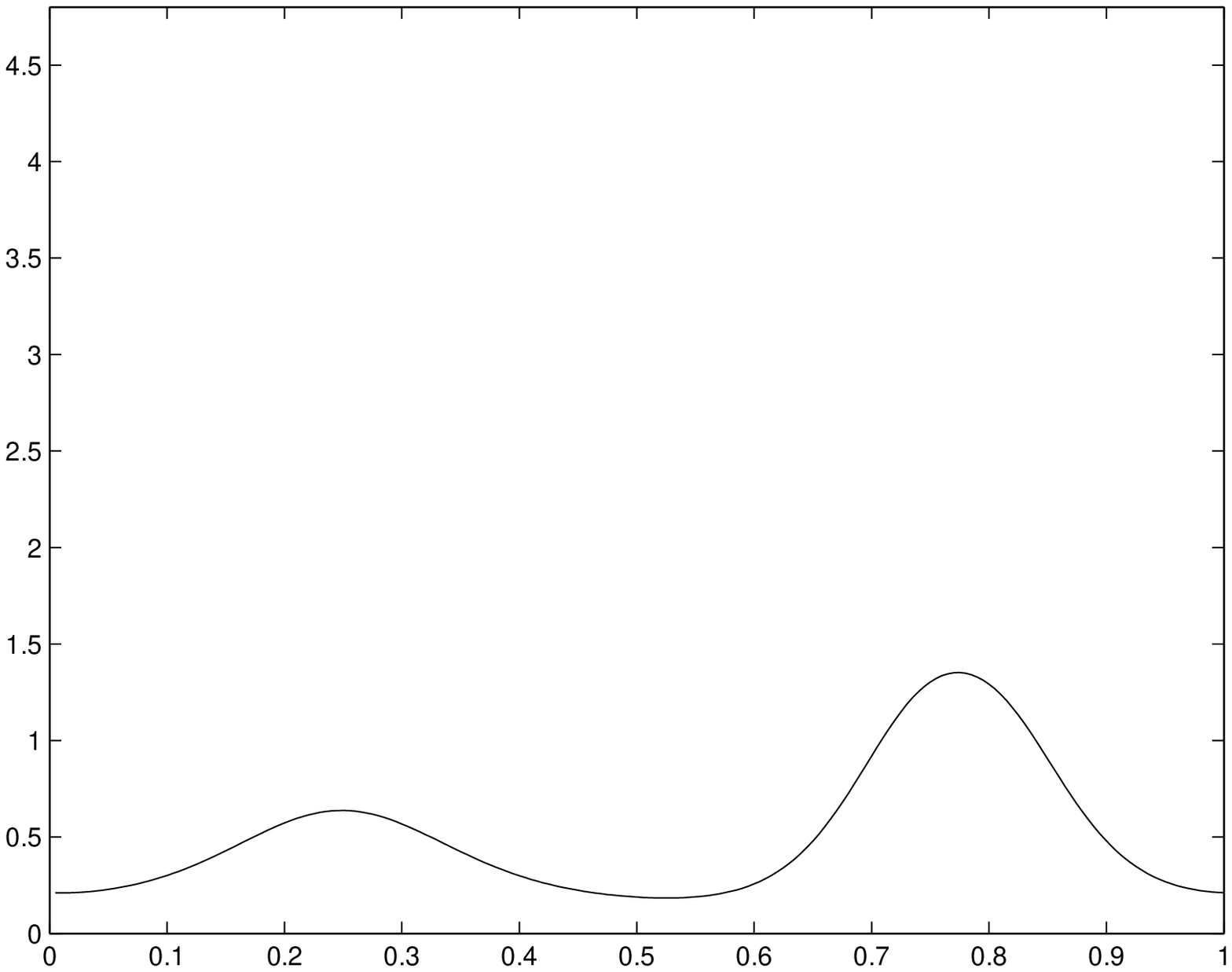}} &
\resizebox*{0.3\linewidth}{!}{\includegraphics{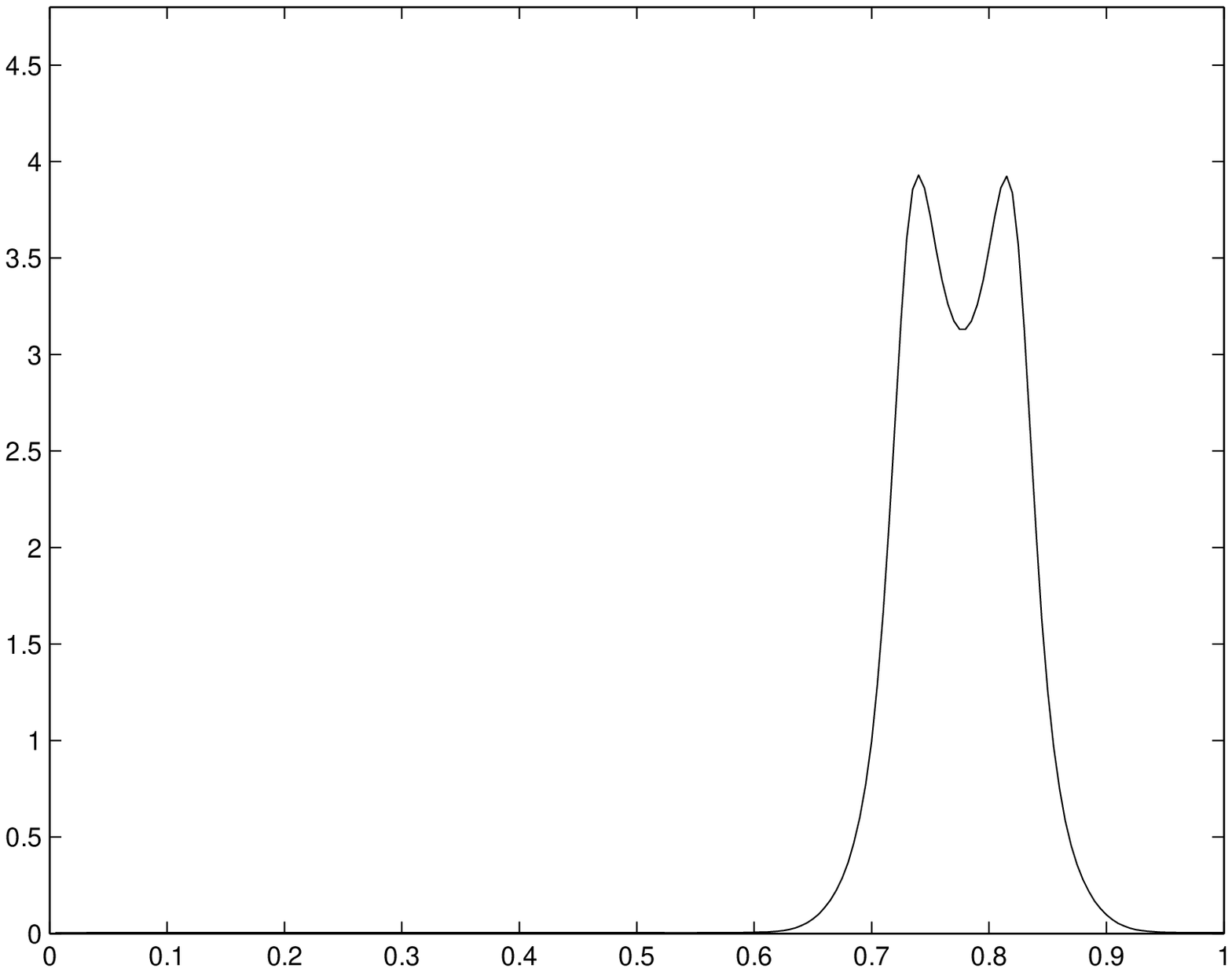}} &
\resizebox*{0.3\linewidth}{!}{\includegraphics{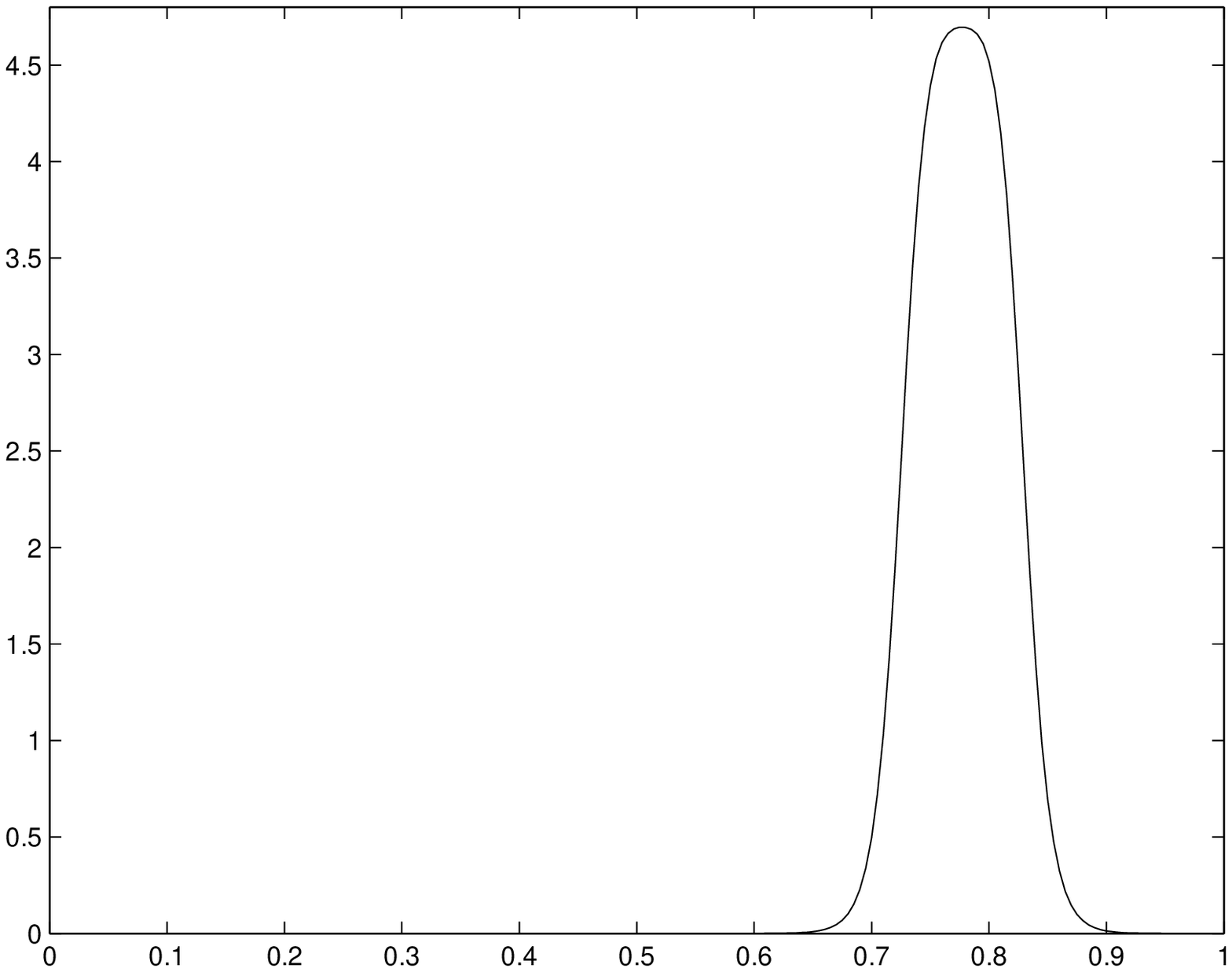}} \\
t=0.8500 & t=2.6500 & t=11.0500
\end{tabular}\par}
\caption{The first-order mean-field model in a periodic 1D setting
with random initial condition (upper left panel).
Steady state was reached at $t\simeq 11$ (lower right panel).
\label{fig:num_mm1}}
\end{figure}

We simulated the first-order mean-field model~\fref{mmodel1}
in the 1D periodic domain $(0,1)$, using semi-implicit
finite difference discretization for the space variable
and first-order forward Euler method for the time variable.
The space grid consisted of $200$ equidistant points,
the time step was $10^{-4}$.
As before, we chose $\W(x) = \w(|x|)$ with $\w$ the characteristic function
of the interval $[0,0.1]$, and $G(s) = \exp(-s/3)$.
We imposed a random initial condition for $\varrho$,
generated such that for every grid point a random number
from the uniform distribution in $[0,1]$ has been drawn.
Snapshots of the evolution are shown in Fig.~\ref{fig:num_mm1}.
We observe that quite a strong smoothing effect takes place on the fast time-scale
(first row in Fig.~\ref{fig:num_mm1}),
while aggregation takes place on a time-scale approximately one order of magnitude slower (second row);
this is explained with the stability analysis in Section~\ref{sec:Stability}.
First, two aggregates of different sizes are created, however, both of them are unstable -
the smaller one is smoothed out, while the larger grows further, until the steady state is reached
(lower right panel). Observe also the characteristic ``fork''-like shape of the profile
in the lower mid panel ($t=2.65$), which is due to the mass arriving from the neighbourhood
with higher diffusivity than the diffusivity in the middle of the profile;
compare also with the ring-shaped aggregate in the right panel of Fig.~\ref{fig:num2}.
However, this fork-like structure is eventually also smoothed out,
to finally obtain the steady profile. Let us note that the steady aggregate, although well localized,
is not compactly supported, i.e., the profile has positive density everywhere on $[0,1]$.

\noindent
\vspace{0.3cm}

\subsection{The first-order mean-field model~\fref{mmodel1} in 2D}
\begin{figure}
{\centering \begin{tabular}[h]{cc}
\resizebox*{0.49\linewidth}{!}{\includegraphics{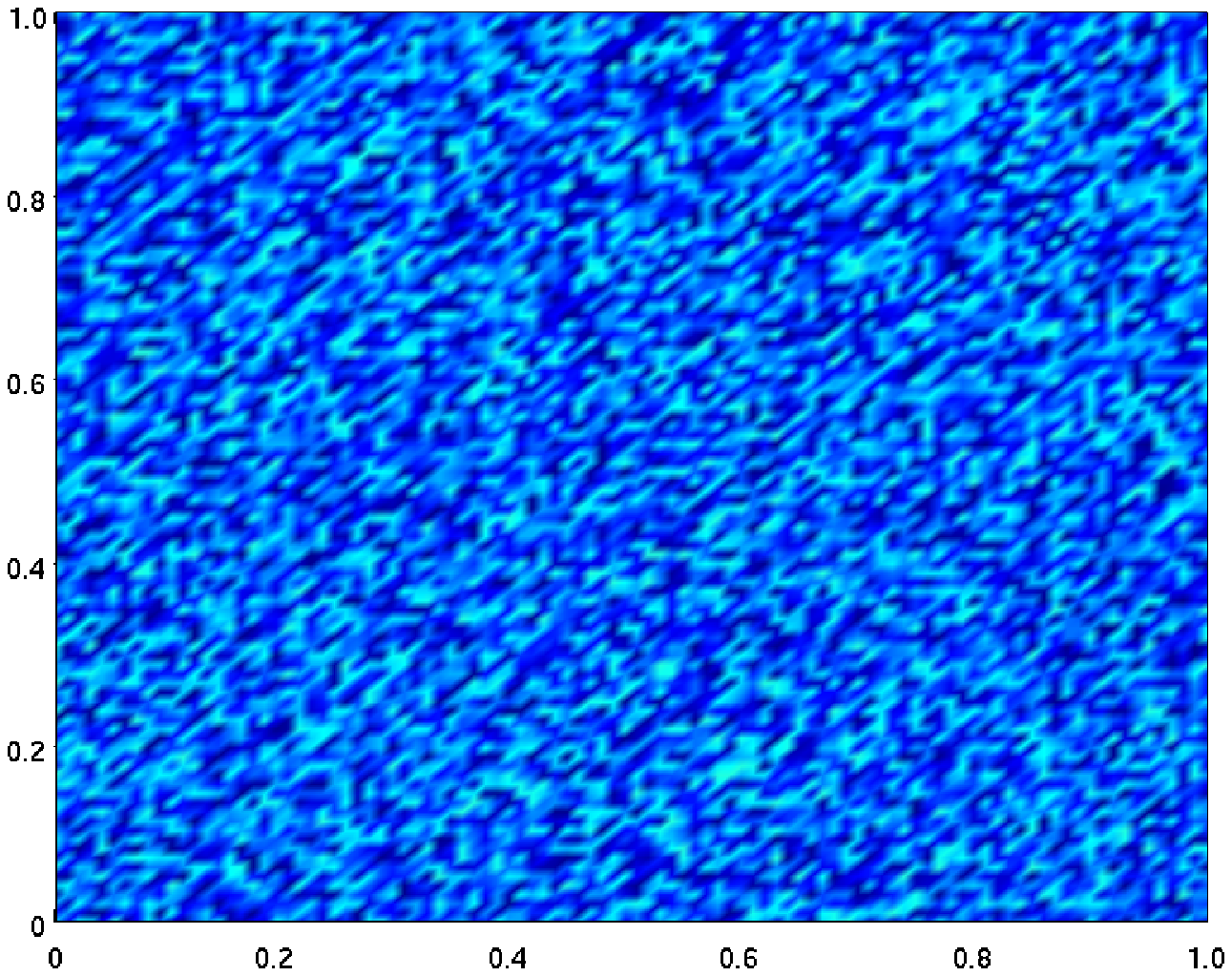}} &
\resizebox*{0.49\linewidth}{!}{\includegraphics{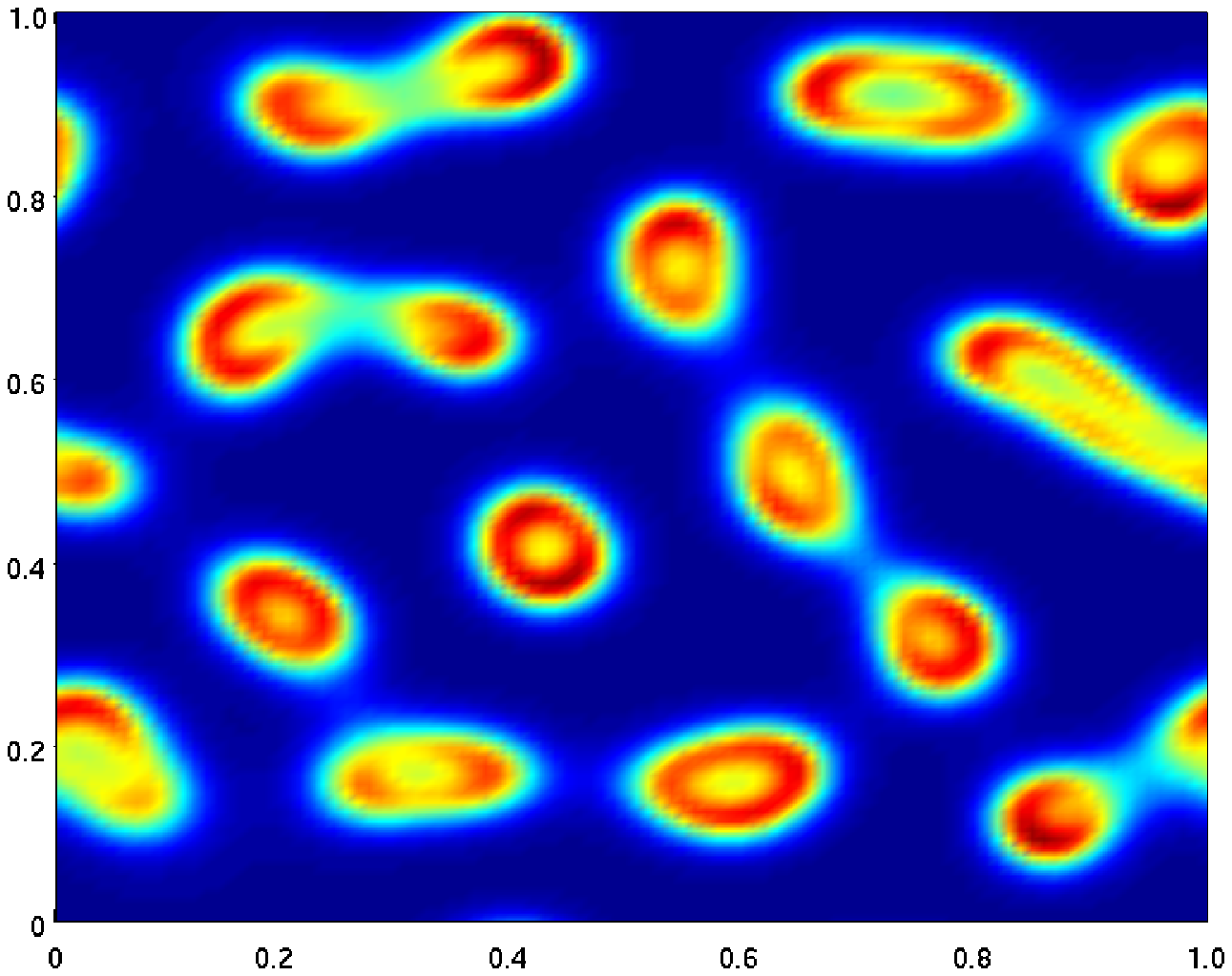}} \\
t=0.000 & t=0.005 \\
\resizebox*{0.49\linewidth}{!}{\includegraphics{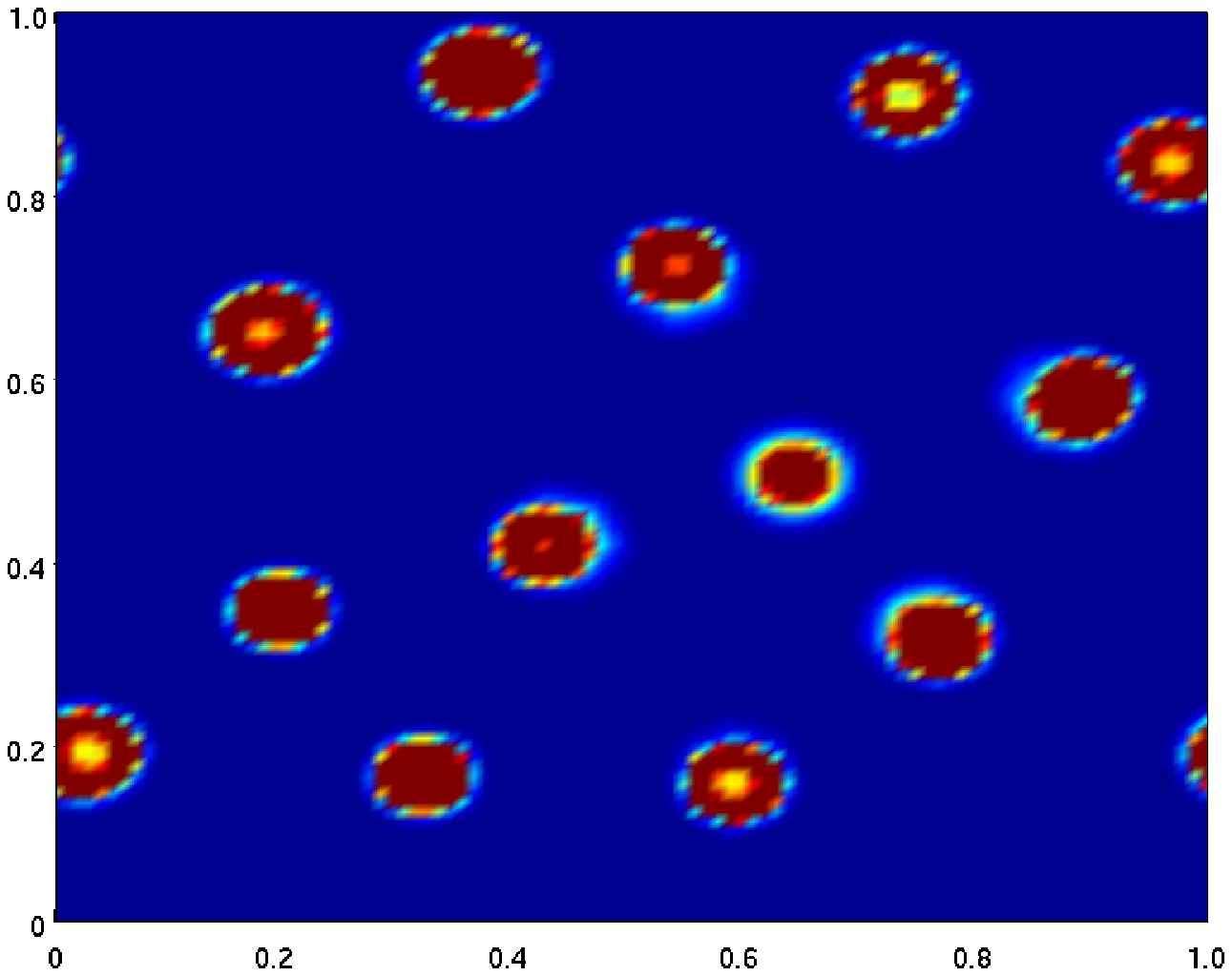}} &
\resizebox*{0.49\linewidth}{!}{\includegraphics{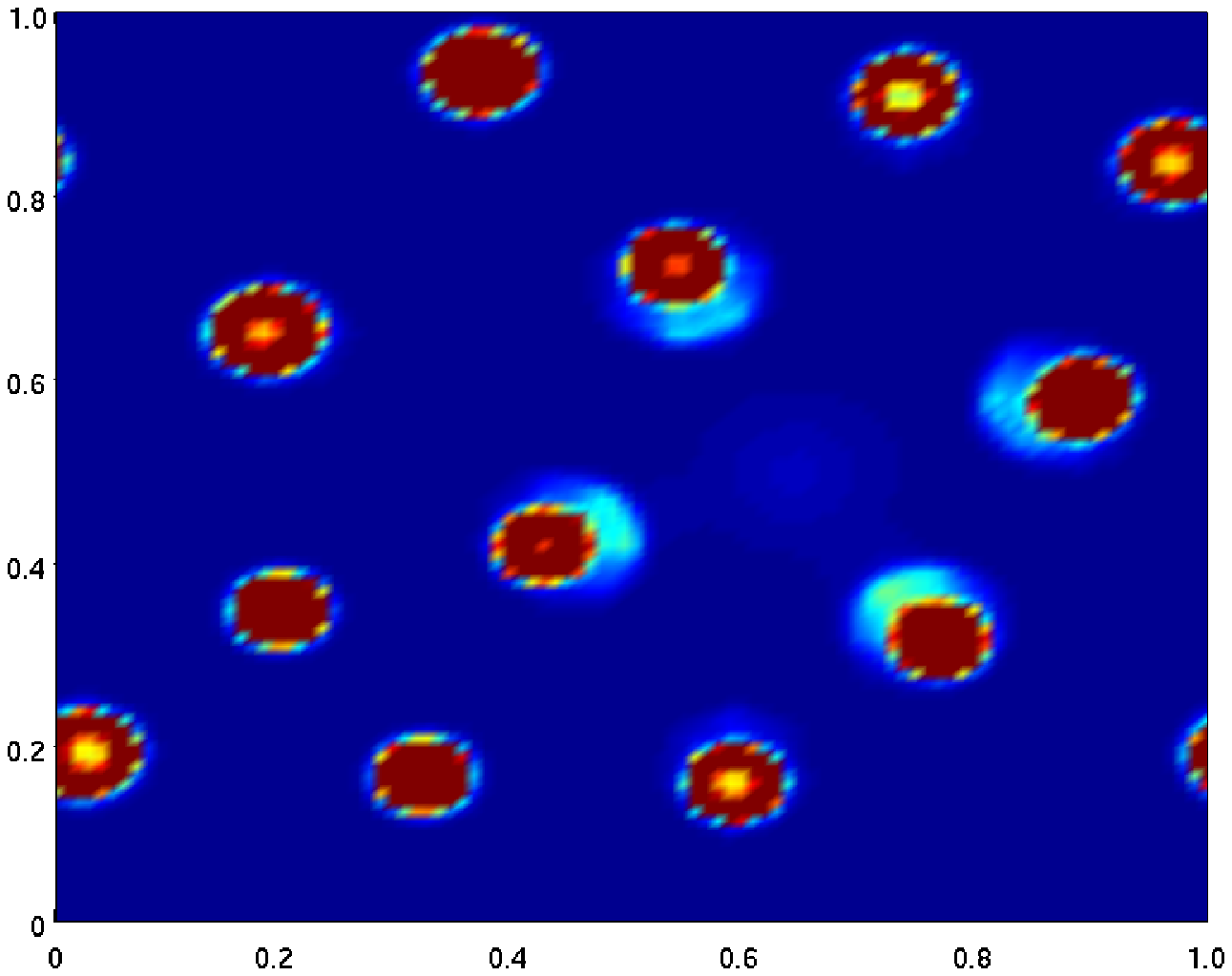}} \\
t=0.900 & t=0.950
\end{tabular}\par}
\caption{The first-order mean-field model in a periodic 2D setting
with random initial condition (upper left panel).
After the initial rapid smoothing of the high-frequency components,
several ring-shaped structures are created (upper right panel),
which eventually turn into an almost regular pattern of well localized aggregates
(lower left panel). However, the smaller aggregates may be unstable
and diffusively disintegrate, and their mass is absorbed by their neighbors,
(lower right panel).
\label{fig:num_mm2}}
\end{figure}

\begin{figure}
{\centering \begin{tabular}[h]{cc}
\resizebox*{0.49\linewidth}{!}{\includegraphics{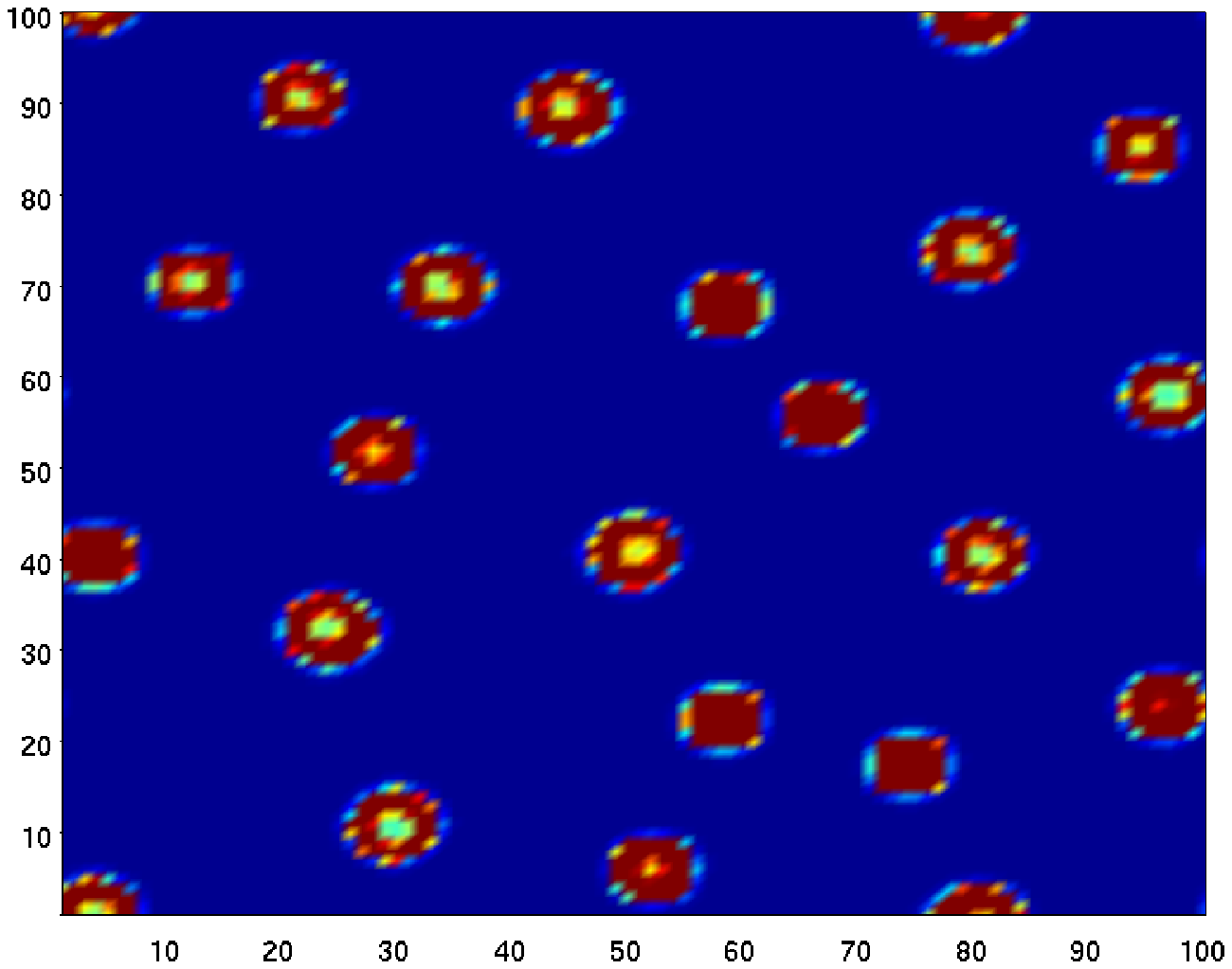}} &
\resizebox*{0.49\linewidth}{!}{\includegraphics{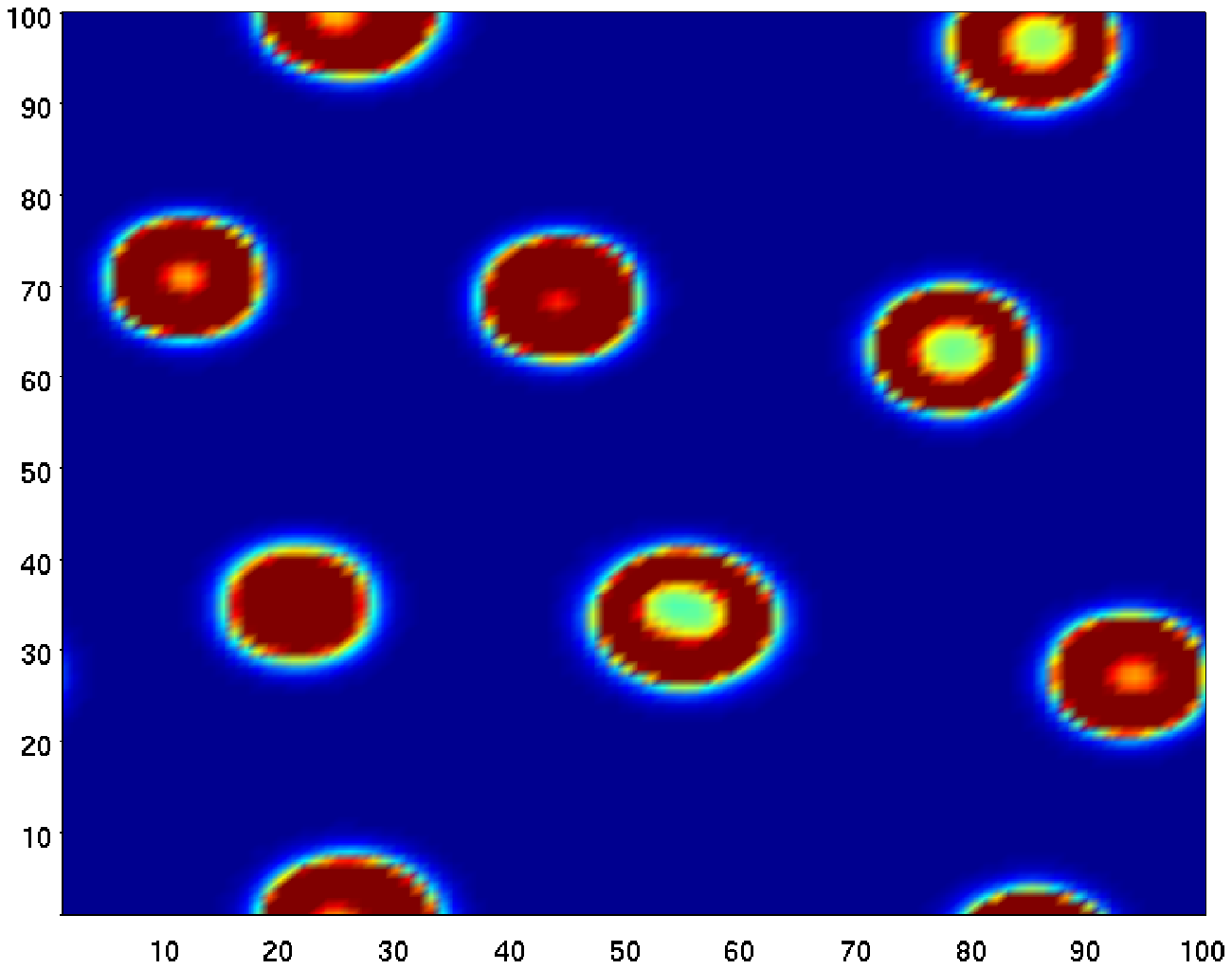}} \\
R=0.06 & R=0.11
\end{tabular}\par}
\caption{An example of patterns produced by the first-order mean-field model in a periodic 2D setting
with the sampling radii $R=0.6$ (left panel) and $R=0.11$ (right panel). 
\label{fig:num_mm3}}
\end{figure}

We simulated the first-order mean-field model~\fref{mmodel1}
in the 2D periodic domain $(0,1)^2$, using the same type
of discretization as in the 1D case.
The space grid consisted of $100\times 100$ equidistant points,
the time step was $10^{-4}$.
We chose the sampling radius $R=0.07$, i.e., $\W(x) = \chi_{[0,0.7]}(|x|)$,
and $G(s) = \exp(-s/3)$.
Snapshots of the evolution are shown in Fig.~\ref{fig:num_mm2}.
Starting again from a random initial condition,
we observed a rapid formation of approximately ring-shaped pre-aggregates,
which eventually turn into an almost regular pattern of well localized (but not compactly supported) clumps.
However, the smaller aggregates may be unstable
and diffusively disintegrate, and their mass is absorbed by their neighbors,
as shown on the lower right panel of Fig.~\ref{fig:num_mm2}.
Finally, in Fig.~\ref{fig:num_mm3} we present examples of patterns
produced with the sampling radii $R=0.6$ (left panel) and $R=0.11$ (right panel).

\noindent

\noindent

\subsection{The second-order mean-field model~\fref{mmodel2}}
\begin{figure}
\begin{center}
\subfigure[Particle distribution function $f(t,x,v)$ at $t=4$]{\includegraphics[width = 0.45 \textwidth]{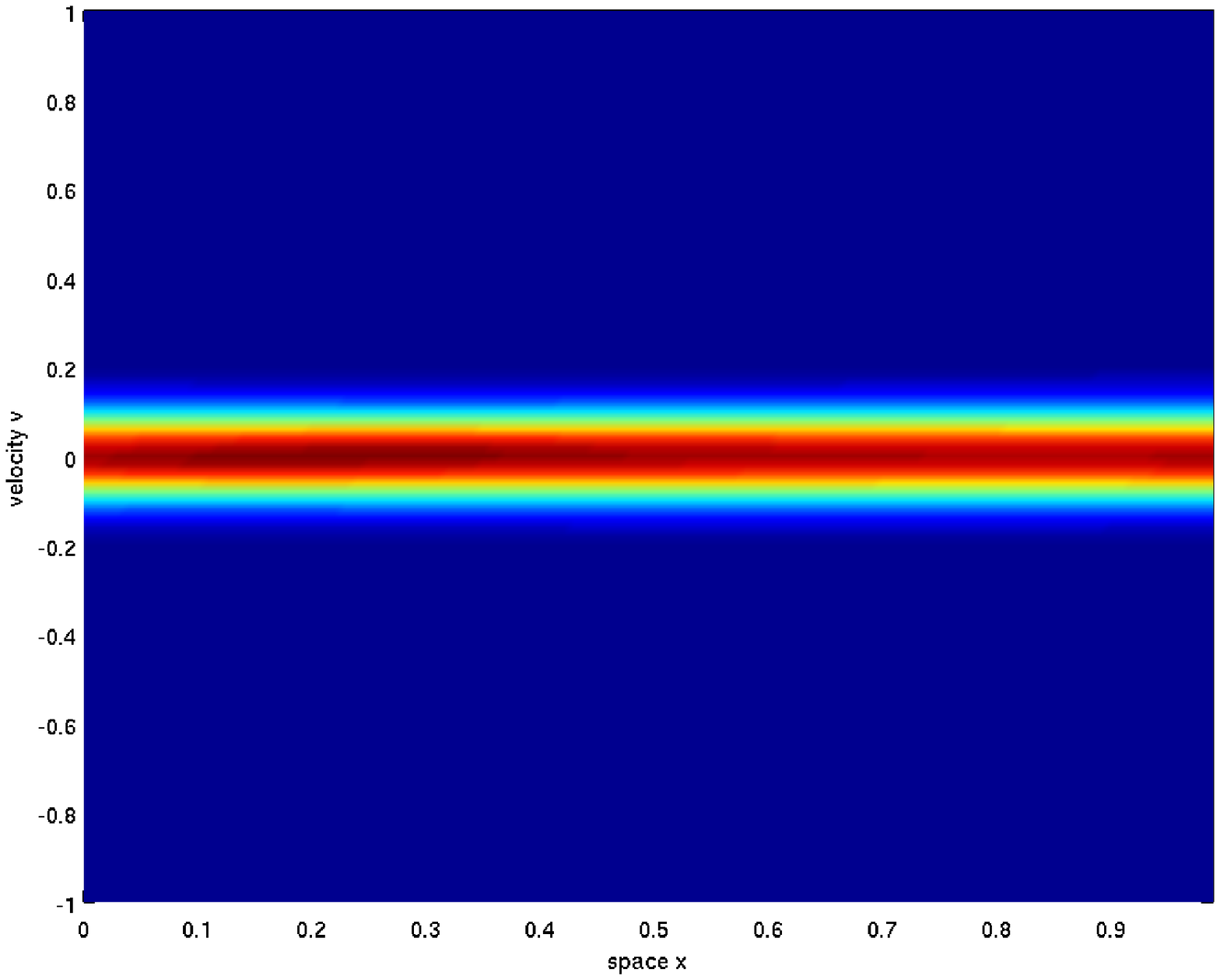}}\hspace*{1cm}
\subfigure[Mass density $\rho(t,x)$ at $t=4$]{\includegraphics[width = 0.45 \textwidth]{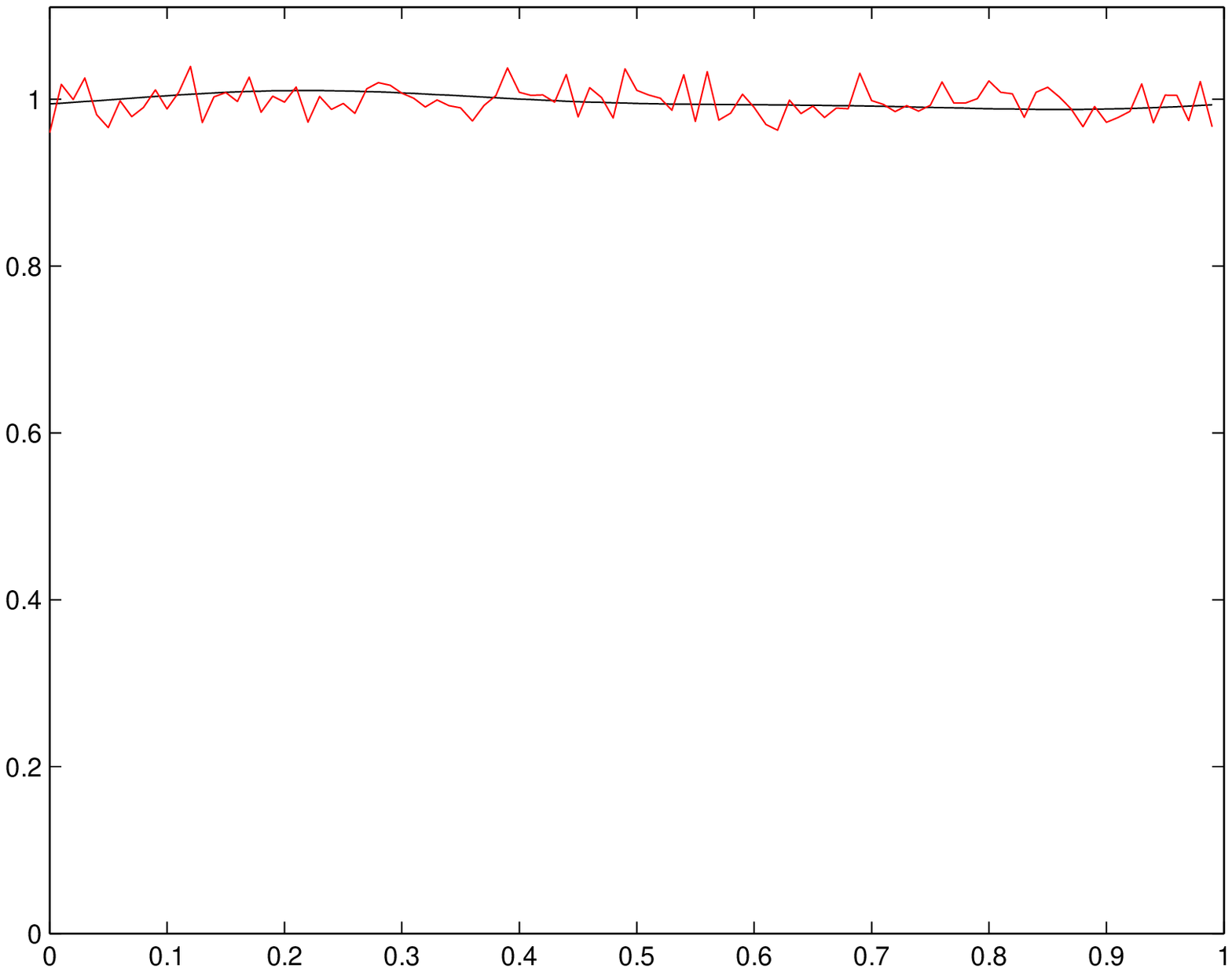}}\\
\subfigure[Particle distribution function $f(t,x,v)$ at $t=8$]{\includegraphics[width = 0.45 \textwidth]{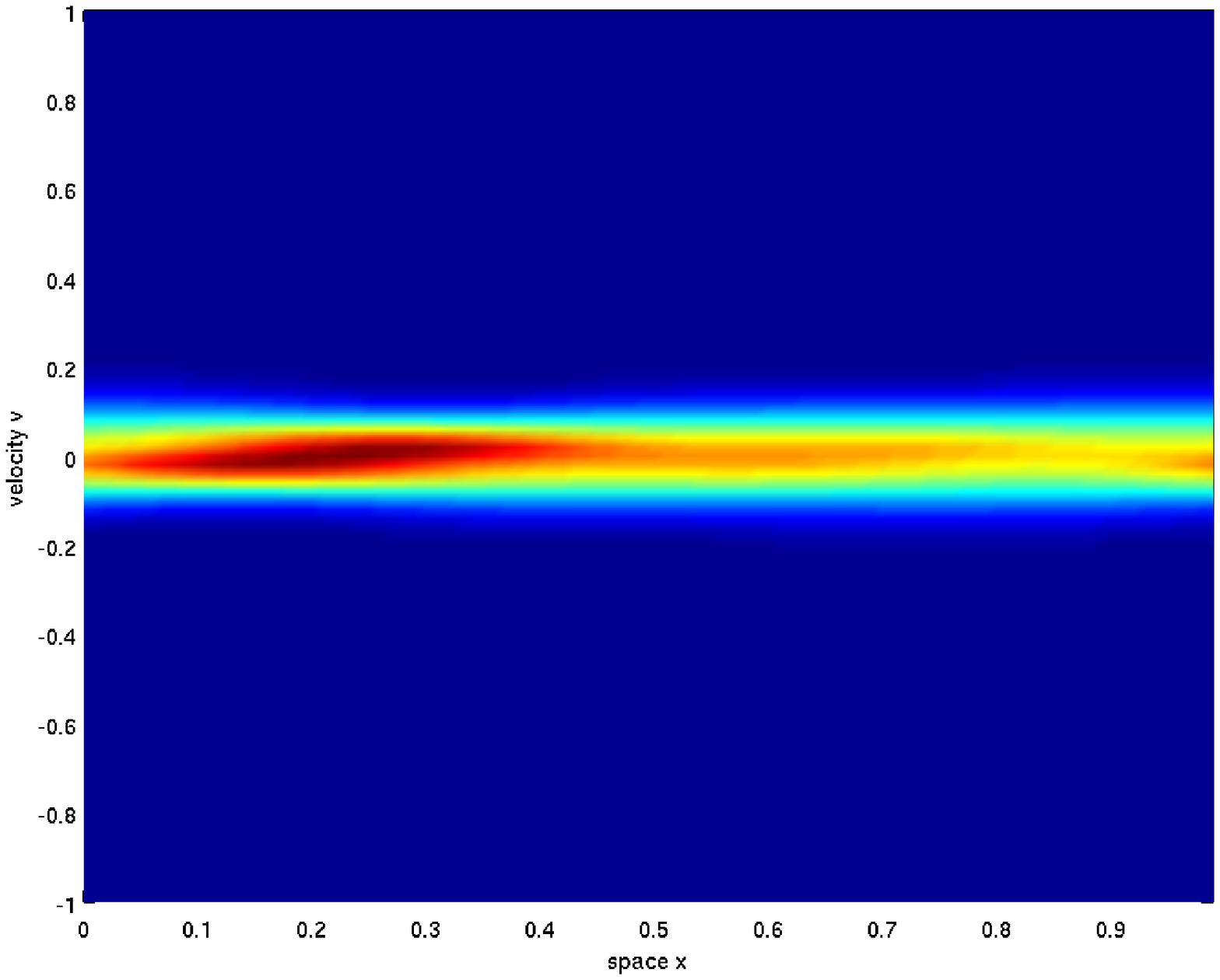}}\hspace*{1cm}
\subfigure[Mass density $\rho(t,x)$ at $t=8$]{\includegraphics[width = 0.45 \textwidth]{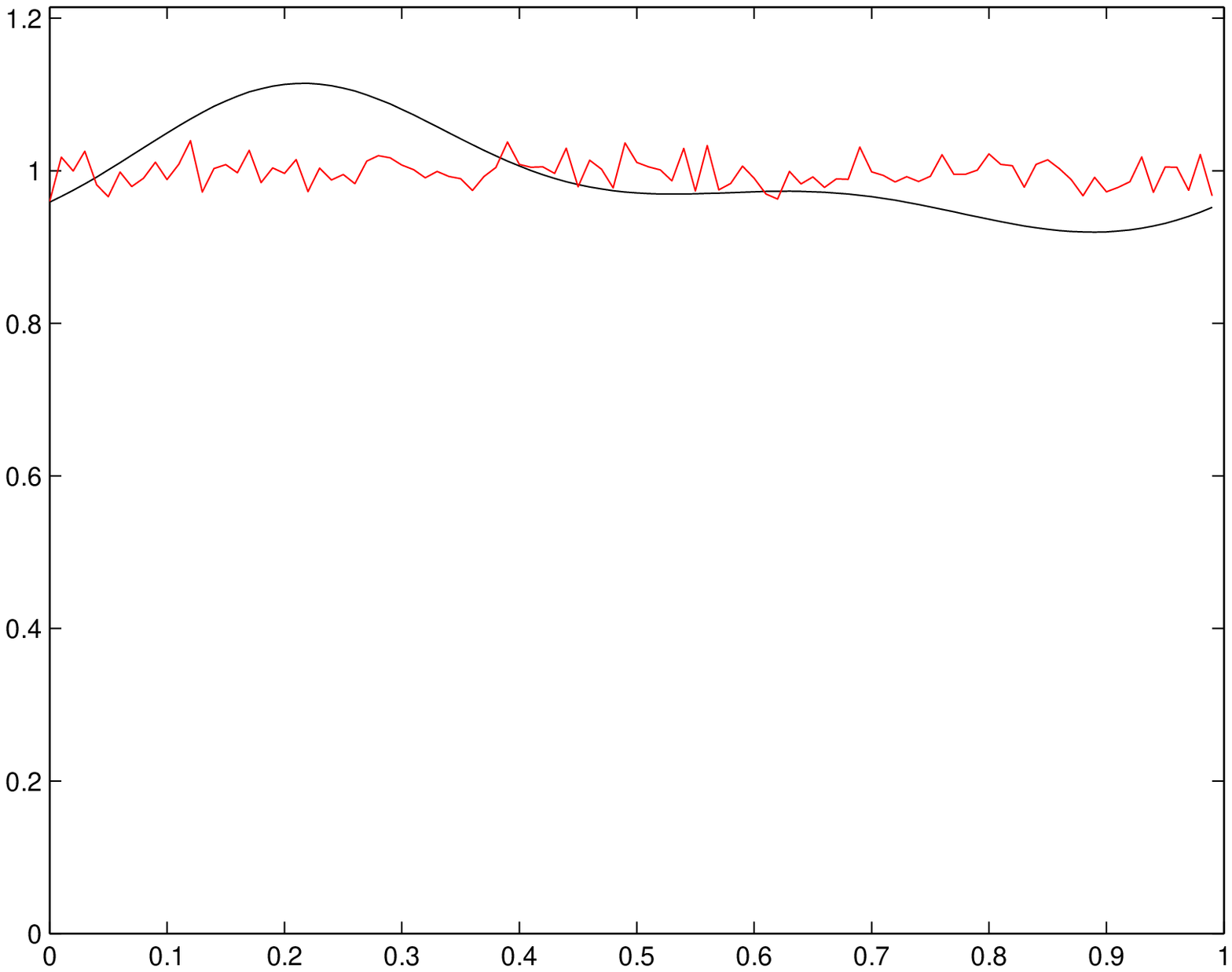}}\\
\subfigure[Particle distribution function $f(t,x,v)$ at $t=20$]{\includegraphics[width = 0.45 \textwidth]{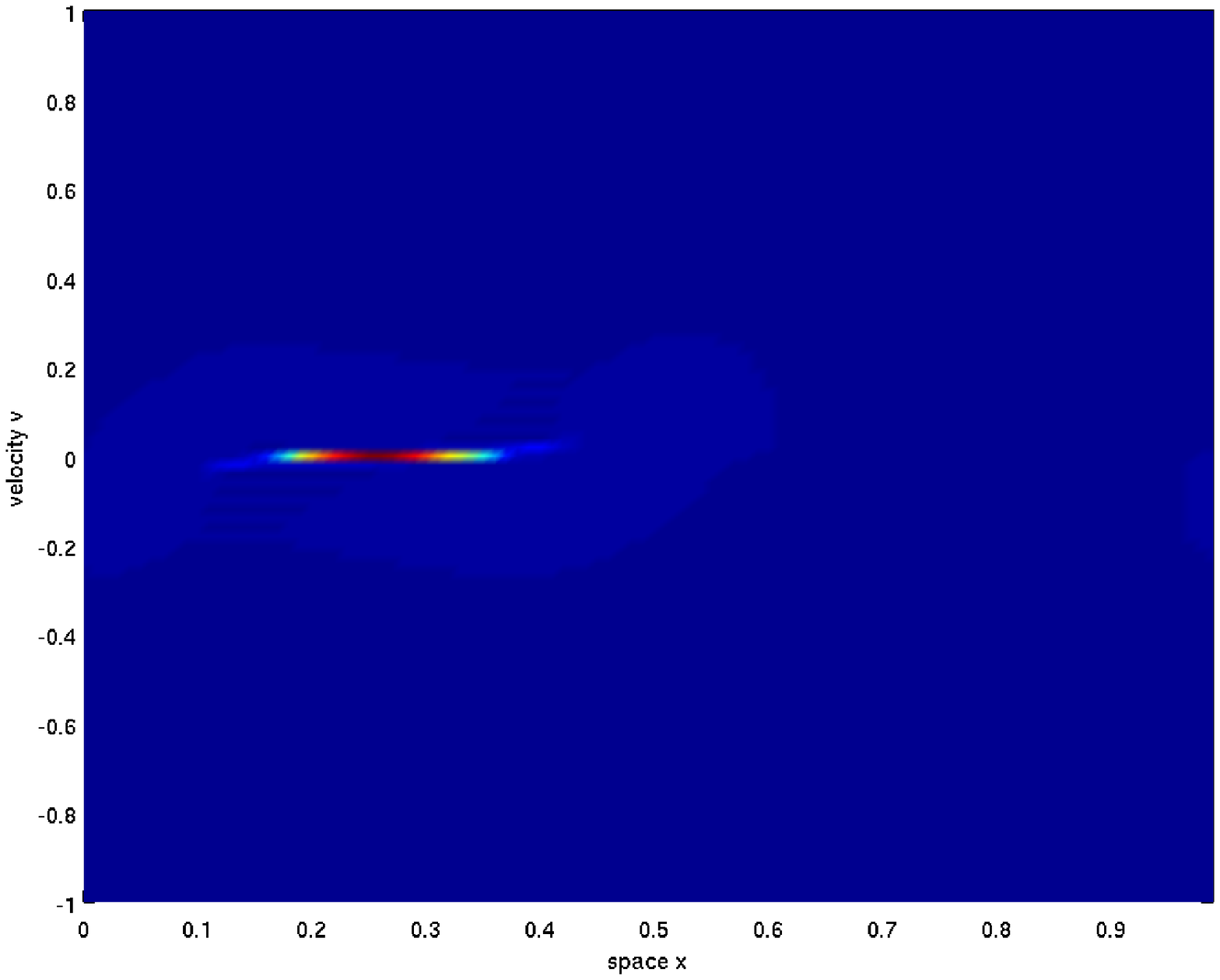}}\hspace*{1cm}
\subfigure[Mass density $\rho(t,x)$ at $t=20$]{\includegraphics[width = 0.45 \textwidth]{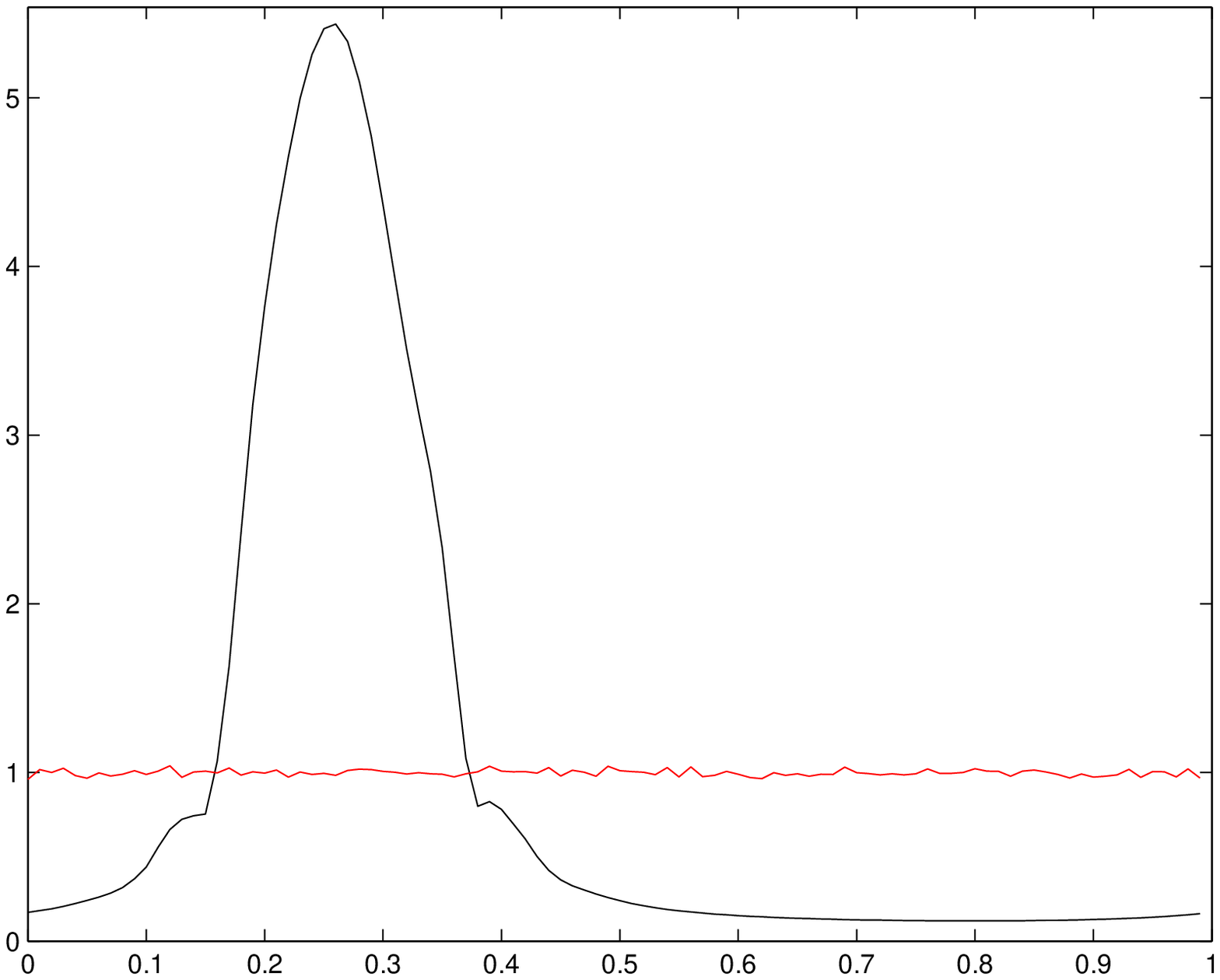}}
\caption{Second order model with limited vision and $R = 0.07$; the red line in the left column indicates the initial datum.}\label{fig:num_secorder1}
\end{center}
\end{figure}

We conclude this section with simulations of the second-order mean-field model \fref{mmodel2} in 1D.
The spatial domain $\Omega = [a,b]$ is discretized using an equidistant mesh $x_i = a + i \Delta x$,
the velocity domain $V = [v_{\min},v_{\max}]$ at grid points $v_j = v_{\min} + j \Delta v$.
The time step $\Delta t$ satisfies the CFL condition for $v_{\max}$, i.e.,
$\Delta t = \frac{\Delta x}{\lvert v_{\max}\rvert}$.
The numerical scheme is based on a splitting method: Given an initial datum $f(x,v,0) = f_0(x,v)$,
we split the system at every time $t_k = k\Delta t$ into the following steps:
\begin{enumerate}
\item Solve transport equation in $x$
\begin{align}\label{e:trans}
\frac{\partial f}{\partial t} + v_i \cdot \nabla_x f &= 0,
\end{align}
subject to periodic boundary conditions on $\Omega$,
for every $v_i$ on the time interval $t \in [t_k, t_k + \frac{1}{2} \Delta t]$,
using an upwind scheme with the superbee flux limiter.
\item Starting with the solution of the transport equation \fref{e:trans}, solve
\begin{align}\label{e:convdiff}
\frac{\partial f}{\partial t} = \nabla_v \cdot \left(H (W \circledast f) v f + \frac{1}{2} \nabla_v(G(W\circledast f)^2f )\right),
\end{align}
with no flux boundary conditions on $V$, on the time interval $[t_k, t_{k+1}]$,
using a semi-implicit time discretization.
\item Finally solve \fref{e:trans} using the solution of \fref{e:convdiff} for another half time step $\frac{1}{2}\Delta t$.
\end{enumerate}
We set the computational domain to $\Omega = [0,1]$, the velocity domain to $V = [-1,1]$
and the mesh sizes to $\Delta x = 10^{-2}$ and $\Delta v = 2 \times 10^{-2}$.
We choose similar conditions as in the individual based model, i.e. a limited cone of vision
$W(x,v) = w(\lvert x \rvert, \frac{ x \cdot v}{\lvert x \rvert \lvert v\rvert})$
with $w(s,z) = \chi_{[0.R]}(s)\chi_{[0,1]}(z)$.
The sampling radius is set to $R = 0.07$, $G(s) = \exp(-2 s)$ and $H \equiv 2$.
The initial datum corresponds to a small perturbation of a uniform distribution with mass one.
Snapshots of the evolution of the particle distribution density $f(t,x,v)$
and the mass density $\rho(t,x)$ at different times are depicted in Fig. \ref{fig:num_secorder1}.
We observe a fast smoothing of $f(t,x,v)$ in time and a subsequent formation of a stable aggregate.

\begin{figure}
\begin{center}
\subfigure[Particle distribution function $f(t,x,v)$ at $t=4$]{\includegraphics[width = 0.45 \textwidth]{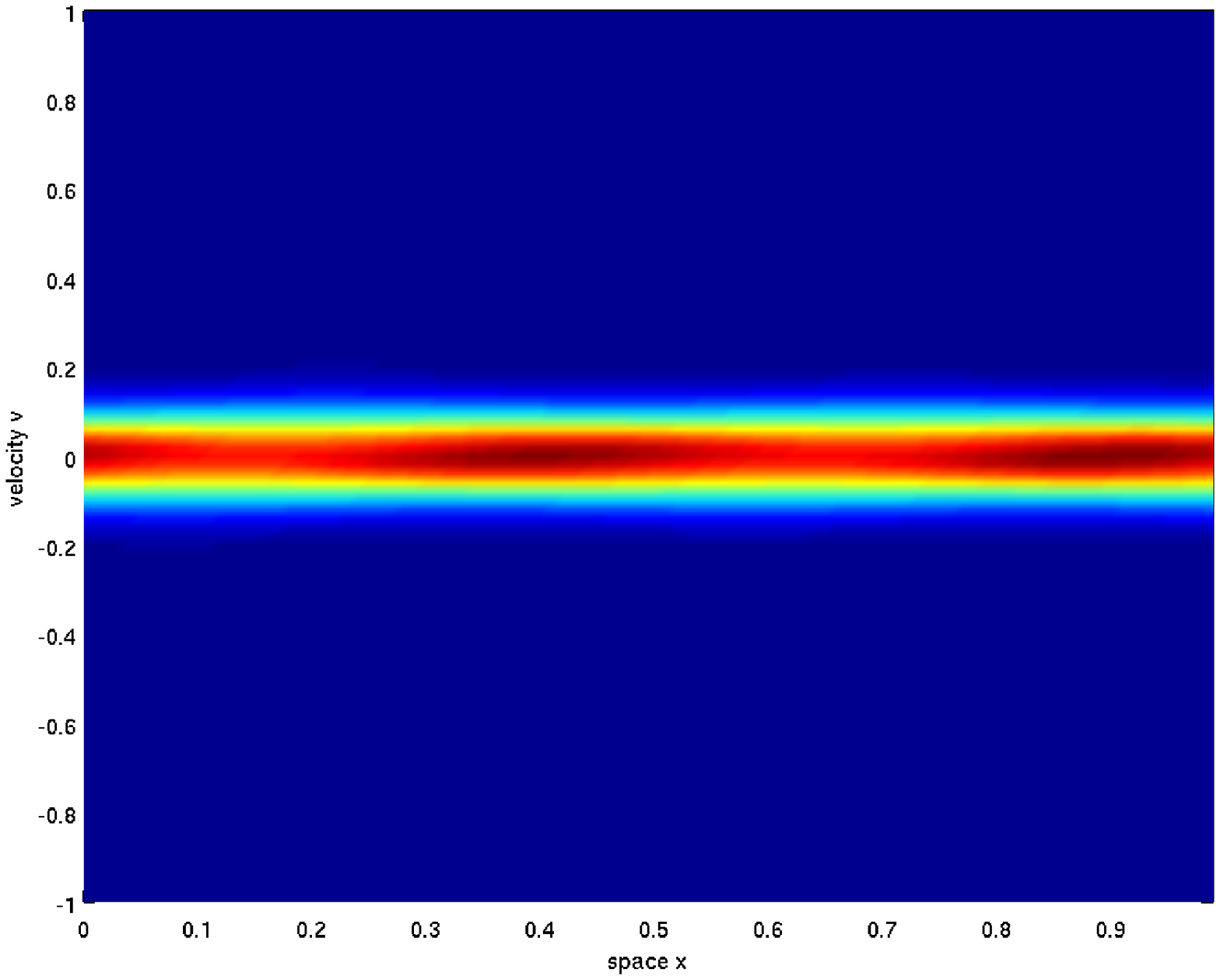}}\hspace*{1cm}
\subfigure[Mass density $\rho(t,x)$ at $t=4$]{\includegraphics[width = 0.45 \textwidth]{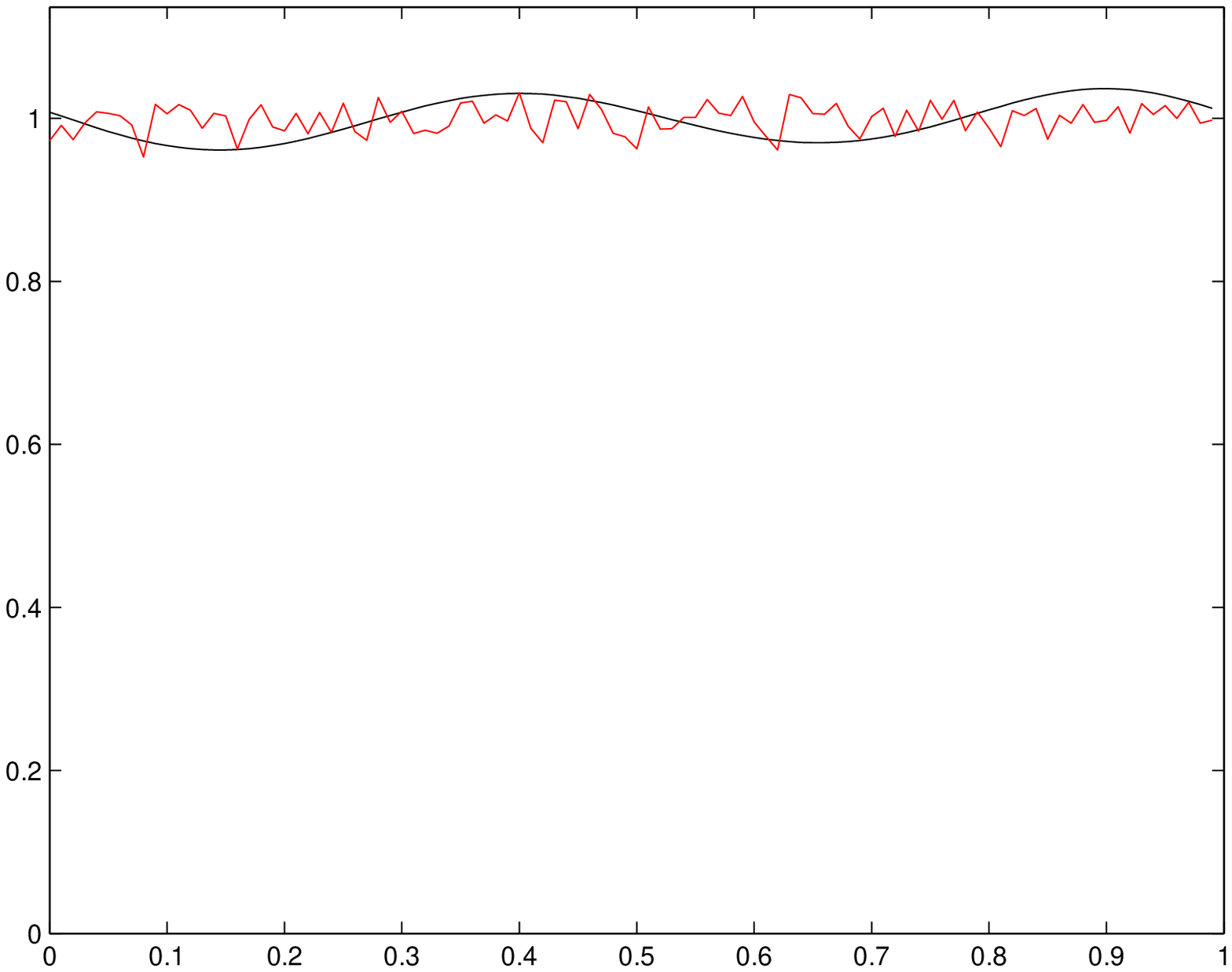}}\\
\subfigure[Particle distribution function $f(t,x,v)$ at $t=8$]{\includegraphics[width = 0.45 \textwidth]{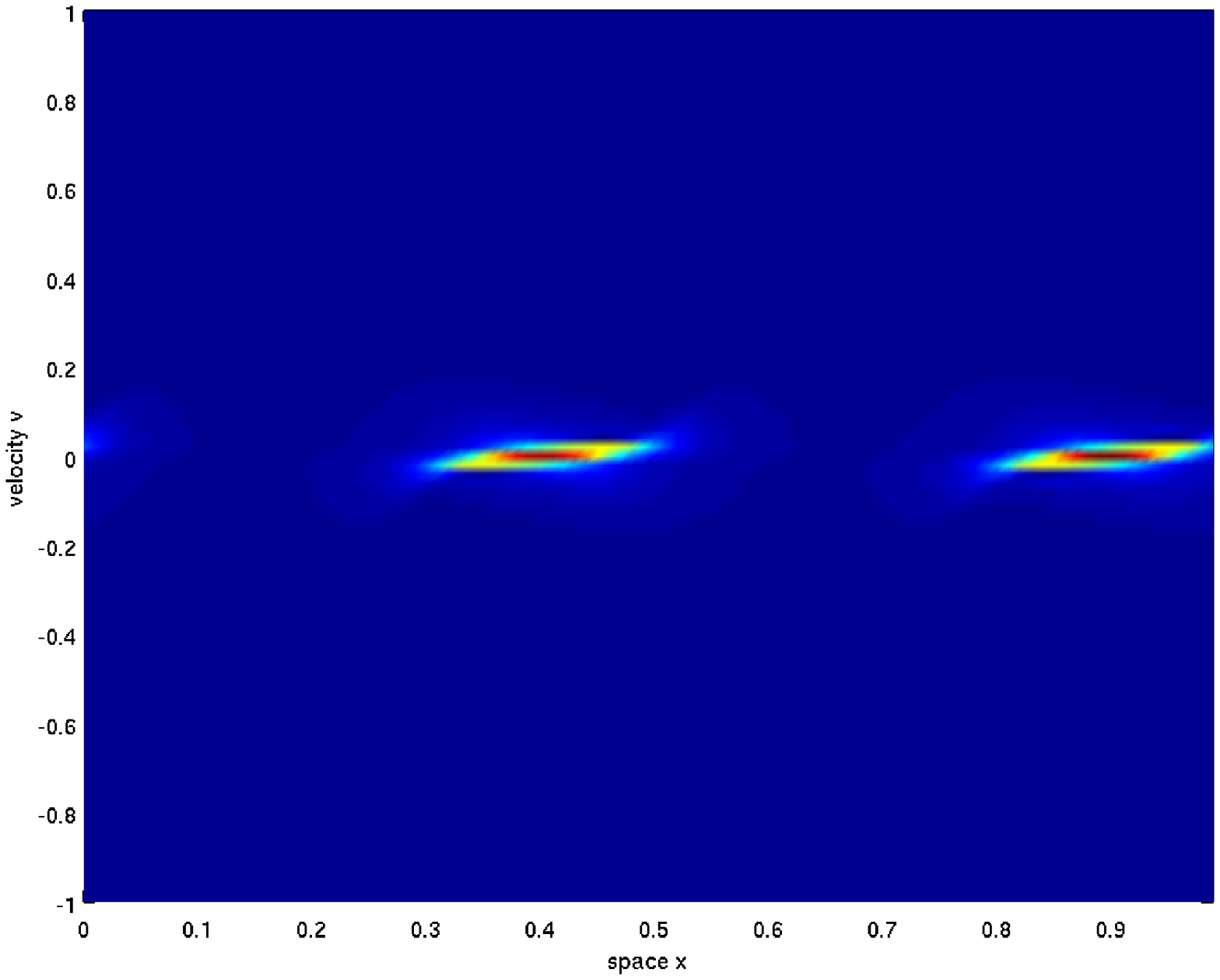}}\hspace*{1cm}
\subfigure[Mass density $\rho(t,x)$ at $t=8$]{\includegraphics[width = 0.45 \textwidth]{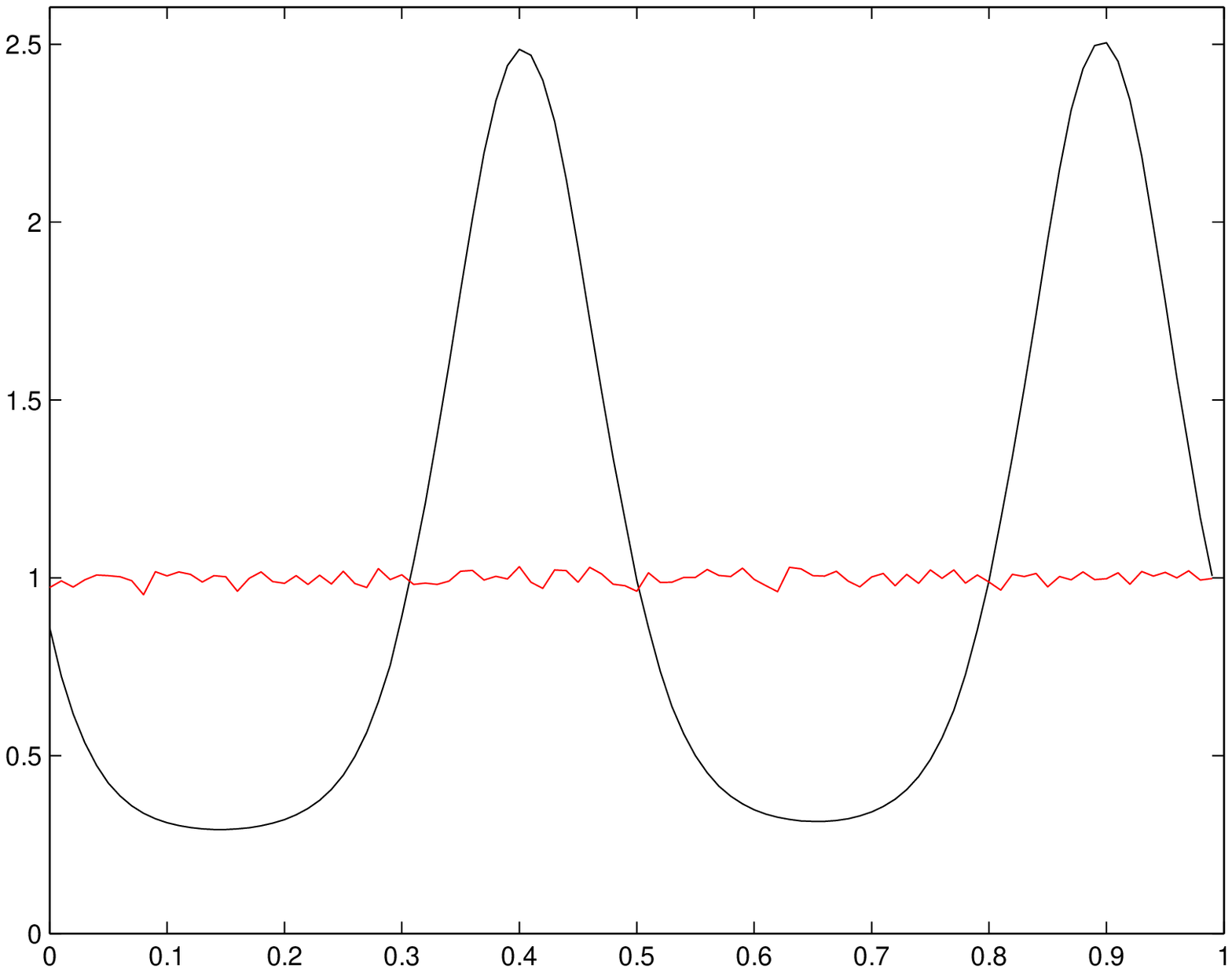}}\\
\subfigure[Particle distribution function $f(t,x,v)$ at $t=20$]{\includegraphics[width = 0.45 \textwidth]{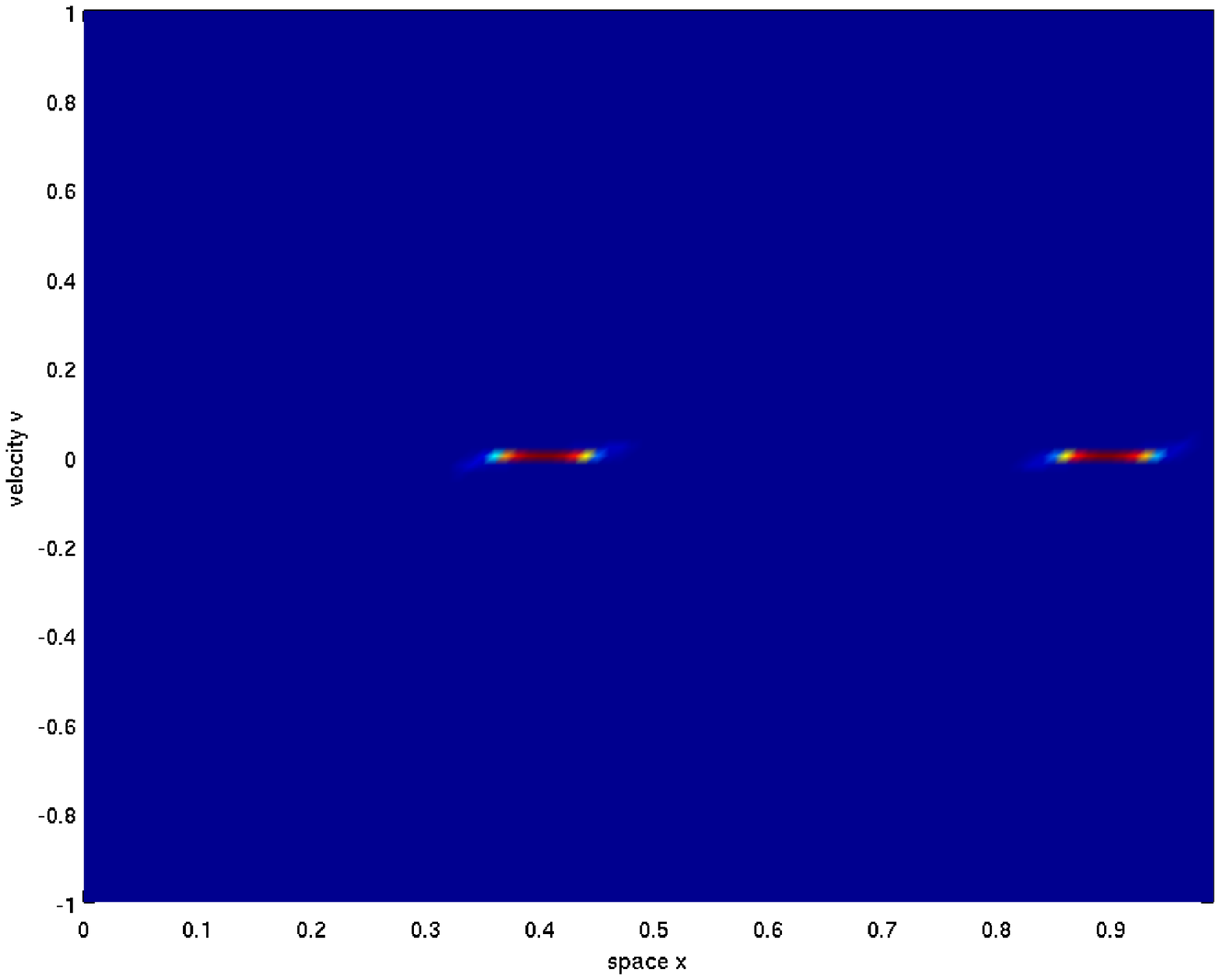}}\hspace*{1cm}
\subfigure[Mass density $\rho(t,x)$ at $t=20$]{\includegraphics[width = 0.45 \textwidth]{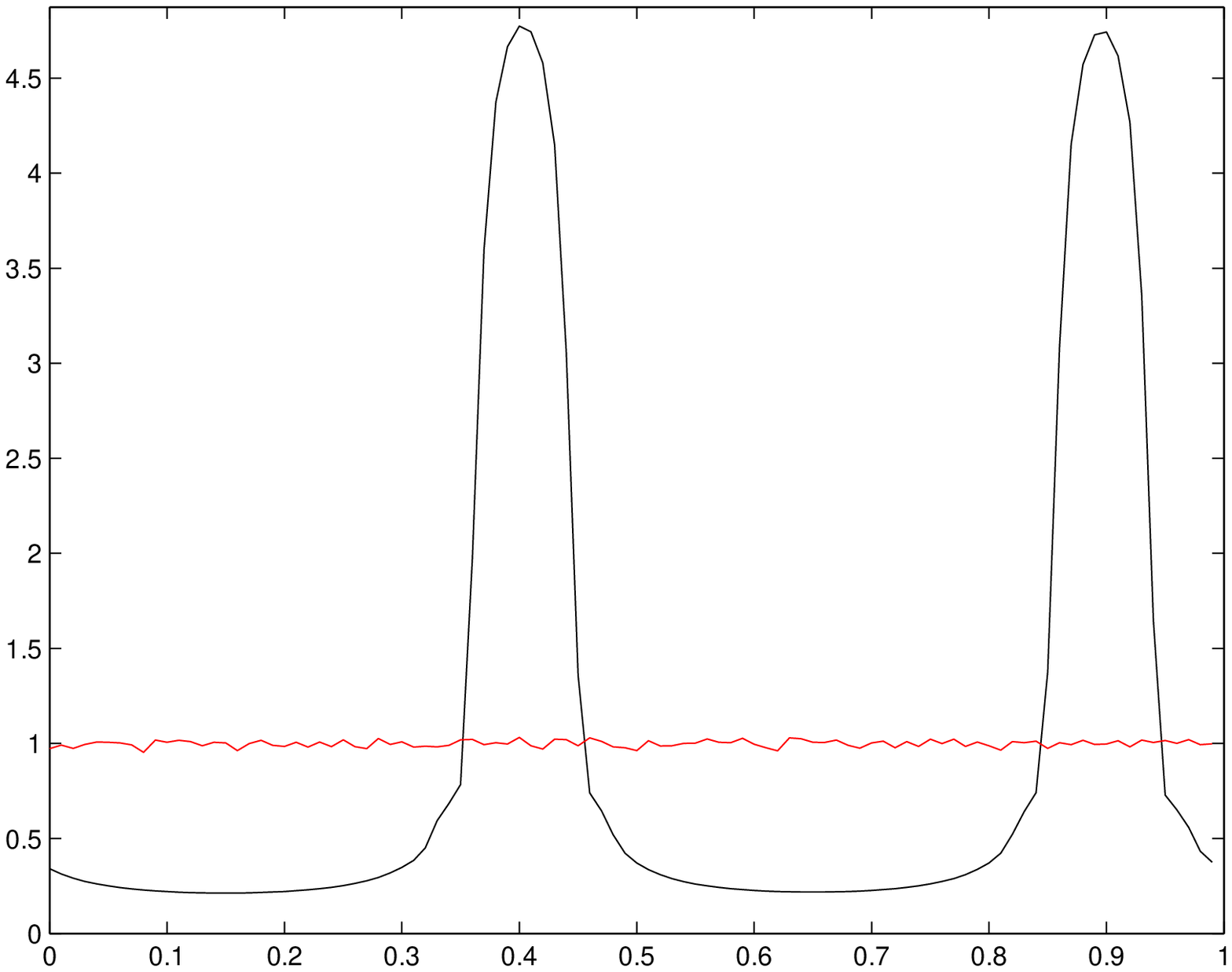}}
\caption{Second order model with limited vision and $R = 0.03$; the red line in the left column indicates the initial datum.}\label{fig:num_secorder2}
\end{center}
\end{figure}

If we decrease the sampling radius to $R = 0.03$, two separate aggregates form,
see Fig. \ref{fig:num_secorder2}.
Again, we observe that with a smaller sampling radius the aggregation
happens on a faster time scale.

\newpage

\noindent{\bf Acknowledgment:}
JH acknowledges 
the hospitality of the Faculty of Mathematics, University of M\"unster
during his stay, where the work leading to this paper has been initiated.
MTW acknowledges financial support from the Austrian Science Foundation (FWF) via
the Hertha Firnberg project T456-N23 and the Award KUK-I1-006-43 by the 
King Abdullah University of Science and Technology (KAUST).


\end{document}